# Higher-order diffusion and Cahn–Hilliard-type models revisited on the half-line

A. Chatziafratis [1], A. Miranville [2], G. Karali [3], A.S. Fokas [4], E.C. Aifantis [5]

[1] Department of Mathematics, National and Kapodistrian University of Athens, Greece
Department of Mathematics and Statistics, School of Pure and Applied Sciences, University of Cyprus
Institute of Applied and Computational Mathematics, FORTH, Crete, Greece

[2] Henan Normal University, School of Mathematics and Statistics, Xinxiang, China
Laboratoire de Mathématiques Appliquées du Havre (LMAH), Université Le Havre Normandie, France

[3] Department of Mathematics, National and Kapodistrian University of Athens, Greece

[4] Department of Applied Mathematics and Theoretical Physics, University of Cambridge, UK
Viterbi School of Engineering, University of Southern California, Los Angeles, California, USA
Mathematics Research Center, Academy of Athens, Greece

[5] Friedrich-Alexander University of Erlangen-Nuremberg, Germany
Laboratory of Mechanics and Materials, College of Engineering, Aristotle University, Greece
College of Engineering, Michigan Technological University, USA

**Abstract.** In this paper, we solve explicitly and analyze rigorously inhomogeneous initial-boundary-value problems (IBVP) for several fourth-order variations of the traditional diffusion equation and the associated linearized Cahn-Hilliard (C-H) model (also Kuramoto-Sivashinsky equation), formulated in the spatiotemporal quarter-plane. Such models are of relevance to heat-mass transfer phenomena, solid-fluid dynamics and the applied sciences. In particular, we derive formally effective solution representations, justifying *a posteriori* their validity. This includes the reconstruction of the prescribed initial and boundary data, which requires careful analysis of the various integral terms appearing in the formulae, proving that they converge in a strictly defined sense. In each IBVP, the novel formula is utilized to rigorously deduce the solution's regularity and asymptotic properties near the boundaries of the domain, including uniform convergence, eventual (long-time) periodicity under (eventually) periodic boundary conditions, and null non-controllability. Importantly, this analysis is indispensable for exploring the (non)uniqueness of the problem's solution and a new counter-example is constructed. Our work is based on the synergy between: (i) the well-known Fokas unified transform method and (ii) a new approach recently introduced for the rigorous analysis of the Fokas method and for investigating qualitative properties of linear evolution partial differential equations (PDE) on semi-infinite strips. Since only up to third-order evolution PDE have been investigated within this novel framework to date, we present our analysis and results in an illustrative manner and in order of progressively greater complexity, for the convenience of readers. The solution formulae established herein are expected to find utility in well-posedness studies for nonlinear counterparts too.

## 1. Introduction

*Overview of Cahn–Hilliard-type model equations*

The standard Cahn–Hilliard equation is a model for the phase separation of a binary alloy proposed in [1]. This model was extended later by [2,3] in order to incorporate thermal fluctuations in the form of an additive noise. The initial value Cahn-Hilliard problem was treated in [4]. The equilibrium theory for the Cahn-Hilliard equation has been studied in one dimension in [5] on the finite interval and in [6] on the real line. It was suggested by heuristic arguments and computational experiments that a near homogeneous state in the spinodal interval rapidly evolves to a highly oscillatory state. The layers have width of O($\varepsilon$) and the oscillations are between small neighborhoods of +1 and -1. Some rigorous results on this direction for the "spinodal decomposition" phenomenon can be found in [7]. Concerning solutions of the Cahn-Hilliard in one-dimension and their slow-motion, results can be found in [8] for the two-layers case, and in [9] for more layers. Therein, it was significant to construct an approximate manifold and provide appropriate spectral estimates. In higher space dimensions, the analog of the layers are curves in two dimensions or surfaces that are either closed or intersect the boundary. Numerous works have been done in these directions. We refer indicatively to the work of Alikakos and Fusco [10] for the two-dimensional case and to [11] for higher space dimensions. In higher dimensions super-slow motion was first noticed and established in [10], where it was shown that single bubble-like solutions persist and evolve exponentially slowly towards the closest point on the boundary. In [12], it was proved that indeed this was true provided the bubble is sufficiently small. There exist several results for the sharp interface limit of the deterministic Cahn-Hilliard equation resulting to coarsening phenomena of alloys; we refer for example to [12] for the two- and three-dimensional cases. Concerning the stochastic Cahn-Hilliard equation, some interesting results may be found in [13,14]. Further results and a systematic survey of recent progress in the mathematical analysis of the C-H model and variants can be found in Miranville's book [15].

More specifically, regarding 1+1 dimensions, the Cahn-Hilliard 4th-order linearized equation, i.e. $u_t = \alpha u_{xx} - \beta u_{xxxx}$ (with $\alpha < 0$, $\beta > 0$), was proposed to study the initial stages of spinodal decomposition [1]: diffusion-controlled phase transformation initially observed in homogeneous solid solution metal alloys (later in other systems as well including polymers and complex fluids), when the temperature is suddenly decreased to the so-called spinodal region. For example, when the temperature of an initially homogeneous solid solution Ni-Cu alloy is rapidly decreased to the (far from equilibrium) spinodal regime, the homogeneously distributed atomic species move from regions of low to high concentrations (as opposed to their usual motion from high to low concentrations near equilibrium), and this "uphill" diffusion process leads eventually to distinct $Ni_x$-$Cu_y$ stable phases. This may be viewed as a "negative" physics (negative diffusion coefficient) phenomenon – similar to others (negative compressibility in van der Waals equation of state, negative elasticity or Poisson's ratio in Hooke's law, negative thermal or hydraulic conductivity in Fourier or Darcy laws) – observed at the macro-scale for material states controlled by micro- or nano- scale processes.

To better understand the physical implications of the linearized C-H equation, we consider small compositional perturbations of the form proportional to $\exp(iqx + \omega t)$ from a uniform average composition $u_0$. A simple calculation shows that these modulations decay for $\omega = -q^2(\alpha + \beta q^2) < 0$ and grow when $\omega > 0$. For the case $\beta = 0$, it follows that they grow exponentially with time when $\alpha < 0$, i.e., for negative diffusion coefficient, suggesting the ill-posedness of the backwards heat equation. By allowing for higher-order gradient effects ($\beta > 0$), well-posedness is restored above a critical wavenumber $q > q_c$ with $q_c = \text{sqrt}(\alpha/\beta)$. Analogous types of pde gradient "regularizations" have been used for elasticity (negative Poisson's ratio) and plasticity (negative hardening). For a recent account of gradient theories, the reader may consult the first chapter of [16] and the final chapter of [17] where a wide spectrum of applications in various scales and disciplines with emphasis on instabilities and pattern formation as well as removal of singularities and size effects are given.

The first to introduce higher-order spatial derivatives in the free energy function was van der Waals [18] to treat fluid states with negative compressibility and obtain "smooth" instead of sharp density profiles of liquid-vapor interfaces. Van der Waals thermodynamic theory was revisited by Aifantis and Serrin [19] within a purely mechanical framework leading to a more general differential equation for the interfacial density. Closed-form solutions of this equation were obtained: (i) transitions i.e. smooth density profiles between the liquid and vapor phases; (ii) reversals i.e. liquid (or vapor) bubbles within a vapor (or liquid) phase; and (iii) periodic ones i.e. periodic liquid-gas layers of equal thickness. Moreover, Maxwell's "equal area rule" line – defining the equilibrium pressure and the corresponding vapor and liquid densities – was relocated as a condition for the existence of such solutions. It was shown, in addition, that the Aifantis-Serrin purely mechanical theory coincides with van der Waals thermodynamic theory only for a special form of the gradient constitutive relation for the interfacial stress. While this special form of the interfacial stress leads back to Maxwell's equal rule, it reveals an unrepairable situation for the Davis-Scriven [20] statistical mechanics formulation of van der Waals theory and their gradient expression for the corresponding chemical potential.

The aforementioned Aifantis-Serrin mechanical theory was motivated by Aifantis' earlier approach to diffusion [21] where the assumption of existence of a thermodynamic "chemical potential" was abandoned in favor of the introduction of a mechanical "diffusive force" arising from the exchange of momentum between the diffusion species and the surrounding matter. This quantity – first envisioned in Maxwell's kinetic theory of gases [22] – enters as an "internal body force" in the corresponding equation of momentum balance which, in conjunction with the mass balance equation for the diffusing species provide a mechanical basis for diffusion more general (economy of further thermodynamic axioms) framework for deriving mass transfer equations. In fact, by following the recipe of continuum mechanics and adopting corresponding constitutive equations for the stress and the diffusive force entering the differential statement of momentum conservation, in conjunction with the differential statement of mass conservation, various classes of PDE were derived [21] to model mass transport in heterogeneous media, including the classical Fick's diffusion equation, Barenblatt's pseudo-parabolic equation of seepage, and C-H spinodal decomposition equation, the mathematics of which will be further analyzed in the sequel.

Before this, however, we provide for completeness its brief derivation within the above mentioned mechanical framework and one-dimensional setting by also partly adopting for convenience the notation used in the original sources. In 1D the differential forms of the mass and momentum balance equations for the diffusing species read $\rho_t + j_x = 0$ and $T_x = \underline{f}$ where $\rho$ and $j = \rho v$ (v stands for the velocity) denote respectively concentration and flux of the diffusing substance, while $T$ denotes the (partial) stress the diffusing species exert on themselves and $\underline{f}$ their exchange of momentum with the surrounding medium. With the simplest constitutive assumptions $T = -\pi \rho$, $\underline{f} = \kappa$ (where $\pi$ and $\kappa$ being constants), the above stated momentum balance leads to Fick's 1st law of diffusion $j = -D\rho_x$ ($D=\pi/\kappa$) which upon use of the mass balance leads to the 2nd Fick's law $\rho_t=D\rho_{xx}$ i.e. the classical parabolic diffusion equation. If instead of the above perfect fluid (ideal gas) assumption the diffusing substance is assumed to behave like a viscous fluid then $T=-\pi\rho+\underline{\pi}\rho_x$ (with $\underline{\pi}$ being a viscosity-like constant) and the corresponding diffusion equation reads $\rho_t=D\rho_{xx}+\underline{D}\rho_{txx}$ ($\underline{D}=\underline{\pi}/\kappa$) i.e. Barenblatt's equation. Finally, if a weak non-locality is assumed for T of the form $T=-\pi\rho+\varepsilon\rho_{xx}$ (with $\varepsilon$ being constant) then the corresponding diffusion equation reads $\rho_t=D\rho_{xx}-E\rho_{xxxx}$ ($E=\varepsilon/\kappa$) i.e. the C-H linearized spinodal decomposition equation. More details on 3D derivation of the above equations and other mass transfer models can be found in [21] and related mathematical analyses in [23].

We hope that the above overview effectively illustrates the importance of these models justifying the necessity for investigating further their mathematical properties. In what follows, we obtain, via rigorous implementation of the Fokas method, and analyze solution formulae for fully non-homogeneous IBVP for various C-H models (also related to the linearized version of the well-known Kuramoto-Sivashinsky equation [24,25] which arises in contexts of mathematical fluid and ion physics).

It is important to note that, while the mathematical analysis of the (linearized) Cahn-Hilliard equation is well established in a *bounded* domain, the situation is less clear in an semi-unbounded domain, even in one space dimension. This constitutes another motivation for our study. We simply mention here the work of Caffarelli and Müller [26] for the Cahn-Hilliard equation in the whole space $R^n$. Again, we refer to [15] for a recent comprehensive review of advances pertaining to C-H equations.

*A modern method approach to initial–boundary-value problems*

In order to place the Fokas method [27] into historical perspective, we note that there are two well-known approaches to the exact analysis of linear PDEs, the method of separation of variables and the approach based on Green's integral identities. The use of separation of variables gives rise to ordinary differential operators, and then, the spectral analysis of the associated operators yields an appropriate transform pair. The prototypical such pair is the Fourier transform; variations include the sine, the cosine, the Laplace, and the Mellin transforms. However, for non-self-adjoint problems such transforms generally do not exist. Regarding the use of the integral representations obtained via Green's functions, we note that these representations involve certain unknown boundary values which can be eliminated only in very special cases (using the method of images).

In the second half of the 20th century it was realized that certain nonlinear evolution PDEs, called integrable, can be formulated as the compatibility condition of two linear eigenvalue equations called a Lax pair. This formulation gives rise to a method for solving the initial value problem for these equations, called the inverse scattering transform method. Actually, this method is based on a deeper form of separation of variables. Motivated by his work in the extension of the inverse scattering transform for solving boundary value problems, one of the authors introduced a completely new method for solving linear PDEs. In the last couple of decades, it has been established, with the aid of many mathematics researchers in hundreds of articles, that this method provides the most efficient approach for explicitly solving initial-boundary-value problems for linear PDEs.

For linear evolution PDEs with polynomial dispersion relation, the unified transform method (UTM) involves three steps: (1) Derive formally an integral representation for the solution, in the complex Fourier plane. This representation, just like the analogous one obtained via Green's integral identities, contains certain transforms of unknown boundary values; thus it is not yet effective. (2) Analyze the global relation, which is an equation coupling the given initial and boundary data with the unknown boundary values. (3) Use certain invariant properties of the global equation and simple algebraic manipulations to eliminate, from the integral representation obtained in step one, the transforms of the unknown boundary values thereby constructing a candidate solution to the given IBVP.

The Fokas method has been spectacularly successful, both analytically and numerically, in a wide class of linear-PDE problems – of any order and in a variety of domains. This includes elliptic PDEs, PDEs with variable coefficients, and even problems involving moving boundaries. It has also been implemented to cases of nonlocal (nonseparable) boundary conditions, interface problems and fractional-order evolution equations too. Incidentally, this method has led to surprising results in spectral theory, inverse problems of medical imaging, as well as to the development of an effective approach to rigorously establishing well-posedness results for nonlinear evolution PDE in low-regularity settings. For all the above, see, e.g., [27-71] (in chronological order) and many references therein.

Recently, a new approach to rigorous aspects of the Fokas method in classical settings has been introduced in [72] and extended in e.g. [73-77]. This new approach has led to valuable new knowledge for a variety of classical IBVPs and PDEs. In addition, it has provided a desirable refinement (extension) and motivated new areas of applicability of the Fokas method, including the analysis of 'distant' sensitivity (or long-range instabilities) of solutions of dispersive PDEs, well-posedness and qualitative theory for linear PDEs on the quarter-plane, etc. More studies along these lines will appear in forthcoming papers, e.g., [78-84].

*Outline of the current paper*

The present work stems from the latter line of inquiry mentioned just above. In more detail, firstly, effective solution representations are formally derived for the spacetime quarter-plane, starting from the one-parameter divergence forms associated with the PDEs under consideration, and then, importantly, detouring into the complex spectral plane. Then, as necessary, the validity of the proposed solution formulae are a posteriori justified analytically, including the reconstruction of prescribed initial and boundary conditions, not only in the evaluation but –crucially– in the limit sense (as the independent variables approach the boundaries) as well. The important and intensive verification stage essentially requires careful interpretation of various integral terms appearing in the proposed solution representations thereby securing that they converge in a strictly defined sense. The novel solution formulae are utilized in order to thoroughly investigate and deduce the solution's regularity properties near the boundaries (the positive semi-axes) of the spatiotemporal domain. Moreover, novel solution-uniqueness results (under certain hypotheses) and the problems' well-posedness are established too, crucially based on the preceding boundary behavior analysis; the study, via the new formulas, of well-posedness for nonlinear analogues with low-regularity data is work in progress. A counter-example –previously unknown for any higher-order linear evolution PDE of the types considered presently– which demonstrates non-uniqueness is also provided.

In addition, we perform a careful and thorough investigation of eventual periodicity properties, at large times, of the solutions. Namely, the important question (which has been addressed elsewhere for linear models of only lower order [41,62,65,85-88]) of whether the solution asymptotically becomes periodic, at any fixed position, under the effect of (eventually) t-periodic boundary data is, in our setting, positively settled. Furthermore, a theorem pertaining to null controllability is presented, again seminal in the setting of 4th-order PDEs. Notably, the latter, in combination with the new analytical approach and [89], has inspired a general framework for generating rigorous controllability results for a broad class of evolution and dispersive PDEs; this will be announced in a forthcoming publication [81].

In conclusion, linear PDEs have traditionally been of great importance, as abstract objects of theoretical interest, as practical models of physical and engineering processes, and as intermediate stages towards investigation of nonlinear problems (in view of e.g. linearizability of certain classes of PDE [90-92], linear estimates necessary for addressing well-posedness questions in low-regularity settings, e.g., with data in Sobolev spaces, and so on). Evidently, the present essay is the first treatise on the rigorous analysis of fully inhomogeneous IBVPs for linear 4th-order evolution PDEs on semi-unbounded domains, even though purely non-linear counterparts, the Cauchy problem on the whole space and IBVPs on bounded domains (for both linear and nonlinear cases) have already been studied in the pertinent literature. Our paper, in addition to expanding the aforementioned ongoing rigorously-minded research programme and achieving new knowledge for celebrated mathematical models of the applied sciences via a modern method, it also opens up a road to future explorations.

*Technical preliminaries and notation*

In what follows, we study analytically the equation $\partial_t u = \alpha \partial_x^2 u - \beta \partial_x^4 u + f$, on the half-line, with $\alpha \in \mathbb{R}$ and $\beta > 0$. More precisely, given $u_0(x)$, $g_0(t)$, $g_1(t)$ and $f(x,t)$, we rigorously solve, under certain assumptions on these data functions, the initial boundary value problem (IBVP)

$$\begin{cases} \dfrac{\partial U}{\partial t} = \alpha \dfrac{\partial^2 U}{\partial x^2} - \beta \dfrac{\partial^4 U}{\partial x^4} + f, \ (x,t) \in Q := \mathbb{R}^+ \times \mathbb{R}^+, \\[2mm] \lim_{t \to 0^+} U(x,t) = u_0(x), \ x \in \mathbb{R}^+, \\[2mm] \lim_{x \to 0^+} U(x,t) = g_0(t), \ t \in \mathbb{R}^+, \\[2mm] \lim_{x \to 0^+} U_x(x,t) = g_1(t), \ t \in \mathbb{R}^+, \end{cases} \qquad (1.1)$$

for $U = U(x,t)$.

Throughtout this paper, we make the following assumptions on the data:

$$u_0(x) \in \mathcal{S}\big([0,\infty)\big), \ g_0(t), g_1(t) \in C^\infty\big([0,\infty)\big) \ \text{ and } \ f = f(x,t) \in C^\infty(\overline{Q}) \ \text{ such that } \ f(\cdot,t) \in \mathcal{S}\big([0,\infty)\big). \qquad (1.2)$$

More precisely, the last assumption on the function $f(x,t)$ means that it is rapidly decreasing with respect to $x$, uniformly for $t$ in compact subsets of $[0,+\infty)$, i.e., for every $\ell, n \in \mathbb{N} \cup \{0\}$ and $t_0 > 0$,

$$\sup \left\{ x^{\ell} \left| \frac{\partial^n f(x,t)}{\partial x^n} \right| : x \geq 0, \, 0 \leq t \leq t_0 \right\} < +\infty \, . \tag{1.3}$$

The technique we employ as a starting point for our purposes is the Fokas' unified transform method.

In order to make the construction and the verification of the solution more transparent, we consider first the cases $[\,\alpha = 0, \, \beta = 1\,]$ and $[\,\alpha = 1, \, \beta = 1\,]$, i.e., the equations

$$\partial_t U = -\partial_x^4 U + f \quad \text{and} \quad \partial_t U = \partial_x^2 U - \partial_x^4 U + f \, .$$

***Notation*** For a function $u_0(x) \in \mathcal{S}\big([0,\infty)\big)$, we define its Fourier transform

$$\hat{u}_0(\lambda) = \int_{y=0}^{\infty} e^{-i\lambda y} u_0(y) dy, \text{ for } \lambda \in \mathbb{C} \text{ with } \operatorname{Im}\lambda \leq 0 \, .$$

Setting

$$\sigma_{\mathrm{M}}(\lambda) := \sum_{j=1}^{\mathrm{M}} \frac{d^{j-1}u_0}{dx^{j-1}}(0) \frac{1}{(i\lambda)^j} \quad (\lambda \in \mathbb{C}, \, \lambda \neq 0)$$

and integrating by parts, we obtain

$$\hat{u}_0(\lambda) - \sigma_{\mathrm{M}}(\lambda) = \frac{1}{(i\lambda)^{\mathrm{M}}} \int_{y=0}^{\infty} e^{-i\lambda y} \frac{d^{\mathrm{M}}u_0(y)}{dy^{\mathrm{M}}} dy = \mathrm{O}(1/\lambda^{\mathrm{M}+1}), \text{ as } \lambda \to \infty \text{ with } \operatorname{Im}\lambda \leq 0 \, . \tag{1.4}$$

In particular,

$$\hat{u}_0(\lambda) = \frac{u_0(0)}{i\lambda} + \frac{(u_0{}')\hat{}(\lambda)}{i\lambda} = \mathrm{O}(1/\lambda) \text{ and } \frac{(u_0{}')\hat{}(\lambda)}{i\lambda} = \mathrm{O}(1/\lambda^2), \text{ as } \lambda \to \infty \text{ with } \lambda \in \mathbb{C}, \, \operatorname{Im}\lambda \leq 0 \, . \tag{1.5}$$

## 2. The equation $\partial_t u = -\partial_x^4 u$

***Problem*** Given $u_0(x) \in \mathcal{S}\big([0,\infty)\big)$ and $g_0(t), g_1(t) \in C^{\infty}\big([0,\infty)\big)$, solve

$$\begin{cases} \partial_t u = -\partial_x^4 u, \, (x,t) \in Q := \mathbb{R}^+ \times \mathbb{R}^+ \\[4pt] \lim_{t \to 0^+} u(x,t) = u_0(x), \, x \in \mathbb{R}^+ \\[4pt] \lim_{x \to 0^+} u(x,t) = g_0(t), \, t \in \mathbb{R}^+, \\[4pt] \lim_{x \to 0^+} u_x(x,t) = g_1(t), \, t \in \mathbb{R}^+, \end{cases} \tag{2.1}$$

for $u = u(x,t)$.

## 2.1 Formal derivation of the solution via the Fokas method

Setting $\omega(\lambda) := \lambda^4$, so that $\omega(\lambda)\big|_{\lambda = -i\partial_x} = \partial_x^4 u$, and in view of the identity

$$i\left( \frac{\lambda^4 - \vartheta^4}{\lambda - \vartheta} \right)\Bigg|_{\vartheta = -i\partial_x} u = i[\vartheta^3 + \vartheta^2\lambda + \vartheta\lambda^2 + \lambda^3]\Big|_{\vartheta = -i\partial_x} u = -u_{xxx} - i\lambda u_{xx} + \lambda^2 u_x + i\lambda^3 u \, ,$$

we are led to the equation

$$\frac{\partial}{\partial t}[e^{-i\lambda x + \omega(\lambda)t} u(x,t)] - \frac{\partial}{\partial x}[e^{-i\lambda x + \omega(\lambda)t}(-u_{xxx} - i\lambda u_{xx} + \lambda^2 u_x + i\lambda^3 u)] = e^{-i\lambda x + \omega(\lambda)t}\left( \frac{\partial u}{\partial t} + \frac{\partial^4 u}{\partial x^4} \right).$$

Fixing $t > 0$ and applying Green's formula in the strip $\Pi_t := \{(x,\tau) : x \geq 0, 0 \leq \tau \leq t\}$ (see fig.1),

A. Chatziafratis, A. Miranville, G. Karali, A.S. Fokas, E.C. Aifantis

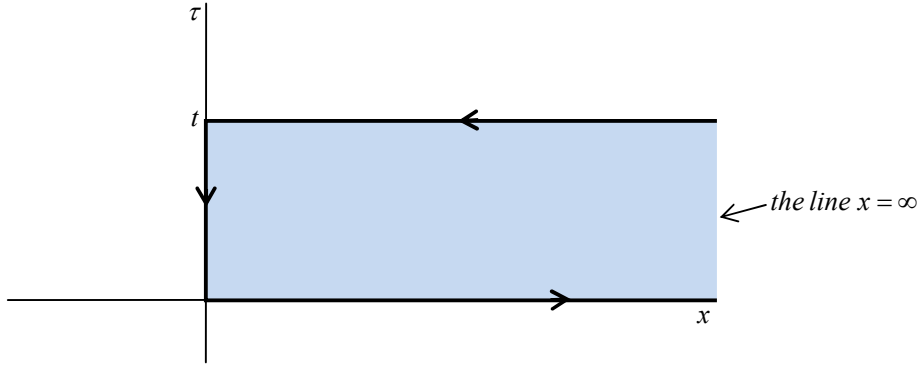

**Fig.1** The semi-strip $\Pi_t$

we obtain

$$\iint_{\partial \Pi_t} [e^{-i\lambda x + \omega(\lambda)\tau} u(x,\tau)]dx + [e^{-i\lambda x + \omega(\lambda)\tau}(-u_{xxx} - i\lambda u_{xx} + \lambda^2 u_x + i\lambda^3 u)]d\tau = -\iint_{\Pi_t} e^{-i\lambda x + \omega(\lambda)t}\left(\frac{\partial u}{\partial \tau} + \frac{\partial^4 u}{\partial x^4}\right)dxd\tau \ . \quad (2.2)$$

Thus, if the function $u = u(x,t)$ satisfies the equation $\partial_t u = -\partial_{xxxx} u$ then

$$\int_{x=0}^{\infty} e^{-i\lambda x} u(x,0)dx - \int_{x=0}^{\infty} e^{-i\lambda x + \omega(\lambda)t} u(x,t)dx - \int_{\tau=0}^{t} e^{\omega(\lambda)\tau}[-u_{xxx}(0,\tau) - i\lambda u_{xx}(0,\tau) + \lambda^2 u_x(0,\tau) + i\lambda^3 u(0,\tau)]d\tau = 0 \ ,$$

where we interpreted the line integral in the LHS in (2.2).
It follows that

$$\hat{u}_0(\lambda) - \hat{u}(\lambda,t)e^{\omega(\lambda)t} + \widetilde{g}_3(\omega(\lambda),t) + i\lambda\widetilde{g}_2(\omega(\lambda),t) - \lambda^2\widetilde{g}_1(\omega(\lambda),t) - i\lambda^3\widetilde{g}_0(\omega(\lambda),t) = 0 \ , \text{ for } \lambda \in \mathbb{C} \text{ with } \operatorname{Im}\lambda \le 0 \ , (2.3)$$

where

$$g_3(t) := u_{xxx}(0,t) \ , \ g_2(t) := u_{xx}(0,t) \ , \ g_1(t) := u_x(0,t) \ , \ g_0(t) := u(0,t) \ ,$$

and

$$\widetilde{g}_j(\omega(\lambda),t) := \int_{\tau=0}^{t} e^{\omega(\lambda)\tau} g_j(\tau)d\tau \ , \ j = 0,1,2,3, \ \lambda \in \mathbb{C} \ .$$

Multiplying (2.3) by $e^{i\lambda x - \omega(\lambda)t}$ and integrating, we obtain

$$\int_{-\infty}^{\infty} e^{i\lambda x - \omega(\lambda)t}\hat{u}_0(\lambda)d\lambda - \int_{-\infty}^{\infty} e^{i\lambda x}\hat{u}(\lambda,t)d\lambda - \int_{-\infty}^{\infty} e^{i\lambda x - \omega(\lambda)t}\lambda^2\widetilde{g}_1(\omega(\lambda),t)d\lambda - i\int_{-\infty}^{\infty} e^{i\lambda x - \omega(\lambda)t}\lambda^3\widetilde{g}_0(\omega(\lambda),t)d\lambda$$

$$+ \int_{-\infty}^{\infty} e^{i\lambda x - \omega(\lambda)t}\widetilde{g}_3(\omega(\lambda),t)d\lambda + i\int_{-\infty}^{\infty} e^{i\lambda x - \omega(\lambda)t}\lambda\widetilde{g}_2(\omega(\lambda),t)d\lambda = 0 \ .$$

It follows that, for $x > 0$ and $t > 0$,

$$2\pi\iota(x,t) = \int_{-\infty}^{\infty} e^{i\lambda x - \omega(\lambda)t}\hat{u}_0(\lambda)d\lambda - i\int_{-\infty}^{\infty} e^{i\lambda x - \omega(\lambda)t}\lambda^3\widetilde{g}_0(\omega(\lambda),t)d\lambda - \int_{-\infty}^{\infty} e^{i\lambda x - \omega(\lambda)t}\lambda^2\widetilde{g}_1(\omega(\lambda),t)d\lambda$$

$$+ i\int_{-\infty}^{\infty} e^{i\lambda x - \omega(\lambda)t}\lambda\widetilde{g}_2(\omega(\lambda),t)d\lambda + \int_{-\infty}^{\infty} e^{i\lambda x - \omega(\lambda)t}\widetilde{g}_3(\omega(\lambda),t)d\lambda \ . \quad (2.4)$$

Now let us consider the sets

$$\Omega_1^- := \{\lambda \in \mathbb{C} : \operatorname{Im}\lambda \ge 0, \ \operatorname{Re}\lambda \le 0 \ and \ \operatorname{Re}\omega(\lambda) = \operatorname{Re}(\lambda^4) \le 0\} = \{\lambda \in \mathbb{C} : \frac{7\pi}{8} \le \lambda \le \frac{5\pi}{8}\}$$

and

$$\Omega_2^- := \{\lambda \in \mathbb{C} : \operatorname{Im}\lambda \ge 0, \ \operatorname{Re}\lambda \ge 0 \ and \ \operatorname{Re}\omega(\lambda) = \operatorname{Re}(\lambda^4) \le 0\} = \{\lambda \in \mathbb{C} : \frac{\pi}{8} \le \lambda \le \frac{3\pi}{8}\} \ ,$$

and their boundaries $\Gamma_1 = \partial\Omega_1^-$ and $\Gamma_2 = \partial\Omega_2^-$ (see fig.2).
Since

$$e^{-\omega(\lambda)t}\widetilde{g}_j(\omega(\lambda),t) = \frac{1}{\omega(\lambda)}g_j(t) - \frac{1}{\omega(\lambda)}e^{-\omega(\lambda)t}g_j(0) - \frac{1}{\omega(\lambda)}e^{-\omega(\lambda)t}(g_j{}')\widetilde{\ }(\omega(\lambda),t) \ , \ \forall\lambda \ne 0, \quad (2.5)$$

we have, uniformly for $t$ in compact subsets of $[0,+\infty)$, that

$$e^{-\omega(\lambda)t}\widetilde{g}_j(\omega(\lambda),t) = \mathcal{O}(1/\lambda^4) \ , \text{ as } \lambda \to \infty \text{ with } \lambda \in \mathbb{C} \text{ and } \operatorname{Re}\omega(\lambda) = \operatorname{Re}(\lambda^4) \ge 0 \ . \quad (2.6)$$

Thus, by Cauchy's theorem and Jordan's lemma,

$$\int_{-\infty}^{\infty} e^{i\lambda x - \omega(\lambda)t} \lambda^{3-j} \widetilde{g}_j(\omega(\lambda),t)d\lambda = \int_{\Gamma_1+\Gamma_2} e^{i\lambda x - \omega(\lambda)t} \lambda^{3-j} \widetilde{g}_j(\omega(\lambda),t)d\lambda \text{ , for } j=0,1,2,3 \text{ ,}$$

and, therefore, (2.4) becomes

$$2\pi u(x,t) = \int_{-\infty}^{\infty} e^{i\lambda x - \omega(\lambda)t} \hat{u}_0(\lambda)d\lambda - i \int_{\Gamma_1+\Gamma_2} e^{i\lambda x - \omega(\lambda)t} \lambda^3 \widetilde{g}_0(\omega(\lambda),t)d\lambda - \int_{\Gamma_1+\Gamma_2} e^{i\lambda x - \omega(\lambda)t} \lambda^2 \widetilde{g}_1(\omega(\lambda),t)d\lambda$$

$$+ i \int_{\Gamma_1+\Gamma_2} e^{i\lambda x - \omega(\lambda)t} \lambda \widetilde{g}_2(\omega(\lambda),t)d\lambda + \int_{\Gamma_1+\Gamma_2} e^{i\lambda x - \omega(\lambda)t} \widetilde{g}_3(\omega(\lambda),t)d\lambda \text{ . (2.7)}$$

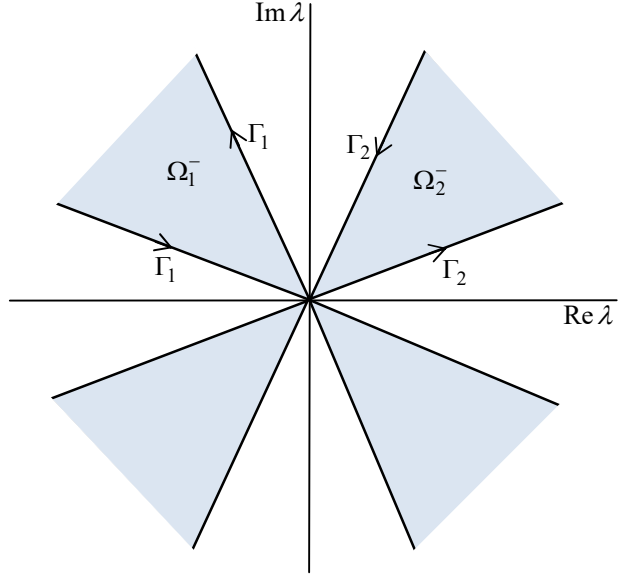

**Fig.2** $\Gamma_1 = \partial\Omega_1^- = \{\lambda : \arg\lambda = \frac{7\pi}{8} \text{ or } \frac{5\pi}{8}\}$ and $\Gamma_2 = \partial\Omega_2^- = \{\lambda : \arg\lambda = \frac{\pi}{8} \text{ or } \frac{3\pi}{8}\}$

In order to express the integrals

$$\int_{\Gamma_1} e^{i\lambda x - \omega(\lambda)t} \lambda \widetilde{g}_2(\omega(\lambda),t)d\lambda \text{ , } \int_{\Gamma_1} e^{i\lambda x - \omega(\lambda)t} \widetilde{g}_3(\omega(\lambda),t)d\lambda \text{ , } \int_{\Gamma_2} e^{i\lambda x - \omega(\lambda)t} \lambda \widetilde{g}_2(\omega(\lambda),t)d\lambda \text{ , } \int_{\Gamma_2} e^{i\lambda x - \omega(\lambda)t} \widetilde{g}_3(\omega(\lambda),t)d\lambda \text{ , (2.8)}$$

which appear in the RHS of (2.7), in terms of the data of problem (2.1), let us observe that

$$\omega(-\lambda) = \omega(i\lambda) = \omega(-i\lambda) = \omega(\lambda), \ \forall\lambda \text{ .}$$

Thus, setting $-\lambda$ (with $\text{Im}(-\lambda) \leq 0$) in place of $\lambda$ in (2.3), we obtain

$$\hat{u}_0(-\lambda) - \hat{u}(-\lambda,t)e^{\omega(\lambda)t} + \widetilde{g}_3(\omega(\lambda),t) - i\lambda\widetilde{g}_2(\omega(\lambda),t) - \lambda^2\widetilde{g}_1(\omega(\lambda),t) + i\lambda^3\widetilde{g}_0(\omega(\lambda),t) = 0 \text{ , for } \text{Im}\lambda \geq 0 \text{ ,}$$

which, in turn, implies that

$$\int_{\Gamma_k} e^{i\lambda x - \omega(\lambda)t} \hat{u}_0(-\lambda)d\lambda + \int_{\Gamma_k} e^{i\lambda x - \omega(\lambda)t} \widetilde{g}_3(\omega(\lambda),t)d\lambda - i\int_{\Gamma_k} e^{i\lambda x - \omega(\lambda)t} \lambda\widetilde{g}_2(\omega(\lambda),t)d\lambda$$

$$- \int_{\Gamma_k} e^{i\lambda x - \omega(\lambda)t} \lambda^2\widetilde{g}_1(\omega(\lambda),t)d\lambda + i\int_{\Gamma_k} e^{i\lambda x - \omega(\lambda)t} \lambda^3\widetilde{g}_0(\omega(\lambda),t)d\lambda = 0 \text{ , } k=1,2 \text{ , (2.9)}$$

where we also used the fact that $\int_{\Gamma_k} e^{i\lambda x} \hat{u}_0(-\lambda,t)d\lambda = 0$ .

Similarly, setting $i\lambda$ (with $\text{Im}(i\lambda) \leq 0$) in place of $\lambda$ in (2.3), we obtain

$$\hat{u}_0(i\lambda) - \hat{u}(i\lambda,t)e^{\omega(\lambda)t} + \widetilde{g}_3(\omega(\lambda),t) - \lambda\widetilde{g}_2(\omega(\lambda),t) + \lambda^2\widetilde{g}_1(\omega(\lambda),t) - \lambda^3\widetilde{g}_0(\omega(\lambda),t) = 0 \text{ , for } \lambda\in\Gamma_1 \text{ ,}$$

and this implies that

$$\int_{\Gamma_1} e^{i\lambda x - \omega(\lambda)t} \hat{u}_0(i\lambda)d\lambda + \int_{\Gamma_1} e^{i\lambda x - \omega(\lambda)t} \widetilde{g}_3(\omega(\lambda),t)d\lambda - \int_{\Gamma_1} e^{i\lambda x - \omega(\lambda)t} \lambda\widetilde{g}_2(\omega(\lambda),t)d\lambda$$

$$+ \int_{\Gamma_1} e^{i\lambda x - \omega(\lambda)t} \lambda^2\widetilde{g}_1(\omega(\lambda),t)d\lambda - \int_{\Gamma_1} e^{i\lambda x - \omega(\lambda)t} \lambda^3\widetilde{g}_0(\omega(\lambda),t)d\lambda = 0 \text{ . (2.10)}$$

Also, setting $-i\lambda$ (with $\mathrm{Im}(-i\lambda) \leq 0$) in place of $\lambda$ in (2.3), we obtain

$$\hat{u}_0(-i\lambda) - \hat{u}(-i\lambda,t)e^{\omega(\lambda)t} + \widetilde{g}_3(\omega(\lambda),t) + \lambda\widetilde{g}_2(\omega(\lambda),t) + \lambda^2\widetilde{g}_1(\omega(\lambda),t) + \lambda^3\widetilde{g}_0(\omega(\lambda),t) = 0 \text{, for } \lambda \in \Gamma_2\text{,}$$

and this implies that

$$\int_{\Gamma_2} e^{i\lambda x - \omega(\lambda)t}\hat{u}_0(-i\lambda)d\lambda + \int_{\Gamma_2} e^{i\lambda x - \omega(\lambda)t}\widetilde{g}_3(\omega(\lambda),t)d\lambda + \int_{\Gamma_2} e^{i\lambda x - \omega(\lambda)t}\lambda\widetilde{g}_2(\omega(\lambda),t)d\lambda$$
$$+ \int_{\Gamma_2} e^{i\lambda x - \omega(\lambda)t}\lambda^2\widetilde{g}_1(\omega(\lambda),t)d\lambda + \int_{\Gamma_2} e^{i\lambda x - \omega(\lambda)t}\lambda^3\widetilde{g}_0(\omega(\lambda),t)d\lambda = 0 \text{. (2.11)}$$

Solving the system of equations {(2.9) (with $k=1$) & (2.10)}, we determine the values of the integrals (2.8) which are taken over $\Gamma_1$, while solving the system of equations {(2.9) (with $k=2$) & (2.10)}, we determine the values of the integrals (2.8) which are taken over $\Gamma_2$. Substituting these values of the integrals (2.8) in (2.7), we obtain the solution of problem (2.1): For $x > 0$ and $t > 0$,

$$2\pi u(x,t) = \int_{-\infty}^{\infty} e^{i\lambda x - \omega(\lambda)t}\hat{u}_0(\lambda)d\lambda$$
$$+ \int_{\Gamma_1} e^{i\lambda x - \omega(\lambda)t}[(2-2i)\lambda^3\widetilde{g}_0(\omega(\lambda),t) + (-2+2i)\lambda^2\widetilde{g}_1(\omega(\lambda),t) - i\hat{u}_0(-\lambda) + (i-1)\hat{u}_0(i\lambda)]d\lambda$$
$$+ \int_{\Gamma_2} e^{i\lambda x - \omega(\lambda)t}[(-2-2i)\lambda^3\widetilde{g}_0(\omega(\lambda),t) + (-2-2i)\lambda^2\widetilde{g}_1(\omega(\lambda),t) + i\hat{u}_0(-\lambda) - (i+1)\hat{u}_0(-i\lambda)]d\lambda \text{. (2.12)}$$

The above formula can be written also in the following form: For fixed $T > 0$,

$$2\pi u(x,t) = \int_{-\infty}^{\infty} e^{i\lambda x - \omega(\lambda)t}\hat{u}_0(\lambda)d\lambda$$
$$+ \int_{\Gamma_1} e^{i\lambda x - \omega(\lambda)t}[(2-2i)\lambda^3\widetilde{g}_0(\omega(\lambda),T) + (-2+2i)\lambda^2\widetilde{g}_1(\omega(\lambda),T) - i\hat{u}_0(-\lambda) + (i-1)\hat{u}_0(i\lambda)]d\lambda$$
$$+ \int_{\Gamma_2} e^{i\lambda x - \omega(\lambda)t}[(-2-2i)\lambda^3\widetilde{g}_0(\omega(\lambda),T) + (-2-2i)\lambda^2\widetilde{g}_1(\omega(\lambda),T) + i\hat{u}_0(-\lambda) - (i+1)\hat{u}_0(-i\lambda)]d\lambda \text{, (2.13)}$$

for $x > 0$ and $0 < t < T$. This expresses the solution $u(x,t)$ in Ehrenpreis form in the rectangle

$$Q_T := \mathrm{int}(\Pi_T) = \{(x,t) \in \mathbb{R} \times \mathbb{R} : x > 0 \text{ and } 0 < t < T\}\text{.}$$

The integrals in the RHS of (2.13), taken over $\Gamma_1$ and $\Gamma_2$, can be differentiated, with respect to $x$ or $t$, under the integral sign, because of the presence of the factor $e^{i\lambda x}$, since

$$\left|e^{i\lambda x}\right| = e^{-[x\sin(\pi/8)]|\lambda|} \text{ or } \left|e^{i\lambda x}\right| = e^{-[x\sin(3\pi/8)]|\lambda|}\text{, for } \lambda \in \Gamma_1 \cup \Gamma_2\text{.} \tag{2.14}$$

Also, the first integral in the RHS of (2.13), taken for $\lambda \in \mathbb{R}$, can be differentiated, with respect to $x$ or $t$, under the integral sign, because of the presence of the factor $e^{-\omega(\lambda)t} = e^{-\lambda^4 t}$.

Hence, for fixed $T > 0$,

$$2\pi\frac{\partial^{n+m}u(x,t)}{\partial x^n\partial t^m} = \int_{-\infty}^{\infty}(i\lambda)^n[-\omega(\lambda)]^m e^{i\lambda x - \omega(\lambda)t}\hat{u}_0(\lambda)d\lambda$$
$$+ \int_{\Gamma_1}(i\lambda)^n[-\omega(\lambda)]^m e^{i\lambda x - \omega(\lambda)t}[(2-2i)\lambda^3\widetilde{g}_0(\omega(\lambda),T) + (-2+2i)\lambda^2\widetilde{g}_1(\omega(\lambda),T) - i\hat{u}_0(-\lambda) + (i-1)\hat{u}_0(i\lambda)]d\lambda$$
$$+ \int_{\Gamma_2}(i\lambda)^n[-\omega(\lambda)]^m e^{i\lambda x - \omega(\lambda)t}[(-2-2i)\lambda^3\widetilde{g}_0(\omega(\lambda),T) + (-2-2i)\lambda^2\widetilde{g}_1(\omega(\lambda),T) + i\hat{u}_0(-\lambda) - (i+1)\hat{u}_0(-i\lambda)]d\lambda \text{, (2.15)}$$

for $x > 0$ and $0 < t < T$.

Furthermore, the integrals in the RHS of (2.15) converge absolutely and uniformly in compact subsets of the rectangle $Q_T$. In particular, the function $u(x,t)$, defined by (2.13) (equivalently by (2.12)), is $C^{\infty}$ and satisfies the equation $\partial_t u = -\partial_{xxxx}u$. Thus, we have proved the first part of the following theorem.

**Theorem 1** $1^{st}$ *The function $u(x,t)$, defined by (2.13) (equivalently by (2.12)) is $C^{\infty}$ (jointly) for $(x,t) \in Q$ and satisfies the differential equation $\partial_t u = -\partial_x^4 u$.*

$2^{nd}$ *For fixed* $x > 0$, $\lim\limits_{t \to 0^+} u(x,t) = u_0(x)$. (2.16)

$3^{rd}$ *For fixed* $t > 0$, $\lim\limits_{x \to 0^+} u(x,t) = g_0(t)$. (2.17)

$4^{th}$ *For fixed* $t > 0$, $\lim\limits_{x \to 0^+} \dfrac{\partial u(x,t)}{\partial x} = g_1(t)$. (2.18)

$5^{th}$ *The convergence in (2.16) is uniform for* $x \geq x_0$ $(\forall x_0 > 0)$ *and the convergence in (2.17) and (2.18) is uniform for* $t \geq t_0$ $(\forall t_0 > 0)$.

**Comment** As we pointed out, the integrals in (2.12), (2.13) and (2.15) converge absolutely, as long as $x > 0$ and $t > 0$. However, if $x = 0$ or $t = 0$, some of these integrals do not converge. Thus, there is some difficulty in taking the limit of these integrals as $x \to 0^+$ or $t \to 0^+$, in the sense that it is not immediately clear that we can interchange the limit with the integral, even in the case this holds. As a matter of fact, in the integrals in the RHS of (1.12), which involve the function $\widetilde{g}_0$, we cannot make this interchange, in the sense that

$$\lim_{x \to 0^+} \int_{\Gamma_k} e^{i\lambda x - \omega(\lambda)t} \lambda^3 \widetilde{g}_0(\omega(\lambda),t) d\lambda \neq \int_{\Gamma_k} e^{-\omega(\lambda)t} \lambda^3 \widetilde{g}_0(\omega(\lambda),t) d\lambda ,$$

in general. (See *remark 2.4*.) As we will see, the last integral, above, converges but not absolutely.

**2.2 Jordan type lemmas** Let $t > 0$ and $\psi(\lambda)$ be a continuous function such that $\lim\limits_{A \to \infty} [\sup\{|\psi(\lambda)| : |\lambda| = A\}] = 0$. Then

(1) $\lim\limits_{R \to \infty} \int\limits_{\{|\lambda| = R \,\&\, -\frac{\pi}{8} \leq \arg \lambda \leq \frac{\pi}{8}\}} \lambda^3 e^{-\omega(\lambda)t} \psi(\lambda) d\lambda = 0$ and (2) $\lim\limits_{R \to \infty} \int\limits_{\{|\lambda| = R \,\&\, \frac{\pi}{8} \leq \arg \lambda \leq \frac{3\pi}{8}\}} \lambda^3 e^{\omega(\lambda)t} \psi(\lambda) d\lambda = 0$.

These follow from the classical Jordan lemma. Indeed, it suffices to notice that, if we set $\kappa = i\lambda^4$, the integral in (1) becomes

$$\frac{1}{4i} \int\limits_{\{|\kappa| = R^4 \,\&\, 0 \leq \arg \kappa \leq \pi\}} e^{i\kappa t} \psi(\sqrt[4]{-i\kappa}) d\kappa ,$$

while, if we set $\kappa = -i\lambda^4$, the integral in (2) becomes

$$\frac{1}{-4i} \int\limits_{\{|\kappa| = R^4 \,\&\, 0 \leq \arg \kappa \leq \pi\}} e^{i\kappa t} \psi(\sqrt[4]{i\kappa}) d\kappa .$$

Similarly,

$$\lim_{R \to \infty} \int\limits_{\{|\lambda| = R \,\&\, \frac{3\pi}{8} \leq \arg \lambda \leq \frac{5\pi}{8}\}} \lambda^3 e^{-\omega(\lambda)t} \psi(\lambda) d\lambda = 0, \quad \lim_{R \to \infty} \int\limits_{\{|\lambda| = R \,\&\, \frac{5\pi}{8} \leq \arg \lambda \leq \frac{7\pi}{8}\}} \lambda^3 e^{\omega(\lambda)t} \psi(\lambda) d\lambda = 0, \quad \lim_{R \to \infty} \int\limits_{\{\lambda = Re^{i\theta}, \, \frac{7\pi}{8} \leq \theta \leq \frac{9\pi}{8}\}} \lambda^3 e^{-\omega(\lambda)t} \psi(\lambda) d\lambda = 0.$$

**2.3 Proof of Theorem 1** *Step 1* We claim that, for fixed $x > 0$,

$$\lim_{t \to 0^+} \int_{-\infty}^{\infty} e^{i\lambda x - \omega(\lambda)t} \hat{u}_0(\lambda) d\lambda = \int_{-\infty}^{\infty} e^{i\lambda x} \hat{u}_0(\lambda) d\lambda = 2\pi u_0(x) . \quad (2.19)$$

Indeed, in view of (1.5),

$$\int_{-\infty}^{\infty} e^{i\lambda x - \omega(\lambda)t} \hat{u}_0(\lambda) d\lambda = \int_{-1}^{1} e^{i\lambda x - \omega(\lambda)t} \hat{u}_0(\lambda) d\lambda + u_0(0) \left( \int_{-\infty}^{-1} + \int_{1}^{\infty} \right) e^{i\lambda x - \omega(\lambda)t} \frac{d\lambda}{i\lambda} + \left( \int_{-\infty}^{-1} + \int_{1}^{\infty} \right) e^{i\lambda x - \omega(\lambda)t} \frac{\widehat{(u_0')}(\lambda)}{i\lambda} d\lambda$$

$$= \int_{-1}^{1} e^{i\lambda x - \omega(\lambda)t} \hat{u}_0(\lambda) d\lambda + u_0(0) \left( \int_{\gamma^-} + \int_{\gamma^+} \right) e^{i\lambda x - \omega(\lambda)t} \frac{d\lambda}{i\lambda} + \left( \int_{-\infty}^{-1} + \int_{1}^{\infty} \right) e^{i\lambda x - \omega(\lambda)t} \frac{\widehat{(u_0')}(\lambda)}{i\lambda} d\lambda ,$$

where

$\gamma^- := (\infty e^{i7\pi/8}, e^{i7\pi/8}] + \{\lambda : |\lambda| = 1 \text{ and } \frac{7\pi}{8} \leq \arg \lambda \leq \pi\}$ and $\gamma^+ := \{\lambda : |\lambda| = 1 \text{ and } 0 \leq \arg \lambda \leq \frac{\pi}{8}\} + [e^{i\pi/8}, \infty e^{i\pi/8})$.

Now we may let $t \to 0^+$ obtaining (2.19):

$$\lim_{t \to 0^+} \int_{-\infty}^{\infty} e^{i\lambda x - \omega(\lambda)t} \hat{u}_0(\lambda) d\lambda = \int_{-1}^{1} e^{i\lambda x} \hat{u}_0(\lambda) d\lambda + u_0(0) \left( \int_{\gamma^-} + \int_{\gamma^+} \right) e^{i\lambda x} \frac{d\lambda}{i\lambda} + \left( \int_{-\infty}^{-1} + \int_{1}^{\infty} \right) e^{i\lambda x} \frac{\widehat{(u_0')}(\lambda)}{i\lambda} d\lambda$$

$$= \int\limits_{-1}^{1} e^{i\lambda x}\hat{u}_0(\lambda)d\lambda + u_0(0)\left(\int\limits_{-\infty}^{-1} + \int\limits_{1}^{\infty}\right)e^{i\lambda x}\frac{d\lambda}{i\lambda} + \left(\int\limits_{-\infty}^{-1} + \int\limits_{1}^{\infty}\right)e^{i\lambda x}\frac{(u_0)\hat{}(\lambda)}{i\lambda}d\lambda = \int\limits_{-\infty}^{\infty} e^{i\lambda x}\hat{u}_0(\lambda)d\lambda = 2\pi u_0(x) \,.$$

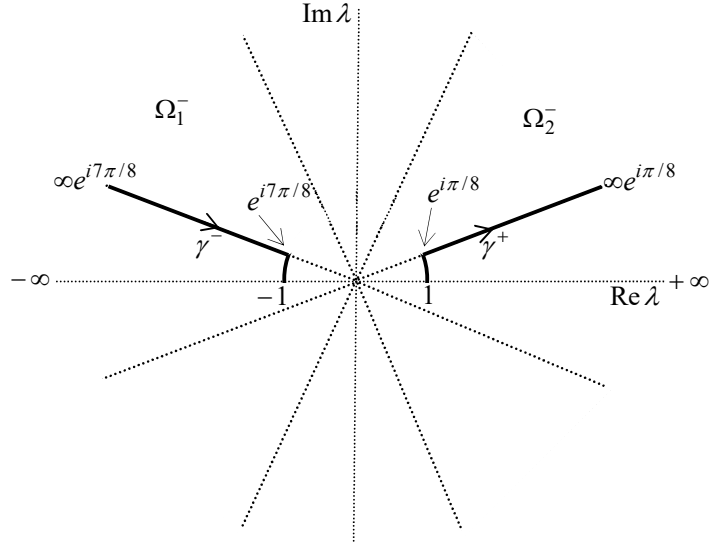

**Fig.3** The contours $\gamma^-$ and $\gamma^+$ used in the deformation procces in the proof of (2.19)

**Step 2** *Proof of $2^{nd}$ part* Using the expression of the integrals $\int_{\Gamma_1} e^{i\lambda x - \omega(\lambda)t} \cdots$ and $\int_{\Gamma_2} e^{i\lambda x - \omega(\lambda)t} \cdots$, as they appear in (2.12), it is immediate that

$$\lim_{t \to 0^+} \int\limits_{\Gamma_1} e^{i\lambda x - \omega(\lambda)t} \cdots = 0 \quad \text{and} \quad \lim_{t \to 0^+} \int\limits_{\Gamma_2} e^{i\lambda x - \omega(\lambda)t} \cdots = 0 \,, \tag{2.20}$$

because of (2.14), (2.6), (1.5), Cauchy's theorem and Jordan's lemma.

Therefore, (2.16) follows from (2.19) and (2.20).

**Step 3** We claim that

$$\lim_{x \to 0^+}\left[\int\limits_{-\infty}^{\infty} e^{i\lambda x - \omega(\lambda)t}\hat{u}_0(\lambda)d\lambda + \int\limits_{\Gamma_1} e^{i\lambda x - \omega(\lambda)t}[-i\hat{u}_0(-\lambda) + (i-1)\hat{u}_0(i\lambda)]d\lambda + \int\limits_{\Gamma_2}[i\hat{u}_0(-\lambda) - (i+1)\hat{u}_0(-i\lambda)]d\lambda\right] = 0 \,. \tag{2.21}$$

Indeed, this follows from the following equations:

$$\lim_{x \to 0^+} \int\limits_{\Gamma_1} e^{i\lambda x - \omega(\lambda)t}[-i\hat{u}_0(-\lambda)]d\lambda = \lim_{x \to 0^+}\left[-i\int\limits_{-\infty}^{0} e^{i\lambda x - \omega(\lambda)t}\hat{u}_0(-\lambda)d\lambda + \int\limits_{0}^{\infty} e^{-\lambda x - \omega(\lambda)t}\hat{u}_0(-i\lambda)d\lambda\right]$$

$$= -i\int\limits_{0}^{\infty} e^{-\omega(\lambda)t}\hat{u}_0(\lambda)d\lambda + \int\limits_{-\infty}^{0} e^{-\omega(\lambda)t}\hat{u}_0(i\lambda)d\lambda \,,$$

$$\lim_{x \to 0^+} \int\limits_{\Gamma_1} e^{i\lambda x - \omega(\lambda)t}[(i-1)\hat{u}_0(i\lambda)]d\lambda = (i-1)\int\limits_{-\infty}^{0} e^{-\omega(\lambda)t}\hat{u}_0(i\lambda)d\lambda + i(i-1)\int\limits_{-\infty}^{0} e^{-\omega(\lambda)t}\hat{u}_0(\lambda)d\lambda \,,$$

$$\lim_{x \to 0^+} \int\limits_{\Gamma_2} e^{i\lambda x - \omega(\lambda)t}[i\hat{u}_0(-\lambda)]d\lambda = i\int\limits_{0}^{\infty} e^{-\omega(\lambda)t}\hat{u}_0(-\lambda)d\lambda + i\int\limits_{0}^{\infty} e^{-\omega(\lambda)t}\hat{u}_0(-i\lambda)d\lambda = i\int\limits_{-\infty}^{0} e^{-\omega(\lambda)t}\hat{u}_0(\lambda)d\lambda + i\int\limits_{-\infty}^{0} e^{-\omega(\lambda)t}\hat{u}_0(i\lambda)d\lambda \,,$$

$$\lim_{x \to 0^+} \int\limits_{\Gamma_2} e^{i\lambda x - \omega(\lambda)t}[-(i+1)\hat{u}_0(-i\lambda)]d\lambda = -(i+1)\int\limits_{0}^{\infty} e^{-\omega(\lambda)t}\hat{u}_0(-i\lambda)d\lambda + (i+1)i\int\limits_{0}^{\infty} e^{-\omega(\lambda)t}\hat{u}_0(\lambda)d\lambda$$

$$= -(i+1)\int\limits_{-\infty}^{0} e^{-\omega(\lambda)t}\hat{u}_0(i\lambda)d\lambda + (i+1)i\int\limits_{0}^{\infty} e^{-\omega(\lambda)t}\hat{u}_0(\lambda)d\lambda \,.$$

(Some of the above equations follow from Cauchy's theorem and Jordan's lemma.)

**Step 4** We claim that

$$\lim_{x \to 0^+}\left[\int_{\Gamma_1}(-2+2i)e^{i\lambda x-\omega(\lambda)t}\lambda^2\widetilde{g}_1(\omega(\lambda),t)d\lambda + \int_{\Gamma_2}(-2-2i)e^{i\lambda x-\omega(\lambda)t}\lambda^2\widetilde{g}_1(\omega(\lambda),t)d\lambda\right]=0 \ . \tag{2.22}$$

Indeed, in view of (2.6), $e^{-\omega(\lambda)t}\lambda^2\widetilde{g}_1(\omega(\lambda),t)=\mathrm{O}(1/\lambda^2)$ for $\lambda \in \Gamma_1 \cup \Gamma_2$ ($\lambda \to \infty$), hence the limit in (2.22) is equal to

$$(-2+2i)\int_{\Gamma_1}e^{-\omega(\lambda)t}\lambda^2\widetilde{g}_1(\omega(\lambda),t)d\lambda + (-2-2i)\int_{\Gamma_2}e^{-\omega(\lambda)t}\lambda^2\widetilde{g}_1(\omega(\lambda),t)d\lambda \ . \tag{2.23}$$

Setting $\mu=-i\lambda$ in the integral $\int_{\Gamma_1}\cdots d\lambda$ which appears in (2.23), we see that $\int_{\Gamma_1}\cdots = -i\int_{\Gamma_2}\cdots$ (for the integrals in (2.23)), and (2.22) follows.

**Step 5** We claim that, for $0<t<T$ ,

$$\left[(2-2i)\int_{\Gamma_1}e^{i\lambda x-\omega(\lambda)t}\lambda^3\widetilde{g}_0(\omega(\lambda),T)d\lambda + (-2-2i)\int_{\Gamma_2}e^{i\lambda x-\omega(\lambda)t}\lambda^3\widetilde{g}_0(\omega(\lambda),T)d\lambda\right]\Bigg|_{x=0}=2\pi g_0(t) \ . \tag{2.24}$$

The LHS of (2.24) is the evaluation of the quantity inside the square brackets and is equal to

$$(2-2i)\int_{\Gamma_1}e^{-\omega(\lambda)t}\lambda^3\widetilde{g}_0(\omega(\lambda),T)d\lambda + (-2-2i)\int_{\Gamma_2}e^{-\omega(\lambda)t}\lambda^3\widetilde{g}_0(\omega(\lambda),T)d\lambda \ . \tag{2.25}$$

(*Comment* The integrals in (2.25) do not – in general – converge absolutely. However, they exist as generalized integrals, i.e., as limits of the form $\lim_{A\to\infty}\int_{\Gamma_k\cap\{\lambda:|\lambda|\le A\}}$, by Cauchy's theorem and Jordan's lemma.)

Setting $\mu=-i\lambda$ , we obtain that

$$\int_{\Gamma_1}e^{-\omega(\lambda)t}\lambda^3\widetilde{g}_0(\omega(\lambda),T)d\lambda = \int_{\Gamma_2}e^{-\omega(\mu)t}\mu^3\widetilde{g}_0(\omega(\mu),T)d\mu \ .$$

Therefore, the quantity (2.25) is equal to

$$-4i\int_{\Gamma_2}e^{-\omega(\lambda)t}\lambda^3\widetilde{g}_0(\omega(\lambda),T)d\lambda = \int_{-\infty}^{\infty}e^{i\kappa t}\left(\int_{\tau=0}^{T}e^{-i\kappa\tau}g_0(\tau)d\tau\right)d\kappa = 2\pi g_0(t) \ ,$$

where we changed the variable in the $d\lambda$-integral, by setting $\kappa=-i\lambda^4$ , and we used Fourier's inversion formula for the function $\varphi(\tau)$ defined as follows:

$$\varphi(\tau):=\begin{cases}g_0(\tau) \ \textit{if} \ 0\le\tau\le T \\ 0 \quad \textit{otherwise}.\end{cases}$$

This proves (2.24).

**Step 6** Fixing $0<t<T$ , we will show that

$$\lim_{x\to 0^+}\int_{\Gamma_k}e^{i\lambda x-\omega(\lambda)t}\lambda^3\widetilde{g}_0(\omega(\lambda),T)d\lambda = \int_{\Gamma_k}e^{-\omega(\lambda)t}\lambda^3\widetilde{g}_0(\omega(\lambda),T)d\lambda \ , \ k=1,2 \ . \tag{2.26}$$

Since

$$e^{-\omega(\lambda)t}\widetilde{g}_j(\omega(\lambda),T) = \frac{1}{\omega(\lambda)}e^{\omega(\lambda)(T-t)}g_j(T) - \frac{1}{\omega(\lambda)}e^{-\omega(\lambda)t}g_j(0) - \frac{1}{\omega(\lambda)}e^{-\omega(\lambda)t}(g_j{}')\widetilde{\phantom{i}}(\omega(\lambda),T) \ , \ \forall\lambda\ne 0 \ , \tag{2.27}$$

and $\mathrm{Re}[\omega(\lambda)]=0$ for $\lambda\in\Gamma_1\cup\Gamma_2$ , it follows that

$$e^{-\omega(\lambda)t}\widetilde{g}_j(\omega(\lambda),T)=\mathrm{O}(1/\lambda^4) \quad \text{and} \quad \frac{1}{\omega(\lambda)}e^{-\omega(\lambda)t}(g_j{}')\widetilde{\phantom{i}}(\omega(\lambda),T)=\mathrm{O}(1/\lambda^8) \ , \ \text{as } \lambda\to\infty \ \text{with} \ \lambda\in\Gamma_1\cup\Gamma_2 \ .$$

Therefore, in order to prove (2.26), it suffices to show that

$$\lim_{x\to 0^+}\int_{\Gamma_k\cap\{|\lambda|\ge 1\}}e^{i\lambda x}e^{\omega(\lambda)(T-t)}\frac{d\lambda}{\lambda} = \int_{\Gamma_k\cap\{|\lambda|\ge 1\}}e^{\omega(\lambda)(T-t)}\frac{d\lambda}{\lambda} \tag{2.28}$$

and

$$\lim_{x\to 0^+}\int_{\Gamma_k\cap\{|\lambda|\ge 1\}}e^{i\lambda x-\omega(\lambda)t}\frac{d\lambda}{\lambda} = \int_{\Gamma_k\cap\{|\lambda|\ge 1\}}e^{-\omega(\lambda)t}\frac{d\lambda}{\lambda} \ . \tag{2.29}$$

*Proof of (2.28)* We will prove it for $k=2$ – the proof for $k=1$ is similar. We have

$$\lim_{x\to 0^+}\int_{\Gamma_2\cap\{|\lambda|\geq 1\}}e^{i\lambda x}e^{\omega(\lambda)(T-t)}\frac{d\lambda}{\lambda}=\lim_{x\to 0^+}\int_{\{\frac{\pi}{8}\leq\arg\lambda\leq\frac{3\pi}{8}\}\cap\{|\lambda|=1\}}e^{i\lambda x}e^{\omega(\lambda)(T-t)}\frac{d\lambda}{\lambda}$$

$$=\int_{\{\frac{\pi}{8}\leq\arg\lambda\leq\frac{3\pi}{8}\}\cap\{|\lambda|=1\}}e^{\omega(\lambda)(T-t)}\frac{d\lambda}{\lambda}=\int_{\Gamma_2\cap\{|\lambda|\geq 1\}}e^{\omega(\lambda)(T-t)}\frac{d\lambda}{\lambda},$$

where we used Jordan's lemma.

*Proof of (2.29)* We will prove it for $k=2$ – the proof for $k=1$ is similar. We have

$$\lim_{x\to 0^+}\int_{\Gamma_2\cap\{|\lambda|\geq 1\}}e^{i\lambda x-\omega(\lambda)t}\frac{d\lambda}{\lambda}$$

$$=\lim_{x\to 0^+}\left[\int_{\{|\lambda|=1\}\cap\{0\leq\arg\lambda\leq\frac{\pi}{8}\}+\{|\lambda|\geq 1\}\cap\{\arg\lambda=0\}}e^{i\lambda x-\omega(\lambda)t}\frac{d\lambda}{\lambda}+\int_{\{|\lambda|=1\}\cap\{\frac{3\pi}{8}\leq\arg\lambda\leq\frac{\pi}{2}\}+\{|\lambda|\geq 1\}\cap\{\arg\lambda=\frac{\pi}{2}\}}e^{i\lambda x-\omega(\lambda)t}\frac{d\lambda}{\lambda}\right]$$

$$=\int_{\{|\lambda|=1\}\cap\{0\leq\arg\lambda\leq\frac{\pi}{8}\}+\{|\lambda|\geq 1\}\cap\{\arg\lambda=0\}}e^{-\omega(\lambda)t}\frac{d\lambda}{\lambda}+\int_{\{|\lambda|=1\}\cap\{\frac{3\pi}{8}\leq\arg\lambda\leq\frac{\pi}{2}\}+\{|\lambda|\geq 1\}\cap\{\arg\lambda=\frac{\pi}{2}\}}e^{-\omega(\lambda)t}\frac{d\lambda}{\lambda}$$

$$=\int_{\{|\lambda|=1\}\cap\{\arg\lambda=\frac{\pi}{8}\}}e^{-\omega(\lambda)t}\frac{d\lambda}{\lambda}+\int_{\{|\lambda|\geq 1\}\cap\{\arg\lambda=\frac{3\pi}{8}\}}e^{-\omega(\lambda)t}\frac{d\lambda}{\lambda}=\int_{\Gamma_2\cap\{|\lambda|\geq 1\}}e^{i\lambda x-\omega(\lambda)t}\frac{d\lambda}{\lambda}.$$

This completes the proof of (2.26).

**Step 7** *Proof of $3^{rd}$ part* It follows from (2.26), (2.24), (2.22), (2.21), (2.19) and (2.13).

**Step 8** By (2.15),

$$2\pi\frac{\partial u(x,t)}{\partial x}=\int_{-\infty}^{\infty}i\lambda e^{i\lambda x-\omega(\lambda)t}\hat{u}_0(\lambda)d\lambda$$

$$+\int_{\Gamma_1}e^{i\lambda x-\omega(\lambda)t}[(2i+2)\lambda^4\widetilde{g}_0(\omega(\lambda),t)+(-2i-2)\lambda^3\widetilde{g}_1(\omega(\lambda),t)+\lambda\hat{u}_0(-\lambda)-(1+i)\lambda\hat{u}_0(i\lambda)]d\lambda$$

$$+\int_{\Gamma_2}e^{i\lambda x-\omega(\lambda)t}[(-2i+2)\lambda^4\widetilde{g}_0(\omega(\lambda),t)+(-2i+2)\lambda^3\widetilde{g}_1(\omega(\lambda),t)-\lambda\hat{u}_0(-\lambda)+(1-i)\lambda\hat{u}_0(-i\lambda)]d\lambda .\ (2.30)$$

Working as in step 3, we show that

$$\lim_{x\to 0^+}\left[\int_{-\infty}^{\infty}i\lambda e^{i\lambda x-\omega(\lambda)t}\hat{u}_0(\lambda)d\lambda+\int_{\Gamma_1}e^{i\lambda x-\omega(\lambda)t}[\lambda\hat{u}_0(-\lambda)-(1+i)\lambda\hat{u}_0(i\lambda)]d\lambda-\int_{\Gamma_2}[\lambda\hat{u}_0(-\lambda)-(1-i)\lambda\hat{u}_0(-i\lambda)]d\lambda\right]=0 .\ (2.31)$$

Working as in steps 5 and 6, we show that

$$\lim_{x\to 0^+}\left[\int_{\Gamma_1}(-2i-2)e^{i\lambda x-\omega(\lambda)t}\lambda^3\widetilde{g}_1(\omega(\lambda),t)d\lambda+(-2i+2)\int_{\Gamma_2}e^{i\lambda x-\omega(\lambda)t}\lambda^3\widetilde{g}_1(\omega(\lambda),t)d\lambda\right]=2\pi g_1(t) .\ (2.32)$$

Finally, we claim that

$$\lim_{x\to 0^+}\left[\int_{\Gamma_1}(2i+2)e^{i\lambda x-\omega(\lambda)t}\lambda^4\widetilde{g}_0(\omega(\lambda),T)d\lambda+\int_{\Gamma_2}(-2i+2)e^{i\lambda x-\omega(\lambda)t}\lambda^4\widetilde{g}_0(\omega(\lambda),T)d\lambda\right]=0 .\ (2.33)$$

Indeed, if we could interchange the order of the limit and the integrals, we would have to show that

$$\int_{\Gamma_1}(2i+2)e^{-\omega(\lambda)t}\lambda^4\widetilde{g}_0(\omega(\lambda),T)d\lambda+\int_{\Gamma_2}(-2i+2)e^{-\omega(\lambda)t}\lambda^4\widetilde{g}_0(\omega(\lambda),T)d\lambda=0 ,\ (2.34)$$

and then we could work as in the proof of (2.22). However, in most cases, this is not immediately clear. In addition, we have to interpret the integrals in (1.34), i.e., to ensure their convergence – in the generalized sense. In order to overcome these difficulties, we use (2.27) and work as in step 6.

**Step 9** *Proof of $4^{th}$ part* It follows from (2.30)−(2.34).

**Step 10** *Proof of $5^{th}$ part* Examining the limiting processes in steps 1-9, we easily obtain the assertion concerning the uniformity of the limits.

**2.4 Remark** Let us point out that

$$\left[(2-2i)\int_{\Gamma_1}e^{i\lambda x-\omega(\lambda)t}\lambda^3\widetilde{g}_0(\omega(\lambda),t)d\lambda+(-2-2i)\int_{\Gamma_2}e^{i\lambda x-\omega(\lambda)t}\lambda^3\widetilde{g}_0(\omega(\lambda),t)d\lambda\right]\Bigg|_{x=0}$$

$$=(2-2i)\int_{\Gamma_1}e^{-\omega(\lambda)t}\lambda^3\widetilde{g}_0(\omega(\lambda),t)d\lambda+(-2-2i)\int_{\Gamma_2}e^{-\omega(\lambda)t}\lambda^3\widetilde{g}_0(\omega(\lambda),t)d\lambda=2\pi\left[\frac{1}{2}g_0(t)\right].$$

Indeed, this follows as (2.24), but this time we apply Fourier's inversion formula for the function

$$\psi(\tau):=\begin{cases}g_0(\tau) & \text{if } 0\le\tau\le t\\ 0 & \text{otherwise}.\end{cases}$$

This, combined with (2.26), implies that, if $g_0(t)\ne0$, then

$$\lim_{x\to0^+}\int_{\Gamma_k}e^{i\lambda x-\omega(\lambda)t}\lambda^3\widetilde{g}_0(\omega(\lambda),t)d\lambda\ne\int_{\Gamma_k}e^{-\omega(\lambda)t}\lambda^3\widetilde{g}_0(\omega(\lambda),t)d\lambda\,,$$

i.e. (2.26) does not hold if we replace $T$ by $t$.

**Theorem 2** *The function $u(x,t)$, defined by (2.12) for $(x,t)\in Q$, satisfies the following:*

*$1^{st}$ The limits*

$$W_{n,m}(x):=\lim_{t\to0^+}\frac{\partial^{n+m}u(x,t)}{\partial x^n\partial t^m}\,,\ n,m=0,1,2,...,$$

*exist and are uniform for $x\ge x_0$ ($\forall x_0>0$). Moreover the functions $W_{n,m}(x)$ are $C^\infty$ for $x\in(0,\infty)$.*

*$2^{nd}$ The limits*

$$V_{n,m}(t):=\lim_{x\to0^+}\frac{\partial^{n+m}u(x,t)}{\partial x^n\partial t^m}\,,\ n,m=0,1,2,...,$$

*exist and are uniform for $t\ge t_0$ ($\forall t_0>0$). Moreover the functions $V_{n,m}(t)$ are $C^\infty$ for $t\in(0,\infty)$.*

*$3^{rd}$ The function $u(x,t)$, extended to $\overline{Q}-\{(0,0)\}$ by setting*

$$u*(x,0):=u_0(x)\text{ for }x>0\text{ and }u*(0,t):=g_0(t)\text{ for }t>0,$$

*is $C^\infty$ (in $\overline{Q}-\{(0,0)\}$).*

*$4^{th}$ The limit $\displaystyle\lim_{t\to0^+}\frac{\partial^n u(x,t)}{\partial x^n}=\frac{d^n u_0(x)}{dx^n}$, uniformly for $x\ge x_0$.*

*$5^{th}$ The limit $\displaystyle\lim_{x\to0^+}\frac{\partial^m u(x,t)}{\partial t^m}=\frac{d^m g_0(t)}{dt^m}$, uniformly for $t\ge t_0$.*

*$6^{th}$ The limit $\displaystyle\lim_{x\to0^+}\frac{\partial^{m+1}u(x,t)}{\partial t^m\partial x}=\frac{d^m g_1(t)}{dt^m}$, uniformly for $t\ge t_0$.*

**Proof** *Step 1* We claim that the limit

$$\lim_{t\to0^+}\frac{\partial^{n+m}}{\partial x^n\partial t^m}\left(\int_{-\infty}^\infty e^{i\lambda x-\omega(\lambda)t}\hat{u}_0(\lambda)d\lambda\right)=\lim_{t\to0^+}\int_{-\infty}^\infty(i\lambda)^n[-\omega(\lambda)]^m e^{i\lambda x-\omega(\lambda)t}\hat{u}_0(\lambda)d\lambda\qquad(2.35)$$

exists and defines a $C^\infty$ function for $x>0$. Indeed, this follows from the formula

$$\int_{-\infty}^\infty e^{i\lambda x-\omega(\lambda)t}\hat{u}_0(\lambda)d\lambda=\int_{-1}^1 e^{i\lambda x-\omega(\lambda)t}\hat{u}_0(\lambda)d\lambda+\int_{\gamma^-+\gamma^+}e^{i\lambda x-\omega(\lambda)t}\sigma_M(\lambda)d\lambda+\int_{|\lambda|\ge1}e^{i\lambda x-\omega(\lambda)t}\frac{(u_0^{(M)})\hat{}(\lambda)}{(i\lambda)^M}d\lambda\,,\ \forall M>0,\quad(2.36)$$

which implies that the limit in (2.35) is equal to

$$\int_{-1}^1(i\lambda)^n[-\omega(\lambda)]^m e^{i\lambda x}\hat{u}_0(\lambda)d\lambda+\int_{\gamma^-+\gamma^+}(i\lambda)^n[-\omega(\lambda)]^m e^{i\lambda x}\sigma_M(\lambda)d\lambda+\int_{|\lambda|\ge1}(i\lambda)^n[-\omega(\lambda)]^m e^{i\lambda x}\frac{(u_0^{(M)})\hat{}(\lambda)}{(i\lambda)^M}d\lambda\,,$$

provided that we choose $M$ sufficiently large.

*Step 2* Proof of $1^{st}$ part Since the limit as $t\to0^+$ can be interchanged with the integrals of the RHS of (2.15), which are taken over $\Gamma_1$ and $\Gamma_2$ (because of (2.14)), the $1^{st}$ part of the theorem follows from step 1.

**Step 3** Using (2.5), Cauchy's theorem and Jordan's lemma, we see that

$$\int_{\Gamma_2} e^{i\lambda x - \omega(\lambda)t}\lambda^3 \widetilde{g}_0(\omega(\lambda),t)d\lambda = \int_{\Gamma_2 \cap \{|\lambda| \le 1\}} e^{i\lambda x - \omega(\lambda)t}\lambda^3 \widetilde{g}_0(\omega(\lambda),t)d\lambda + \int_{\{\frac{\pi}{8} \le \arg \lambda \le \frac{3\pi}{8}\} \cap \{|\lambda|=1\}} e^{i\lambda x}\frac{d\lambda}{\lambda}$$

$$+ \int_{\{|\lambda|=1\} \cap \{0 \le \arg \lambda \le \frac{\pi}{8}\} + \{|\lambda| \ge 1\} \cap \{\arg \lambda = 0\}} e^{i\lambda x - \omega(\lambda)t}\frac{d\lambda}{\lambda} + \int_{\{|\lambda|=1\} \cap \{\frac{3\pi}{8} \le \arg \lambda \le \frac{\pi}{2}\} + \{|\lambda| \ge 1\} \cap \{\arg \lambda = \frac{\pi}{2}\}} e^{i\lambda x - \omega(\lambda)t}\frac{d\lambda}{\lambda}$$

$$+ \int_{\Gamma_2 \cap \{|\lambda| \ge 1\}} e^{i\lambda x - \omega(\lambda)t}(g_j \, {}'\widetilde{)}(\omega(\lambda),t)\frac{d\lambda}{\lambda} \, . \quad (2.37)$$

It follows from (2.37) that the limits

$$\lim_{x \to 0^+} \frac{\partial^n}{\partial x^n}\left(\int_{\Gamma_2} e^{i\lambda x - \omega(\lambda)t}\lambda^3 \widetilde{g}_0(\omega(\lambda),t)d\lambda\right), \text{ for } n = 0, 1, 2, 3 \, ,$$

exist and define $C^\infty$ functions for $t \in (0, +\infty)$.

In order to deal with higher order derivatives (including derivatives with respect to $t$), we have to continue the process which led to (2.37), rewriting the last integral in (2.37) using (2.5), and so on.

**Step 4** *Proof of $2^{nd}$ part* Working as in the previous step, we can deal with all the integrals in the RHS of (2.12), taken on $\Gamma_1$ and $\Gamma_2$, as far as the $2^{nd}$ part is concerned. The first integral in the RHS of (2.12) is easier to deal with.

**Completion of the proof of theorem** Since the limits, which define the functions $U_{n,m}(x)$ and $V_{n,m}(t)$, are uniform in $x$ and $t$, respectively, as stated in $1^{st}$ and $2^{nd}$ part, the $3^{rd}$ assertion follows, also in view of Theorem 1, and, in particular, from the formulas (2.16) and (2.17). Finally, the assertions in $4^{th}$, $5^{th}$ and $6^{th}$ parts follow from $3^{rd}$ part in combination with equations (2.16), (2.17) and (2.18).

**Theorem 3** *For fixed $t_0 > 0$, the function $u(x,t)$, defined by (2.12), satisfies*

$$\lim_{x \to +\infty}\left(x^\ell \frac{\partial^{n+m}u(x,t)}{\partial x^n \partial t^m}\right) = 0, \quad (2.38)$$

*for every nonnegative integers $n$, $m$ and $\ell$, uniformly for $0 < t \le t_0$.*

**Proof** *Step 1* We claim that

$$\lim_{x \to +\infty} \frac{\partial^{n+m}}{\partial x^n \partial t^m}\left(\int_{-\infty}^{\infty} e^{i\lambda x - \omega(\lambda)t}\hat{u}_0(\lambda)d\lambda\right) = 0, \quad (2.39)$$

uniformly for $0 < t \le t_0$.

To prove this, let us write the following variation of (2.36):

$$\int_{-\infty}^{\infty} e^{i\lambda x - \omega(\lambda)t}\hat{u}_0(\lambda)d\lambda = \int_{-1}^{1} e^{i\lambda x - \omega(\lambda)t}\hat{u}_0(\lambda)d\lambda + \int_{\gamma_0^- + \gamma_0^+} e^{i\lambda x - \omega(\lambda)t}\sigma_{\mathrm{M}}(\lambda)d\lambda + \int_{|\lambda| \ge 1} e^{i\lambda x - \omega(\lambda)t}[\hat{u}_0(\lambda) - \sigma_{\mathrm{M}}(\lambda)]d\lambda \, , \quad (2.40)$$

where

$$\gamma_0^- := [\infty e^{i7\pi/8}, \frac{1}{\cos(\pi/8)}e^{i7\pi/8}] + [\frac{1}{\cos(\pi/8)}e^{i7\pi/8}, -1] \text{ and } \gamma_0^+ := [1, \frac{1}{\cos(\pi/8)}e^{i\pi/8}] + [\frac{1}{\cos(\pi/8)}e^{i\pi/8}, \infty e^{i\pi/8}] \, .$$

Working with $\mathrm{M}$ sufficiently large, we have

$$\frac{\partial^{n+m}}{\partial x^n \partial t^m}\left(\int_{-\infty}^{\infty} e^{i\lambda x - \omega(\lambda)t}\hat{u}_0(\lambda)d\lambda\right) = \int_{-1}^{1}(i\lambda)^n[-\omega(\lambda)]^m e^{i\lambda x - \omega(\lambda)t}\hat{u}_0(\lambda)d\lambda$$

$$+ \int_{\gamma_0^- + \gamma_0^+}(i\lambda)^n[-\omega(\lambda)]^m e^{i\lambda x - \omega(\lambda)t}\sigma_{\mathrm{M}}(\lambda)d\lambda + \int_{|\lambda| \ge 1}(i\lambda)^n[-\omega(\lambda)]^m e^{i\lambda x - \omega(\lambda)t}[\hat{u}_0(\lambda) - \sigma_{\mathrm{M}}(\lambda)]d\lambda$$

$$= \frac{1}{ix}\int_{-1}^{1}e^{i\lambda x}\frac{d}{d\lambda}\{(i\lambda)^n[-\omega(\lambda)]^m e^{-\omega(\lambda)t}\hat{u}_0(\lambda)\}d\lambda + \frac{1}{ix}\int_{\gamma_0^- + \gamma_0^+}e^{i\lambda x}\frac{d}{d\lambda}\{(i\lambda)^n[-\omega(\lambda)]^m e^{-\omega(\lambda)t}\sigma_{\mathrm{M}}(\lambda)\}d\lambda$$

$$+ \frac{1}{ix}\int_{|\lambda| \ge 1}e^{i\lambda x}\frac{d}{d\lambda}\{(i\lambda)^n[-\omega(\lambda)]^m e^{-\omega(\lambda)t}[\hat{u}_0(\lambda) - \sigma_{\mathrm{M}}(\lambda)]\}d\lambda \, . \quad (2.41)$$

The last equation in (2.41) follows by integration by parts and the fact that the intermediate boundary terms cancel each other.

Letting $x \to \infty$, we easily see that (2.41) implies (2.39).

**Step 2** We claim that

$$\lim_{x \to +\infty}\left[ x^{\ell} \frac{\partial^{n+m}}{\partial x^n \partial t^m} \left( \int_{-\infty}^{\infty} e^{i\lambda x - \omega(\lambda)t} \hat{u}_0(\lambda) d\lambda \right) \right] = 0\,, \tag{2.42}$$

uniformly for $0 < t \le t_0$.

Indeed, the case $\ell = 1$ of (2.42) follows by multiplying (2.41) by $x$ and integrating – twice – by parts as in step 1. The general case of any positive integer $\ell$ can be proved similarly by repeated integration by parts.

**Step 3** We claim that, for fixed $t_0 > 0$,

$$\lim_{x \to +\infty}\left[ x^{\ell} \frac{\partial^{n+m}}{\partial x^n \partial t^m} \left( \int_{\Gamma_1} e^{i\lambda x - \omega(\lambda)t} \lambda^3 \widetilde{g}_0(\omega(\lambda), t) d\lambda \right) \right] = 0\,, \tag{2.43}$$

uniformly for $0 < t \le t_0$.

To prove this, let us write the integral in (2.43) as follows:

$$\int_{\Gamma_1} \cdots d\lambda = \int_{\Gamma_1 \cap \{|\lambda| \ge 1\}} \cdots d\lambda + \int_{\Gamma_1 \cap \{|\lambda| \le 1\}} \cdots d\lambda\,. \tag{2.44}$$

In order to deal with the integral over $\Gamma_1 \cap \{|\lambda| \ge 1\}$, we use the inequalities

$$x^{\ell}\left|\lambda^N e^{i\lambda x}\right| \le x^{\ell}\left|\lambda\right|^N e^{-[x\sin(\pi/8)]|\lambda|} \text{ or } x^{\ell}\left|\lambda^N e^{i\lambda x}\right| \le x^{\ell}\left|\lambda\right|^N e^{-[x\sin(3\pi/8)]|\lambda|}\,, \text{ for } \lambda \in \Gamma_1 \text{ with } |\lambda| \ge 1\,,$$

and

$$\left|e^{-\omega(\lambda)t}\widetilde{g}_0(\omega(\lambda), t)\right| \le \int_0^t |g_0(\tau)|d\tau\,, \text{ for } \lambda \in \Gamma_1\,,$$

to obtain that, for $x \ge 1$,

$$x^{\ell}\left| \int_{\Gamma_1 \cap \{|\lambda| \ge 1\}} \lambda^N e^{i\lambda x - \omega(\lambda)t} \lambda^3 \widetilde{g}_0(\omega(\lambda), t) d\lambda \right| \le \left( \int_0^t |g_0(\tau)|d\tau \right) x^{\ell} e^{-[x\sin(\pi/8)]/2} \int_{\Gamma_1 \cap \{|\lambda| \ge 1\}} |\lambda|^{N+3} e^{-[\sin(\pi/8)]|\lambda|/2} d|\lambda|\,.$$

In order to deal with the integral over $\Gamma_1 \cap \{|\lambda| \le 1\}$, it suffices to integrate by parts several times, depending on $\ell$.

**Step 4** *Completion of the proof* First, we prove analogues of (2.43) for all the integrals in the RHS of (2.12) which are taken on $\Gamma_1$ or $\Gamma_2$. Then combining these with (2.39), we obtain (2.38).

# 3. The equation $\partial_t U = -\partial_x^4 U + f$

**Problem** *Solve* (1.1) *in the case of the equation* $\partial_t U = -\partial_x^4 U + f$, *i.e., when* $\alpha = 0$ *and* $\beta = 1$. $\tag{3.1}$

### 3.1 *The Fokas method solution* Defining

$$\hat{f}(\lambda, t) = \int_{y=0}^{\infty} e^{-i\lambda y} f(y, t) dy \text{ and } \widetilde{f}(\lambda, \omega(\lambda), t) = \int_{\tau=0}^{t} e^{\omega(\lambda)\tau} \hat{f}(\lambda, \tau) d\tau \text{ (where } \omega(\lambda) = \lambda^4),$$

for $\lambda \in \mathbb{C}$ with $\operatorname{Im}\lambda \le 0$, and setting

$$\Phi_{\mathbb{R}}(x, t) = \int_{\lambda=-\infty}^{\infty} e^{i\lambda x - \omega(\lambda)t} \widetilde{f}(\lambda, \omega(\lambda), t) d\lambda$$

$$\Phi_{\Gamma_1}(x, t) = \int_{\Gamma_1} e^{i\lambda x - \omega(\lambda)t} [-i\widetilde{f}(-\lambda, \omega(\lambda), t) + (i-1)\widetilde{f}(i\lambda, \omega(\lambda), t)]d\lambda$$

and

$$\Phi_{\Gamma_2}(x,t) = \int_{\Gamma_2} e^{i\lambda x - \omega(\lambda)t} [i\widetilde{\widehat{f}}(-\lambda, \omega(\lambda), t) - (i+1)\widetilde{\widehat{f}}(-i\lambda, \omega(\lambda), t)]d\lambda ,$$

for $x > 0$ and $t > 0$, the solution $U(x,t)$ of (3.1) has the following integral representation:

$$U(x,t) = u(x,t) + \frac{1}{2\pi}\Phi_{\mathbb{R}}(x,t) + \frac{1}{2\pi}\Phi_{\Gamma_1}(x,t) + \frac{1}{2\pi}\Phi_{\Gamma_2}(x,t), \tag{3.2}$$

where $u(x,t)$ is the solution of the homogeneous equation given by (2.12).

***Some computations*** **(1)** Integrating by parts we obtain

$$\hat{f}(\lambda,t) = \int_{y=0}^{\infty} e^{-i\lambda y}f(y,t)dy = \frac{1}{i\lambda}f(0,t) + \frac{1}{i\lambda}\int_{y=0}^{\infty}e^{-i\lambda y}\frac{\partial f(y,t)}{\partial y}dy \quad (\lambda \in \mathbb{C} - \{0\}, \operatorname{Im}\lambda \le 0, \ t \ge 0) \tag{3.4}$$

and

$$\widetilde{\hat{f}}(\lambda, \omega(\lambda), t) = \int_{\tau=0}^{t} e^{\omega(\lambda)\tau}\hat{f}(\lambda,\tau)d\tau$$

$$= \frac{1}{i\lambda}\int_{\tau=0}^{t}\frac{\partial}{\partial\tau}\left(\frac{e^{\omega(\lambda)\tau}}{\omega(\lambda)}\right)f(0,\tau)d\tau + \frac{1}{i\lambda}\int_{\tau=0}^{t}\frac{\partial}{\partial\tau}\left(\frac{e^{\omega(\lambda)\tau}}{\omega(\lambda)}\right)\left(\int_{y=0}^{\infty}e^{-i\lambda y}\frac{\partial f(y,\tau)}{\partial y}dy\right)d\tau$$

$$= \frac{1}{i\lambda}\frac{e^{\omega(\lambda)t}}{\omega(\lambda)}f(0,t) - \frac{1}{i\lambda}\frac{1}{\omega(\lambda)}f(0,0) - \frac{1}{i\lambda}\frac{1}{\omega(\lambda)}\int_{\tau=0}^{t}e^{\omega(\lambda)\tau}\frac{\partial}{\partial\tau}[f(0,\tau)]d\tau$$

$$+ \frac{1}{i\lambda}\frac{e^{\omega(\lambda)t}}{\omega(\lambda)}\left(\int_{y=0}^{\infty}e^{-i\lambda y}\frac{\partial f(y,t)}{\partial y}dy\right) - \frac{1}{i\lambda}\frac{1}{\omega(\lambda)}\left(\int_{y=0}^{\infty}e^{-i\lambda y}\frac{\partial f(y,0)}{\partial y}dy\right)$$

$$- \frac{1}{i\lambda}\frac{1}{\omega(\lambda)}\int_{\tau=0}^{t}e^{\omega(\lambda)\tau}\left(\int_{y=0}^{\infty}e^{-i\lambda y}\frac{\partial^2 f(y,t)}{\partial\tau\partial y}dy\right)d\tau . \tag{3.5}$$

**(2)** It follows from (3.5), in view of (3.2), that

$$e^{-\omega(\lambda)t}\widetilde{\hat{f}}(\lambda, \omega(\lambda), t) = \mathrm{O}(1/\lambda^5) \ \text{ for } \lambda \to \infty, \text{ with } \lambda \in \mathbb{C} \text{ and } \operatorname{Re}\omega(\lambda) = \operatorname{Re}(\lambda^4) \ge 0 , \tag{3.6}$$

in the sense that

$$\sup\left\{\left|\lambda^5 \cdot e^{-\omega(\lambda)t}\widetilde{\hat{f}}(\lambda, \omega(\lambda), t)\right| : |\lambda| \ge 1, \operatorname{Re}\omega(\lambda) \ge 0, 0 \le t \le t_0\right\} < +\infty .$$

**(3)** Generalizing (3.4) and similarly to equation (1.4), integration by parts gives

$$\hat{f}(\lambda,t) = h_{\mathrm{N}}(\lambda,t) + \frac{1}{(i\lambda)^{\mathrm{N}}}\int_{y=0}^{\infty}e^{-i\lambda y}\frac{\partial^{\mathrm{N}}f(y,t)}{\partial y^{\mathrm{N}}}dy \text{ where } h_{\mathrm{N}}(\lambda,t) := \sum_{m=1}^{\mathrm{N}}\frac{1}{(i\lambda)^m}\frac{\partial^{m-1}f(y,t)}{\partial y^{m-1}}\bigg|_{y=0} .$$

Therefore,

$$\delta_{\mathrm{N}}(\lambda,t) := \hat{f}(\lambda,t) - h_{\mathrm{N}}(\lambda,t) = \mathrm{O}(1/\lambda^{\mathrm{N}+1}), \text{ for } \lambda \to \infty, \text{ with } \operatorname{Im}\lambda \le 0 . \tag{3.7}$$

Now

$$\widetilde{\hat{f}}(\lambda, \omega(\lambda), t) = \int_{\tau=0}^{t} e^{\omega(\lambda)\tau}\hat{f}(\lambda,\tau)d\tau$$

$$= \int_{\tau=0}^{t} e^{\omega(\lambda)\tau}\delta_{\mathrm{N}}(\lambda,\tau)d\tau + \int_{\tau=0}^{t} e^{\omega(\lambda)\tau}h_{\mathrm{N}}(\lambda,\tau)d\tau = \widetilde{\delta}_{\mathrm{N}}(\lambda, \omega(\lambda), t) + \widetilde{h}_{\mathrm{N}}(\lambda, \omega(\lambda), t) . \tag{3.8}$$

Let us notice that, although (3.8) holds for $\lambda \in \mathbb{C} - \{0\}$ with $\operatorname{Im}\lambda \le 0$, the last integral in (3.8), i.e. $\widetilde{h}_{\mathrm{N}}(\lambda, \omega(\lambda), t)$, is defined for every $\lambda \in \mathbb{C} - \{0\}$.

Further integration by parts gives

$$\widetilde{h}_{\mathrm{N}}(\lambda, \omega(\lambda), t) = e^{\omega(\lambda)t}\mu_{\mathrm{N},\mathrm{M}}(\lambda,t) - \mu_{\mathrm{N},\mathrm{M}}(\lambda,0) + \frac{(-1)^{\mathrm{M}}}{[\omega(\lambda)]^{\mathrm{M}}}\int_{\tau=0}^{t}e^{\omega(\lambda)\tau}h_{\mathrm{N}}^{(\mathrm{M})}(\lambda,\tau)d\tau \tag{3.9}$$

where

$$\mu_{\mathrm{N},\mathrm{M}}(\lambda,t) := \frac{h_{\mathrm{N}}(\lambda,t)}{\omega(\lambda)} - \frac{h_{\mathrm{N}}{}'(\lambda,t)}{\omega^2(\lambda)} + \cdots + (-1)^{\mathrm{M}-1}\frac{h_{\mathrm{N}}^{(\mathrm{M}-1)}(\lambda,t)}{\omega^{\mathrm{M}}(\lambda)} \quad (\lambda \in \mathbb{C} - \{0\}) .$$

(The derivatives of $h_N(\lambda, t)$, in the above quantity, are taken with respect to $t$.)

**(4)** Similarly to (3.6), we have, for $\lambda \to \infty$ with $\lambda \in \mathbb{C}$ and $\operatorname{Re} \omega(\lambda) = \operatorname{Re}(\lambda^4) \geq 0$,

$$e^{-\omega(\lambda)t}\widetilde{h}_N(\lambda, \omega(\lambda), t) = \mathrm{O}(1/\lambda^5), \ e^{-\omega(\lambda)t}\widetilde{\delta}_N(\lambda, \omega(\lambda), t) = \mathrm{O}(1/\lambda^{N+5}) \tag{3.10}$$

and

$$\frac{1}{[\omega(\lambda)]^M} e^{-\omega(\lambda)t} \int_{\tau=0}^{t} e^{\omega(\lambda)\tau} h_N^{(M)}(\lambda, \tau) d\tau = \mathrm{O}(1/\lambda^{4M+5}) . \tag{3.11}$$

**Theorem 4** *The function $U(x,t)$, defined by (3.2), is $C^\infty$ for $(x,t) \in Q$ and satisfies the following:*

$1^{st}$ *The differential equation $\partial_t U = -\partial_x^4 U + f$, for $x > 0$ and $t > 0$.*

$2^{nd}$ *The limit $\lim_{t\to 0^+} U(x,t) = u_0(x)$, uniformly for $x \geq x_0$ ($\forall x_0 > 0$).*

$3^{rd}$ *The limits $\lim_{x\to 0^+} U(x,t) = g_0(t)$ and $\lim_{x\to 0^+} \frac{\partial U(x,t)}{\partial x} = g_0(t)$, uniformly for $t \geq t_0$ ($\forall t_0 > 0$).*

$4^{th}$ *The limit $\lim_{t\to 0^+} \frac{\partial^n U(x,t)}{\partial x^n} = \frac{d^n u_0(x)}{dx^n}$, uniformly for $x \geq x_0$.*

$5^{th}$ *The limits $\lim_{x\to 0^+} \frac{\partial^m U(x,t)}{\partial t^m} = \frac{d^m g_0(t)}{dt^m}$ and $\lim_{x\to 0^+} \frac{\partial^{m+1} U(x,t)}{\partial t^m \partial x} = \frac{d^m g_1(t)}{dt^m}$, uniformly for $t \geq t_0$.*

$6^{th}$ *The limit $\lim_{x\to +\infty}\left( x^\ell \frac{\partial^{n+m} U(x,t)}{\partial x^n \partial t^m}\right) = 0$, for every nonnegative integers $n$, $m$ and $\ell$, uniformly for $0 < t \leq t_0$.*

**Proof** *Step 1* First, the function $U(x,t)$ is $C^\infty$ for $(x,t) \in Q$, since $u(x,t)$ is $C^\infty$ and so are the functions defined by the integrals $\Phi_{\mathbb{R}}(x,t)$ and $\Phi_{\Gamma_1}(x,t) + \Phi_{\Gamma_2}(x,t)$, because of the presence of the factors $e^{-\omega(\lambda)t}$ and $e^{i\lambda x}$, respectively. Now,

$$\left(\frac{\partial}{\partial t} + \frac{\partial^4}{\partial x^4}\right)[\Phi_{\mathbb{R}}(x,t) + \Phi_{\Gamma_1}(x,t) + \Phi_{\Gamma_2}(x,t)] = \int_{-\infty}^{\infty} e^{i\lambda x}\hat{f}(\lambda, t)d\lambda$$

$$+ \int_{\Gamma_1} e^{i\lambda x}[-i\hat{f}(-\lambda,t) + (i-1)\hat{f}(i\lambda,t)]d\lambda + \int_{\Gamma_2} e^{i\lambda x}[i\hat{f}(-\lambda,t) - (i+1)\hat{f}(-i\lambda,t)]d\lambda = 2\pi f(x) ,$$

where the last equation follows from Fourier's inversion formula, Cauchy's theorem and Jordan's lemma.

*Step 2* It follows from (3.6) that
$$\lim_{t\to 0^+} \Phi_{\mathbb{R}}(x,t) = 0, \ \lim_{t\to 0^+} \Phi_{\Gamma_1}(x,t) = 0, \ \lim_{t\to 0^+} \Phi_{\Gamma_2}(x,t) = 0 ,$$

and

$$\lim_{x\to 0^+}[\Phi_{\mathbb{R}}(x,t) + \Phi_{\Gamma_1}(x,t) + \Phi_{\Gamma_2}(x,t)] = \int_{-\infty}^{\infty} e^{-\omega(\lambda)t}\widetilde{\widehat{f}}(\lambda, \omega(\lambda), t)d\lambda$$

$$+ \int_{\Gamma_1} e^{-\omega(\lambda)t}[-i\widehat{\widetilde{f}}(-\lambda,\omega(\lambda),t) + (i-1)\widehat{\widetilde{f}}(i\lambda,\omega(\lambda),t)]d\lambda + \int_{\Gamma_2} e^{-\omega(\lambda)t}[i\widehat{\widetilde{f}}(-\lambda,\omega(\lambda),t) - (i+1)\widehat{\widetilde{f}}(-i\lambda,\omega(\lambda),t)]d\lambda = 0 ,$$

where the last equation follows as (2.21).
Similarly, also in view of (3.6),

$$\lim_{x\to 0^+}\left\{\frac{\partial}{\partial x}[\Phi_{\mathbb{R}}(x,t) + \Phi_{\Gamma_1}(x,t) + \Phi_{\Gamma_2}(x,t)]\right\} = \int_{\lambda=-\infty}^{\infty} e^{-\omega(\lambda)t}(i\lambda)\widehat{\widetilde{f}}(\lambda, \omega(\lambda), t)d\lambda$$

$$+ \int_{\Gamma_1} e^{-\omega(\lambda)t}(i\lambda)[-i\widehat{\widetilde{f}}(-\lambda,\omega(\lambda),t) + (i-1)\widehat{\widetilde{f}}(i\lambda,\omega(\lambda),t)]d\lambda + \int_{\Gamma_2} e^{-\omega(\lambda)t}(i\lambda)[i\widehat{\widetilde{f}}(-\lambda,\omega(\lambda),t) - (i+1)\widehat{\widetilde{f}}(-i\lambda,\omega(\lambda),t)]d\lambda = 0 ,$$

as (2.31).
The above results, combined with the conclusions of Theorem 1, imply the $1^{st}$, $2^{nd}$ and $3^{rd}$ assertions.

*Step 3* We claim that the limit

$$\lim_{t \to 0^+} \frac{\partial^{n+m} U(x,t)}{\partial x^n \partial t^m} \tag{3.12}$$

exists, is uniform for $x \geq x_0$ ($\forall x_0 > 0$) and defines a $C^\infty$ function of $x \in (0, +\infty)$.

To prove this, let us write first, in view of (3.8), the equation

$$\Phi_{\mathbb{R}}(x,t) = \int_{-1}^{1} e^{i\lambda x - \omega(\lambda)t} \widetilde{\widetilde{f}}(\lambda, \omega(\lambda), t) d\lambda + \left( \int_{-\infty}^{-1} + \int_{1}^{\infty} \right) e^{i\lambda x - \omega(\lambda)t} \widetilde{\delta}_{\mathrm{N}}(\lambda, \omega(\lambda), t) d\lambda + \int_{\gamma^- + \gamma^+} e^{i\lambda x - \omega(\lambda)t} \widetilde{h}_{\mathrm{N}}(\lambda, \omega(\lambda), t) d\lambda . \tag{3.13}$$

It follows that

$$\frac{\partial^{n+m} \Phi_{\mathbb{R}}(x,t)}{\partial x^n \partial t^m} = \int_{-1}^{1} (i\lambda)^n [-\omega(\lambda)]^m e^{i\lambda x - \omega(\lambda)t} \widetilde{\widetilde{f}}(\lambda, \omega(\lambda), t) d\lambda$$

$$+ \left( \int_{-\infty}^{-1} + \int_{1}^{\infty} \right) (i\lambda)^n [-\omega(\lambda)]^m e^{i\lambda x - \omega(\lambda)t} \widetilde{\delta}_{\mathrm{N}}(\lambda, \omega(\lambda), t) d\lambda + \int_{\gamma^- + \gamma^+} (i\lambda)^n [-\omega(\lambda)]^m e^{i\lambda x - \omega(\lambda)t} \widetilde{h}_{\mathrm{N}}(\lambda, \omega(\lambda), t) d\lambda . \tag{3.14}$$

Choosing N sufficiently large and letting $t \to 0^+$, in (3.14), we easily obtain the existence of the limit (3.12). (The integrals $\Phi_{\Gamma_1}(x,t)$ and $\Phi_{\Gamma_2}(x,t)$, which are parts of $U(x,t)$, are easy to deal with, as far as the limit $\lim_{t \to 0^+}$ is conserned, because of the presence of the factor $e^{i\lambda x}$ and (2.14).)

***Step 4*** We claim that the limit

$$\lim_{x \to 0^+} \frac{\partial^{n+m} U(x,t)}{\partial x^n \partial t^m}$$

exists, is uniform for $t \geq t_0$ ($\forall t_0 > 0$) and defines a $C^\infty$ function of $t \in (0, +\infty)$.

To prove this, first, we substitute (3.9) in (3.13) and we obtain

$$\Phi_{\mathbb{R}}(x,t) = \int_{-1}^{1} e^{i\lambda x - \omega(\lambda)t} \widetilde{\widetilde{f}}(\lambda, \omega(\lambda), t) d\lambda + \left( \int_{-\infty}^{-1} + \int_{1}^{\infty} \right) e^{i\lambda x - \omega(\lambda)t} \widetilde{\delta}_{\mathrm{N}}(\lambda, \omega(\lambda), t) d\lambda$$

$$+ \int_{\gamma^- + \gamma^+} e^{i\lambda x} \mu_{\mathrm{N,M}}(\lambda, t) d\lambda - \left( \int_{-\infty}^{-1} + \int_{1}^{\infty} \right) e^{i\lambda x - \omega(\lambda)t} \mu_{\mathrm{N,M}}(\lambda, 0) d\lambda + \left( \int_{-\infty}^{-1} + \int_{1}^{\infty} \right) \left[ e^{i\lambda x - \omega(\lambda)t} \frac{(-1)^{\mathrm{M}}}{[\omega(\lambda)]^{\mathrm{M}}} \int_{\tau=0}^{t} e^{\omega(\lambda)\tau} h_{\mathrm{N}}^{(\mathrm{M})}(\lambda, \tau) d\tau \right] d\lambda$$

$$= \int_{-1}^{1} e^{i\lambda x - \omega(\lambda)t} \widetilde{\widetilde{f}}(\lambda, \omega(\lambda), t) d\lambda + \left( \int_{-\infty}^{-1} + \int_{1}^{\infty} \right) e^{i\lambda x - \omega(\lambda)t} \widetilde{\delta}_{\mathrm{N}}(\lambda, \omega(\lambda), t) d\lambda$$

$$+ \int_{\{|\lambda|=1 \& 0 \leq \arg \lambda \leq \pi\}} e^{i\lambda x} \mu_{\mathrm{N,M}}(\lambda, t) d\lambda - \left( \int_{-\infty}^{-1} + \int_{1}^{\infty} \right) e^{i\lambda x - \omega(\lambda)t} \mu_{\mathrm{N,M}}(\lambda, 0) d\lambda$$

$$\left( \int_{-\infty}^{-1} + \int_{1}^{\infty} \right) \left[ e^{i\lambda x - \omega(\lambda)t} \frac{(-1)^{\mathrm{M}}}{[\omega(\lambda)]^{\mathrm{M}}} \int_{\tau=0}^{t} e^{\omega(\lambda)\tau} h_{\mathrm{N}}^{(\mathrm{M})}(\lambda, \tau) d\tau \right] d\lambda . \tag{3.15}$$

Choosing N and M sufficiently large, differentiating (3.15) and letting $x \to 0^+$, we easily obtain the required conclusion concerning the limit, as $x \to 0^+$, of the derivatives of $\Phi_{\mathbb{R}}(x,t)$.

To deal with the part $\Phi_{\Gamma_1}(x,t) + \Phi_{\Gamma_2}(x,t)$ of $U(x,t)$, let us consider the integral

$$\Psi(x,t) := \int_{\Gamma_2} e^{i\lambda x - \omega(\lambda)t} \widetilde{\widetilde{f}}(-\lambda, \omega(\lambda), t) d\lambda ,$$

and write it as follows:

$$\Psi(x,t) = \int_{\Gamma_2 \cap \{|\lambda| \leq 1\}} e^{i\lambda x - \omega(\lambda)t} \widetilde{\widetilde{f}}(-\lambda, \omega(\lambda), t) d\lambda + \int_{\Gamma_2 \cap \{|\lambda| \geq 1\}} e^{i\lambda x - \omega(\lambda)t} \widetilde{\widetilde{f}}(-\lambda, \omega(\lambda), t) d\lambda$$

$$= \int_{\Gamma_2 \cap \{|\lambda| \leq 1\}} e^{i\lambda x - \omega(\lambda)t} \widetilde{\widetilde{f}}(-\lambda, \omega(\lambda), t) d\lambda + \int_{\Gamma_2 \cap \{|\lambda| \geq 1\}} e^{i\lambda x - \omega(\lambda)t} \widetilde{\delta}_{\mathrm{N}}(-\lambda, \omega(\lambda), t) d\lambda + \int_{\{|\lambda|=1 \& \frac{\pi}{8} \leq \arg \lambda \leq \frac{3\pi}{8}\}} e^{i\lambda x} \mu_{\mathrm{N,M}}(-\lambda, t) d\lambda$$

$$- \left( \int_{\{|\lambda|=1 \& 0 \leq \arg \lambda \leq \frac{\pi}{8}\}} + \int_{[1, +\infty)} \right) e^{i\lambda x - \omega(\lambda)t} \mu_{\mathrm{N,M}}(-\lambda, 0) d\lambda - \left( \int_{\{|\lambda|=1 \& \frac{3\pi}{8} \leq \arg \lambda \leq \frac{\pi}{2}\}} + \int_{[i, \infty i)} \right) e^{i\lambda x - \omega(\lambda)t} \mu_{\mathrm{N,M}}(-\lambda, 0) d\lambda$$

$$+ \int_{\Gamma_2 \cap \{|\lambda| \geq 1\}} \left[ e^{i\lambda x - \omega(\lambda)t} \frac{(-1)^{\mathrm{M}}}{[\omega(\lambda)]^{\mathrm{M}}} \int_{\tau=0}^{t} e^{\omega(\lambda)\tau} h_{\mathrm{N}}^{(\mathrm{M})}(-\lambda, \tau) d\tau \right] d\lambda . \quad (3.16)$$

Choosing N and M sufficiently large, differentiating (3.16) and letting $x \to 0^+$, we easily obtain the required conclusion concerning the limit, as $x \to 0^+$, of the derivatives of $\Psi(x,t)$.

Dealing with the other parts of $\Phi_{\Gamma_1}(x,t) + \Phi_{\Gamma_2}(x,t)$ is similar. This concludes the proof of the claim.

**Step 5** *Completion of the proof of the theorem* The 1$^{\mathrm{st}}$, 2$^{\mathrm{nd}}$ and 3$^{\mathrm{rd}}$ assertion follow from the calculations in steps 1 and 2, and the corresponding results of Theorem 1. These, combined with the conclusions of steps 3 and 4, imply that the function $U(x,t)$, extended to $\overline{Q} - \{(0,0)\}$ by setting

$$U*(x,0) := u_0(x) \text{ for } x > 0 \text{ and } U*(0,t) := g_0(t) \text{ for } t > 0 ,$$

is $C^\infty$ (in $\overline{Q} - \{(0,0)\}$). This, in turn, easily gives the 4$^{\mathrm{th}}$ and 5$^{\mathrm{th}}$ conclusions of the theorem.

Finally, the proof of 6$^{\mathrm{th}}$ conclusion is similar to the proof of Theorem 3.

## 4. The equation $\partial_t u = \partial_x^2 u - \partial_x^4 u$

**Problem** *Solve* (1.1) *in the case of the equation* $\partial_t u = \partial_x^2 u - \partial_x^4 u$ *, i.e., when* $\alpha = 1$, $\beta = 1$ *and* $f \equiv 0$. $\quad (4.1)$

**4.1 Jordan type lemmas in the case** $\omega(\lambda) = \lambda^2 + \lambda^4$. Let $t > 0$ and $\psi(\lambda)$ be a continuous function such that

$$\lim_{A \to \infty} [\sup \{|\psi(\lambda)| : |\lambda| = A\}] = 0 .$$

Let us also consider, for each $A > 0$, an arc $\mathcal{T}_A$ on the circle centered at $0$, of radius $A$. Then

$$\lim_{A \to \infty} \int_{\mathcal{T}_A \cap \{\lambda : \operatorname{Re}\omega(\lambda) \geq 0\}} \lambda^{3-j} e^{-\omega(\lambda)t} \psi(\lambda) d\lambda = 0 , \text{ for } j = 0, 1, 2, \dots .$$

For example, if $\Omega_k^+$, $k = 1, 2, 3, 4$, are the four connected components of the set

$$\Omega^+ := \{\lambda \in \mathbb{C} : \operatorname{Re}\omega(\lambda) = \operatorname{Re}(\lambda^2 + \lambda^4) \geq 0\} ,$$

then

$$\lim_{A \to \infty} \int_{\Omega_k^+ \cap \{\lambda : |\lambda| = A\}} \lambda^{3-j} e^{-\omega(\lambda)t} \psi(\lambda) d\lambda = 0 \text{ and } \lim_{A \to \infty} \int_{\mathcal{T}_{0,A}} \lambda^{3-j} e^{\omega(\lambda)t} \psi(\lambda) d\lambda = 0 \ ( j = 0, 1, 2, \dots ),$$

where $\mathcal{T}_{0,A}$ one of the arcs $\{\lambda : \pm \operatorname{Re}\lambda \leq 0, \operatorname{Im}\lambda \geq 0, \operatorname{Re}\omega(\lambda) \leq 0, |\lambda| = A \}$.

**4.2 Derivation of the solution** In this section, we set $\omega(\lambda) := \lambda^2 + \lambda^4$, so that $-\omega(\lambda)\big|_{\lambda = -i\partial_x} = \partial_x^2 u - \partial_x^4 u$. The identity

$$i \left[ \frac{(\lambda^2 + \lambda^4) - (\vartheta^2 + \vartheta^4)}{\lambda - \vartheta} \right]\Bigg|_{\vartheta = -i\partial_x} u = i[\lambda + \vartheta + \vartheta^3 + \vartheta^2\lambda + \vartheta\lambda^2 + \lambda^3]\big|_{\vartheta = -i\partial_x} u = -u_{xxx} - i\lambda u_{xx} + (\lambda^2 + 1)u_x + i\lambda(\lambda^2 + 1)u ,$$

leads to the equation

$$\frac{\partial}{\partial t}[e^{-i\lambda x + \omega(\lambda)t} u(x,t)] - \frac{\partial}{\partial x}\{e^{-i\lambda x + \omega(\lambda)t}[-u_{xxx} - i\lambda u_{xx} + (\lambda^2 + 1)u_x + i\lambda(\lambda^2 + 1)u]\} = e^{-i\lambda x + \omega(\lambda)t}\left( \frac{\partial u}{\partial t} - \frac{\partial^2 u}{\partial x^2} + \frac{\partial^4 u}{\partial x^4} \right).$$

Thus, if the function $u = u(x,t)$ solves (4.1) then

$$\hat{u}_0(\lambda) - \hat{u}(\lambda,t)e^{\omega(\lambda)t} + \widetilde{g}_3(\omega(\lambda),t) + i\lambda\widetilde{g}_2(\omega(\lambda),t) - (\lambda^2 + 1)\widetilde{g}_1(\omega(\lambda),t) - i\lambda(\lambda^2 + 1)\widetilde{g}_0(\omega(\lambda),t) = 0 , \ \operatorname{Im}\lambda \leq 0 . \quad (4.2)$$

Multiplying (4.2) by $e^{i\lambda x - \omega(\lambda)t}$ and integrating we obtain

$$\int_{-\infty}^{\infty} e^{i\lambda x - \omega(\lambda)t} \hat{u}_0(\lambda) d\lambda - \int_{-\infty}^{\infty} e^{i\lambda x} \hat{u}(\lambda,t) d\lambda - \int_{-\infty}^{\infty} e^{i\lambda x - \omega(\lambda)t} (\lambda^2 + 1)\widetilde{g}_1(\omega(\lambda),t) d\lambda - i\int_{-\infty}^{\infty} e^{i\lambda x - \omega(\lambda)t} \lambda(\lambda^2 + 1)\widetilde{g}_0(\omega(\lambda),t) d\lambda$$
$$+ \int_{-\infty}^{\infty} e^{i\lambda x - \omega(\lambda)t} \widetilde{g}_3(\omega(\lambda),t) d\lambda + i\int_{-\infty}^{\infty} e^{i\lambda x - \omega(\lambda)t} \lambda\widetilde{g}_2(\omega(\lambda),t) d\lambda = 0 .$$

Therefore, for $x > 0$ and $t > 0$,

$$2\pi u(x,t) = \int_{-\infty}^{\infty} e^{i\lambda x - \omega(\lambda)t}\hat{u}_0(\lambda)d\lambda - i\int_{-\infty}^{\infty} e^{i\lambda x - \omega(\lambda)t}\lambda(\lambda^2+1)\widetilde{g}_0(\omega(\lambda),t)d\lambda - \int_{-\infty}^{\infty} e^{i\lambda x - \omega(\lambda)t}(\lambda^2+1)\widetilde{g}_1(\omega(\lambda),t)d\lambda$$

$$+ i\int_{-\infty}^{\infty} e^{i\lambda x - \omega(\lambda)t}\lambda\widetilde{g}_2(\omega(\lambda),t)d\lambda + \int_{-\infty}^{\infty} e^{i\lambda x - \omega(\lambda)t}\widetilde{g}_3(\omega(\lambda),t)d\lambda \ . \ (4.3)$$

Next, we consider the set

$$\Omega := \{\lambda \in \mathbb{C} : \operatorname{Im}\lambda \ge 0 \ \text{ and } \operatorname{Re}\omega(\lambda) = \operatorname{Re}(\lambda^2+\lambda^4) \le 0\}\ .$$

With $\lambda = \xi + i\eta$, in the $\xi\eta$-plane, we have

$$\Omega = \{(\xi,\eta)\in\mathbb{R}^2 : \eta \ge 0 \ \text{ and } \eta^4 - (6\xi^2+1)\eta^2 + (\xi^2+\xi^4) \le 0\}\ .$$

By Cauchy's theorem and Jordan's lemma, (4.3) becomes

$$2\pi u(x,t) = \int_{-\infty}^{\infty} e^{i\lambda x - \omega(\lambda)t}\hat{u}_0(\lambda)d\lambda - i\int_{\partial\Omega} e^{i\lambda x - \omega(\lambda)t}\lambda(\lambda^2+1)\widetilde{g}_0(\omega(\lambda),t)d\lambda - \int_{\partial\Omega} e^{i\lambda x - \omega(\lambda)t}(\lambda^2+1)\widetilde{g}_1(\omega(\lambda),t)d\lambda$$

$$+ i\int_{\partial\Omega} e^{i\lambda x - \omega(\lambda)t}\lambda\widetilde{g}_2(\omega(\lambda),t)d\lambda + \int_{\partial\Omega} e^{i\lambda x - \omega(\lambda)t}\widetilde{g}_3(\omega(\lambda),t)d\lambda \ , \ (4.4)$$

where $\partial\Omega = \Gamma_1 + \Gamma_2 + \Gamma_3 + \Gamma_4$ (See fig. 4 – The notation in this section is different from the notation in sections 2 and 3, as far as the contours $\Gamma$'s are concerned.)

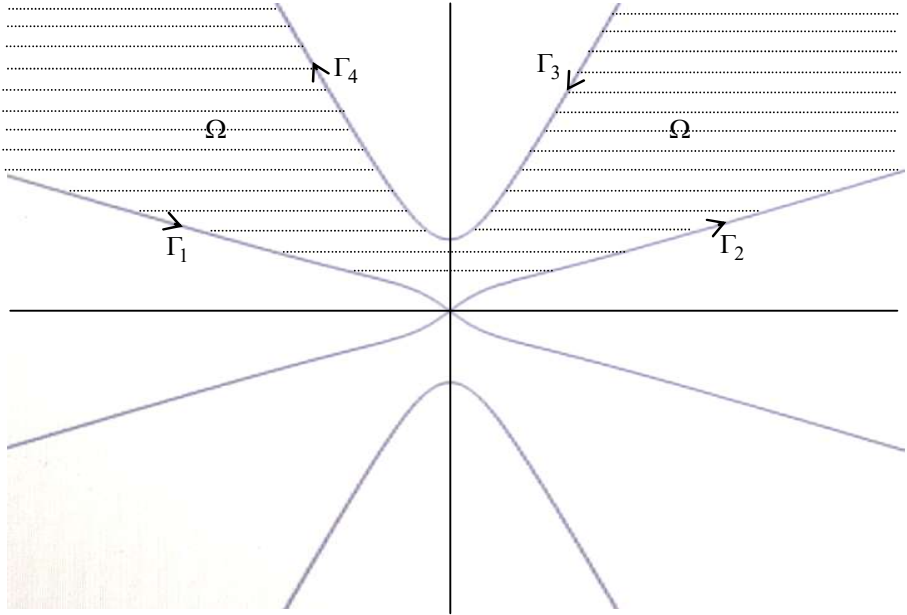

**Fig. 4** $\Omega$ is the shaded set and $\partial\Omega = \Gamma_1 + \Gamma_2 + \Gamma_3 + \Gamma_4$.

Next, we observe that if $\rho(\lambda) = \sqrt{\lambda^2+1}$ then

$$\omega(-\lambda) = \omega(i\rho(\lambda)) = \omega(-i\rho(\lambda)) = \omega(\lambda)\ , \ \forall\lambda\ . \tag{4.5}$$

We choose the following holomorphic branch of the square root of $\lambda^2+1$, $\rho:\Theta\to\mathbb{C}$, in the set

$$\Theta := \mathbb{C} - \{\lambda = \xi + i\eta \in \mathbb{C} : \xi = 0 \ \& \ |\eta| \ge 1\} = \mathbb{C} - [(-i\infty, -i]\cup[i,i\infty)]\ ,$$

defined by

$$\rho(\lambda) = \sqrt{|\lambda - i|\cdot|\lambda + i|}\exp\{i[\theta_i(\lambda) + \theta_{-i}(\lambda) + \pi]/2\}\ , \text{ for } \lambda\in\Theta \text{ (see fig.5)}.$$

($\sqrt{|\lambda - i|\cdot|\lambda + i|}$ denotes the positive square root of the positive number $|\lambda - i|\cdot|\lambda + i|$, for $\lambda \ne \pm i$.)

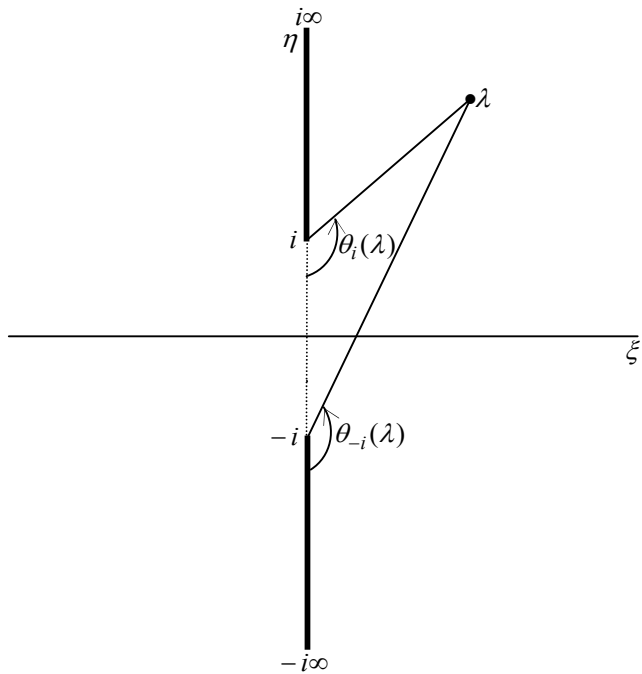

**Fig.5** With $-\pi < \theta_i(\lambda) < \pi$ and $0 < \theta_{-i}(\lambda) < 2\pi$, the function $\rho(\lambda)$ is analytic for $\lambda \in \Theta$ and $\rho(-\lambda) = \rho(\lambda)$.

With this choice of the square root, we have

$$\lambda \in \partial\Omega \implies \text{Im}[i\rho(\lambda)] \le 0 . \text{ (See fig6.)} \tag{4.6}$$

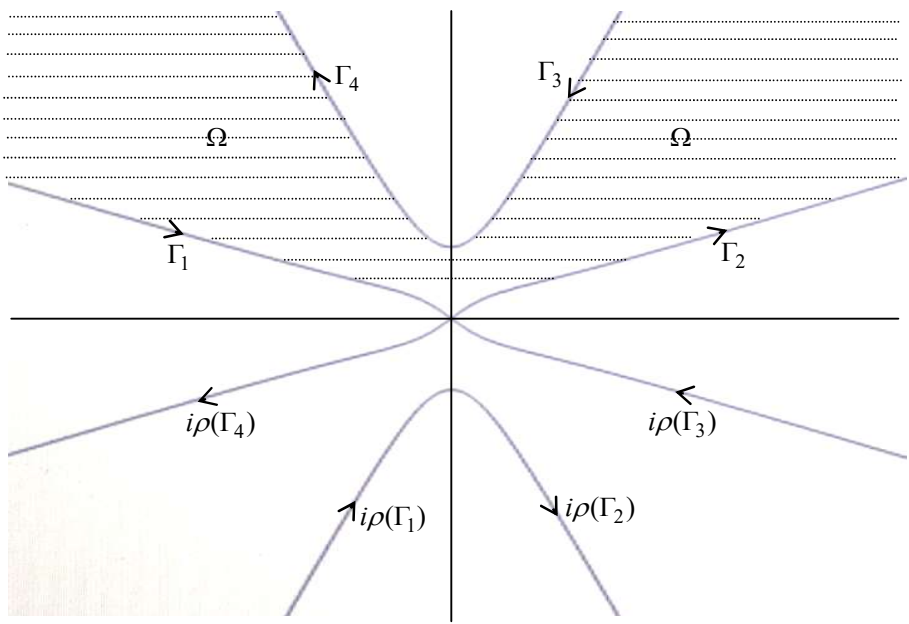

**Fig6** The images of the contours $\Gamma_1$, $\Gamma_2$, $\Gamma_3$, $\Gamma_4$ under the map $\lambda \rightarrow i\rho(\lambda)$. Also $\Gamma_3 = -i\rho(\Gamma_1)$, $\Gamma_2 = -i\rho(\Gamma_4)$, $\Gamma_1 = -i\rho(\Gamma_3)$ and $\Gamma_4 = -i\rho(\Gamma_2)$.

It follows from (4.2) and (4.5) that

$$\hat{u}_0(-\lambda) - \hat{u}(-\lambda,t)e^{\omega(\lambda)t} + \widetilde{g}_3(\omega(\lambda),t) - i\lambda\widetilde{g}_2(\omega(\lambda),t) - (\lambda^2+1)\widetilde{g}_1(\omega(\lambda),t) + i\lambda(\lambda^2+1)\widetilde{g}_0(\omega(\lambda),t) = 0 , \ \text{Im}\,\lambda \ge 0 , \tag{4.7}$$

and

$$\hat{u}_0(i\rho(\lambda)) - \hat{u}(i\rho(\lambda),t)e^{\omega(\lambda)t} + \widetilde{g}_3(\omega(\lambda),t) - \rho(\lambda)\widetilde{g}_2(\omega(\lambda),t)$$

$$+ \lambda^2 \widetilde{g}_1(\omega(\lambda),t) - \lambda^2 \rho(\lambda)\widetilde{g}_0(\omega(\lambda),t) = 0 \,, \ \mathrm{Im}[i\rho(\lambda)] \le 0 \,. \quad (4.8)$$

Keeping in mind (4.6), solving the system of equations (4.7) and (4.8) for $\widetilde{g}_2(\omega(\lambda),t)$ and $\widetilde{g}_3(\omega(\lambda),t)$, with $\lambda \in \partial\Omega$, substituting the results in equation (4.4) and taking into consideration that the integrals which involve the terms $\hat{u}(-\lambda,t)$ and $\hat{u}(i\rho(\lambda),t)$ vanish (by Cauchy's theorem, Jordan's lemma and the analyticity of $\rho$ in $\Theta$), we obtain the following integral representation of the solution of problem (4.1): For $x > 0$ and $t > 0$,

$$2\pi u(x,t) = \int_{-\infty}^{\infty} e^{i\lambda x - \omega(\lambda)t}\hat{u}_0(\lambda)d\lambda + \int_{\partial\Omega} e^{i\lambda x - \omega(\lambda)t}\left[\frac{i\lambda + \rho(\lambda)}{i\lambda - \rho(\lambda)}\hat{u}_0(-\lambda) - \frac{2i\lambda}{i\lambda - \rho(\lambda)}\hat{u}_0(i\rho(\lambda))\right]d\lambda$$
$$- 2i\int_{\partial\Omega} e^{i\lambda x - \omega(\lambda)t}\lambda[i\lambda + \rho(\lambda)][\rho(\lambda)\widetilde{g}_0(\omega(\lambda),t) - \widetilde{g}_1(\omega(\lambda),t)]d\lambda \,. \quad (4.9)$$

The above formula can be written also in the following equivalent way: For fixed $T > 0$,

$$2\pi u(x,t) = \int_{-\infty}^{\infty} e^{i\lambda x - \omega(\lambda)t}\hat{u}_0(\lambda)d\lambda + \int_{\partial\Omega} e^{i\lambda x - \omega(\lambda)t}\left[\frac{i\lambda + \rho(\lambda)}{i\lambda - \rho(\lambda)}\hat{u}_0(-\lambda) - \frac{2i\lambda}{i\lambda - \rho(\lambda)}\hat{u}_0(i\rho(\lambda))\right]d\lambda$$
$$- 2i\int_{\partial\Omega} e^{i\lambda x - \omega(\lambda)t}\lambda[i\lambda + \rho(\lambda)][\rho(\lambda)\widetilde{g}_0(\omega(\lambda),T) - \widetilde{g}_1(\omega(\lambda),T)]d\lambda \,, \quad (4.10)$$

for $x > 0$ and $0 < t < T$.

**Theorem 5** $1^{st}$ *The function $u(x,t)$, defined by (4.9) (equivalently by (4.10) is $C^{\infty}$ (jointly) for $(x,t) \in Q$ and satisfies the differential equation $\partial_t u = \partial_x^2 u - \partial_x^4 u$.*

$2^{nd}$ *For fixed $x_0 > 0$, $\lim\limits_{t \to 0^+} u(x,t) = u_0(x)$, uniformly for $x \ge x_0$.*

$3^{rd}$ *For fixed $t_0 > 0$, $\lim\limits_{x \to 0^+} u(x,t) = g_0(t)$ and $\lim\limits_{x \to 0^+} \dfrac{\partial u(x,t)}{\partial x} = g_1(t)$, uniformly for $t \ge t_0$.*

$4^{th}$ *The function $u(x,t)$, extended to $\overline{Q} - \{(0,0)\}$ by setting*

$$u*(x,0) := u_0(x) \ for \ x > 0 \ and \ u*(0,t) := g_0(t) \ for \ t > 0 \,,$$

*is $C^{\infty}$ (in $\overline{Q} - \{(0,0)\}$).*

$5^{th}$ *For fixed $x_0 > 0$, $\lim\limits_{t \to 0^+} \dfrac{\partial^n u(x,t)}{\partial x^n} = \dfrac{d^n u_0(x)}{dx^n}$, uniformly for $x \ge x_0$.*

$6^{th}$ *For fixed $t_0 > 0$, $\lim\limits_{x \to 0^+} \dfrac{\partial^m u(x,t)}{\partial t^m} = \dfrac{d^m g_0(t)}{dt^m}$ and $\lim\limits_{x \to 0^+} \dfrac{\partial^{m+1} u(x,t)}{\partial t^m \partial x} = \dfrac{d^m g_1(t)}{dt^m}$, uniformly for $t \ge t_0$.*

$7^{th}$ *For fixed $t_0 > 0$, the function $u(x,t)$, satisfies*

$$\lim\limits_{x \to +\infty}\left(x^{\ell} \frac{\partial^{n+m} u(x,t)}{\partial x^n \partial t^m}\right) = 0 \,,$$

*for every nonnegative integers $n$, $m$ and $\ell$, uniformly for $0 < t \le t_0$.*

**Proof** Since the proofs of the assertions of this theorem are similar to the proofs of Theorems 1-3, we will give only the necessary modifications, and these are the results in steps 1-3, below.

**Step 1** We claim that

$$\left[-2i\int_{\partial\Omega} e^{i\lambda x - \omega(\lambda)t}\lambda(i\lambda + \rho)[\rho\widetilde{g}_0(\omega(\lambda),T) - \widetilde{g}_1(\omega(\lambda),T)]d\lambda\right]\Bigg|_{x=0} =$$
$$- 2i\int_{\partial\Omega} e^{-\omega(\lambda)t}\lambda(i\lambda + \rho)[\rho\widetilde{g}_0(\omega(\lambda),T) - \widetilde{g}_1(\omega(\lambda),T)]d\lambda = 2\pi g_0(t) \,. \quad (4.11)$$

To prove this, first, we set $\mu = -i\rho(\lambda)$, in the integral $\int_{\Gamma_1 + \Gamma_4} \cdots d\lambda$, which is part of the above integral $\int_{\partial\Omega} \cdots d\lambda$, and we obtain

$$-2i\int_{\partial\Omega} e^{-\omega(\lambda)t}\lambda(i\lambda + \rho)[\rho\widetilde{g}_0(\omega(\lambda),T)]d\lambda$$

$$= -2i \int_{\Gamma_1 + \Gamma_4} e^{-\omega(\lambda)t} \lambda(i\lambda + \rho)\rho \widetilde{g_0}(\omega(\lambda),T)d\lambda - 2i \int_{\Gamma_2 + \Gamma_3} e^{-\omega(\lambda)t} \lambda(i\lambda + \rho)\rho \widetilde{g_0}(\omega(\lambda),T)d\lambda$$

$$= \int_{\Gamma_2 + \Gamma_3} e^{-\omega(\lambda)t}\{[-2\lambda^2\rho(\lambda) - 2i\lambda^3] - 2i\lambda[i\lambda + \rho(\lambda)]\rho(\lambda)\}\widetilde{g_0}(\omega(\lambda),T)d\lambda = -\int_{\Gamma_2 + \Gamma_3} e^{-\omega(\lambda)t}(4i\lambda^3 + 2i\lambda)\widetilde{g_0}(\omega(\lambda),T)d\lambda .$$

$$(4.12)$$

Indeed, with this change of variable, the differential

$$-2i\lambda[i\lambda + \rho(\lambda)]\rho(\lambda)d\lambda\Big|_{\lambda \in \Gamma_1 + \Gamma_4} = [-2\mu^2\rho(\mu) - 2i\mu^3]d\mu\Big|_{\mu \in \Gamma_3 + \Gamma_4} ,$$

since, for $\lambda \in \Gamma_1 + \Gamma_4$, we have

$$\mu = -i\rho(\lambda) \in \Gamma_2 + \Gamma_3 \text{ (see fig6)}, \ \lambda = -i\rho(\mu), \ \lambda^2 = -\mu^2 - 1 \text{ and } \lambda d\lambda = -\mu d\mu .$$

Setting $\kappa = i\omega(\lambda) = i(\lambda^4 + \lambda^2)$, the last integral in (4.12) becomes

$$integral \ (4.12) = \int_{\kappa = -\infty}^{\infty} e^{i\kappa t}\left(\int_{\tau = 0}^{T} e^{-i\kappa\tau} g_0(\tau)d\tau\right)d\kappa = 2\pi g_0(t) . \qquad (4.13)$$

Similarly, since $2i\lambda[i\lambda + \rho(\lambda)]d\lambda\Big|_{\lambda \in \Gamma_1 + \Gamma_4} = -2i\mu[i\mu + \rho(\mu)]d\mu\Big|_{\mu = -i\rho(\lambda) \in \Gamma_2 + \Gamma_3}$, we obtain

$$-2i \int_{\partial\Omega} e^{-\omega(\lambda)t} \lambda(i\lambda + \rho)[-\widetilde{g_1}(\omega(\lambda),T)]d\lambda$$

$$= -2i \int_{\Gamma_1 + \Gamma_4} e^{-\omega(\lambda)t} \lambda(i\lambda + \rho)[-\widetilde{g_1}(\omega(\lambda),T)]d\lambda - 2i \int_{\Gamma_2 + \Gamma_3} e^{-\omega(\lambda)t} \lambda(i\lambda + \rho)[-\widetilde{g_1}(\omega(\lambda),T)]d\lambda = 0 . \quad (4.14)$$

Now (4.11) follows from (4.12), (4.13) and (4.14).

***Step 2*** We claim that

$$\left[\int_{-\infty}^{\infty} e^{i\lambda x - \omega(\lambda)t} \hat{u}_0(\lambda)d\lambda + \int_{\partial\Omega} e^{i\lambda x - \omega(\lambda)t}\left[\frac{i\lambda + \rho}{i\lambda - \rho}\hat{u}_0(-\lambda) - \frac{2i\lambda}{i\lambda - \rho}\hat{u}_0(i\rho(\lambda))\right]d\lambda\right]\Bigg|_{x=0}$$

$$= \int_{-\infty}^{\infty} e^{-\omega(\lambda)t} \hat{u}_0(\lambda)d\lambda + \int_{\partial\Omega} e^{-\omega(\lambda)t}\left[\frac{i\lambda + \rho}{i\lambda - \rho}\hat{u}_0(-\lambda) - \frac{2i\lambda}{i\lambda - \rho}\hat{u}_0(i\rho(\lambda))\right]d\lambda = 0 . \quad (4.15)$$

Indeed, since $\mu = i\rho(\lambda)\Big|_{\lambda \in \Gamma_3 + \Gamma_4} \Leftrightarrow \lambda = -i\rho(\mu)\Big|_{\mu \in i\rho(\Gamma_3) + i\rho(\Gamma_4)}$, we have

$$\int_{\Gamma_1 + \Gamma_2} e^{-\omega(\lambda)t} \frac{i\lambda + \rho(\lambda)}{i\lambda - \rho(\lambda)}\hat{u}_0(-\lambda)d\lambda - \int_{\Gamma_3 + \Gamma_4} e^{-\omega(\lambda)t} \frac{2i\lambda}{i\lambda - \rho(\lambda)}\hat{u}_0(i\rho(\lambda))d\lambda$$

$$= \int_{-\infty}^{\infty} e^{-\omega(\lambda)t} \frac{i\lambda + \rho(\lambda)}{i\lambda - \rho(\lambda)}\hat{u}_0(-\lambda)d\lambda + \int_{i\rho(\Gamma_3) + i\rho(\Gamma_4)} e^{-\omega(\mu)t} \frac{2i\mu}{\rho(\mu) + i\mu}\hat{u}_0(\mu)d\mu$$

$$= \int_{-\infty}^{\infty} e^{-\omega(\lambda)t} \frac{-i\lambda + \rho(\lambda)}{i\lambda - \rho(\lambda)}\hat{u}_0(\lambda)d\lambda - \int_{-\infty}^{\infty} e^{-\omega(\mu)t} \frac{2i\mu}{\rho(\mu) + i\mu}\hat{u}_0(\mu)d\mu = -\int_{-\infty}^{\infty} e^{-\omega(\lambda)t}\hat{u}_0(\lambda)d\lambda . \quad (4.16)$$

Also, since $\lambda = -i\rho(\mu)\Big|_{\mu \in \Gamma_1 + \Gamma_2} \Leftrightarrow \mu = -i\rho(\lambda)\Big|_{\lambda \in \Gamma_3 + \Gamma_4}$,

$$\int_{\Gamma_3 + \Gamma_4} e^{-\omega(\lambda)t} \frac{i\lambda + \rho(\lambda)}{i\lambda - \rho(\lambda)}\hat{u}_0(-\lambda)d\lambda - \int_{\Gamma_1 + \Gamma_2} e^{-\omega(\lambda)t} \frac{2i\lambda}{i\lambda - \rho}\hat{u}_0(i\rho(\lambda))d\lambda$$

$$= \int_{\Gamma_1 + \Gamma_2} e^{-\omega(\mu)t} \frac{\rho(\mu) + i\mu}{\rho(\mu) - i\mu}\hat{u}_0(i\rho(\mu))\frac{\mu d\mu}{i\rho(\mu)} - \int_{-\infty}^{\infty} e^{-\omega(\lambda)t} \frac{2i\lambda}{i\lambda - \rho(\lambda)}\hat{u}_0(i\rho(\lambda))d\lambda$$

$$= \int_{-\infty}^{\infty} e^{-\omega(\lambda)t} \frac{\rho(\lambda) + i\lambda}{\rho(\lambda) - i\lambda}\hat{u}_0(i\rho(\lambda))\frac{\lambda d\lambda}{i\rho(\lambda)} - \int_{-\infty}^{\infty} e^{-\omega(\lambda)t} \frac{2i\lambda}{i\lambda - \rho(\lambda)}\hat{u}_0(i\rho(\lambda))d\lambda = \int_{-\infty}^{\infty} e^{-\omega(\lambda)t} \frac{-\lambda}{i\rho(\lambda)}d\lambda = 0 . \quad (4.17)$$

Now (4.15) follows from (4.16) and (4.17).

***Step 3*** Differentiating (4.10), we obtain

$$2\pi\frac{\partial u(x,t)}{\partial x}\Bigg|_{x=0} = \int_{-\infty}^{\infty} (i\lambda)e^{-\omega(\lambda)t} \hat{u}_0(\lambda)d\lambda + \int_{\partial\Omega} e^{-\omega(\lambda)t}\left[(i\lambda)\frac{i\lambda + \rho}{i\lambda - \rho}\hat{u}_0(-\lambda) - (i\lambda)\frac{2i\lambda}{i\lambda - \rho}\hat{u}_0(i\rho(\lambda))\right]d\lambda$$

$$-2i \int_{\partial\Omega} (i\lambda) e^{-\omega(\lambda)t} \lambda(i\lambda+\rho) [\rho \widetilde{g}_0(\omega(\lambda),T) - \widetilde{g}_1(\omega(\lambda),T)] d\lambda \; . \tag{4.18}$$

We claim that

$$-2i \int_{\partial\Omega} (i\lambda) e^{-\omega(\lambda)t} \lambda(i\lambda+\rho) [\rho \widetilde{g}_0(\omega(\lambda),t) - \widetilde{g}_1(\omega(\lambda),t)] d\lambda = 2\pi g_1(t) \; . \tag{4.19}$$

Indeed, setting $\mu = -i\rho(\lambda)$, in the integral $\int_{\Gamma_1+\Gamma_4} \cdots d\lambda$, which is part of the above integral, and taking into consideration that

$$-2i(i\lambda)\lambda[i\lambda+\rho(\lambda)]\rho(\lambda)d\lambda \Big|_{\lambda\in\Gamma_1+\Gamma_4} = 2i(i\mu)\mu[\rho(\mu)+i\mu]\rho(\mu)d\mu \Big|_{\mu\in\Gamma_3+\Gamma_4}$$

and

$$-2i(i\lambda)\lambda[i\lambda+\rho(\lambda)](-1)d\lambda \Big|_{\lambda\in\Gamma_1+\Gamma_4} = [-2i\mu^3 - 2i\mu + 2\mu^2\rho(\mu)]d\mu \Big|_{\mu\in\Gamma_3+\Gamma_4} \; ,$$

we see that

$$integral \; (4.19) = -\int_{\Gamma_2+\Gamma_3} e^{-\omega(\lambda)t}(4i\lambda^3+2i\lambda)\widetilde{g}_1(\omega(\lambda),T)d\lambda = \int_{\kappa=-\infty}^{\infty} e^{i\kappa t}\left(\int_{\tau=0}^{T} e^{-i\kappa\tau} g_1(\tau)d\tau\right)d\kappa = 2\pi g_1(t) \; , \tag{4.20}$$

where we also set $\kappa = i\omega(\lambda) = i(\lambda^4+\lambda^2)$.

**Step 4** We claim that

$$\int_{-\infty}^{\infty} (i\lambda) e^{-\omega(\lambda)t} \hat{u}_0(\lambda)d\lambda + \int_{\partial\Omega} e^{-\omega(\lambda)t}\left[(i\lambda)\frac{i\lambda+\rho}{i\lambda-\rho}\hat{u}_0(-\lambda) - (i\lambda)\frac{2i\lambda}{i\lambda-\rho}\hat{u}_0(i\rho(\lambda))\right]d\lambda = 0 \; . \tag{4.21}$$

The proof is similar to the proof of (4.15).

### 4.3 The inhomogeneous equation

**Problem** Solve (1.1) in the case of the equation $\partial_t U = \partial_x^2 U - \partial_x^4 U + f$, i.e., when $\alpha = 1$ and $\beta = 1$. $\tag{4.22}$

**Solution** With

$$\Phi_{\mathbb{R}}(x,t) = \int_{-\infty}^{\infty} e^{i\lambda x - \omega(\lambda)t} \widetilde{\widetilde{f}}(\lambda,\omega(\lambda),t)d\lambda$$

and

$$\Phi_{\partial\Omega}(x,t) = \int_{\partial\Omega} e^{i\lambda x - \omega(\lambda)t}\left[\frac{i\lambda+\rho(\lambda)}{i\lambda-\rho(\lambda)}\widetilde{\widetilde{f}}(-\lambda,\omega(\lambda),t) - \frac{2i\lambda}{i\lambda-\rho(\lambda)}\widetilde{\widetilde{f}}(i\rho(\lambda),\omega(\lambda),t))\right]d\lambda \; ,$$

the solution of (4.22) is given by the formula

$$U(x,t) = u(x,t) + \frac{1}{2\pi}\Phi_{\mathbb{R}}(x,t) + \frac{1}{2\pi}\Phi_{\partial\Omega}(x,t) \; , \tag{4.23}$$

where $u(x,t)$ is the solution of the homogeneous equation given by (4.10).

**Verification and boundary behavior of the solution** Having in mind the results and the techniques of the previous paragraphs, it is easy to state and prove a theorem analogous to Theorem 4, also for the solution (4.23) of problem (4.22).

## 5. The equation $\partial_t U = \alpha\partial_x^2 U - \beta\partial_x^4 U + f$

**Problem** With $\alpha > 0$ and $\beta > 0$, solve

$$\begin{cases} \partial_t U = \alpha \partial_x^2 U - \beta \partial_x^4 U + f, \ (x,t) \in Q := \mathbb{R}^+ \times \mathbb{R}^+, \\ U(x,0) = u_0(x), \ x \in \mathbb{R}^+, \\ U(0,t) = g_0(t), \ t \in \mathbb{R}^+, \\ U_x(0,t) = g_1(t), \ t \in \mathbb{R}^+. \end{cases} \qquad (5.1)$$

***Notation*** In order to write the integral representation of the solution of problem (5.1), we ***adjust*** the notation used in section 4, as far as the dispersion relation $\omega(\lambda)$, the sets $\Omega$ and $\Theta$, and the function $\rho(\lambda)$ are concerned. More precisely, in this section, we define

$$\omega(\lambda) = \omega_{\alpha,\beta}(\lambda) = \alpha\lambda^2 + \beta\lambda^4,$$

$$\Omega = \Omega_{\alpha,\beta} = \{\lambda : \operatorname{Im}\lambda \geq 0 \ \& \ \operatorname{Re}\omega_{\alpha,\beta}(\lambda) = \operatorname{Re}(\alpha\lambda^2 + \beta\lambda^4) \leq 0\} \ (\text{see fig7 and fig8}),$$

$$\rho(\lambda) = \rho_{\alpha,\beta}(\lambda) = \sqrt{\left|\lambda^2 + \frac{\alpha}{\beta}\right|} \exp\{i[\theta_1(\lambda) + \theta_2(\lambda) + \pi]/2\} \ (\text{see fig9})$$

and

$$\Theta = \Theta_{\alpha,\beta} = \mathbb{C} - \{\lambda = \xi + i\eta \in \mathbb{C} : \xi = 0 \ \& \ |\eta| \geq \sqrt{\alpha/\beta}\} = \mathbb{C} - \{(-i\infty, -\sqrt{\tfrac{\alpha}{\beta}}i] \cup [\sqrt{\tfrac{\alpha}{\beta}}i, i\infty)\}.$$

If we write $\lambda = \xi + i\eta$ then

$$\Omega = \Omega_{\alpha,\beta} = \{\xi + i\eta : \psi_{\alpha,\beta}^-(\xi) \leq \eta \leq \psi_{\alpha,\beta}^+(\xi), \ \forall \xi \in \mathbb{R}\}$$

and

$$\partial\Omega = \partial\Omega_{\alpha,\beta} = \{\xi + i\eta : \eta = \psi_{\alpha,\beta}^-(\xi) \ or \ \eta = \psi_{\alpha,\beta}^+(\xi), \ for \ some \ \xi \in \mathbb{R}\},$$

where

$$\psi_{\alpha,\beta}^-(\xi) = \sqrt{\frac{6\beta\xi^2 + \alpha - \sqrt{32\beta^2\xi^4 + 8\alpha\beta\xi^2 + \alpha^2}}{2\beta}} \ \text{ and } \ \psi_{\alpha,\beta}^+(\xi) = \sqrt{\frac{6\beta\xi^2 + \alpha + \sqrt{32\beta^2\xi^4 + 8\alpha\beta\xi^2 + \alpha^2}}{2\beta}} \ .$$

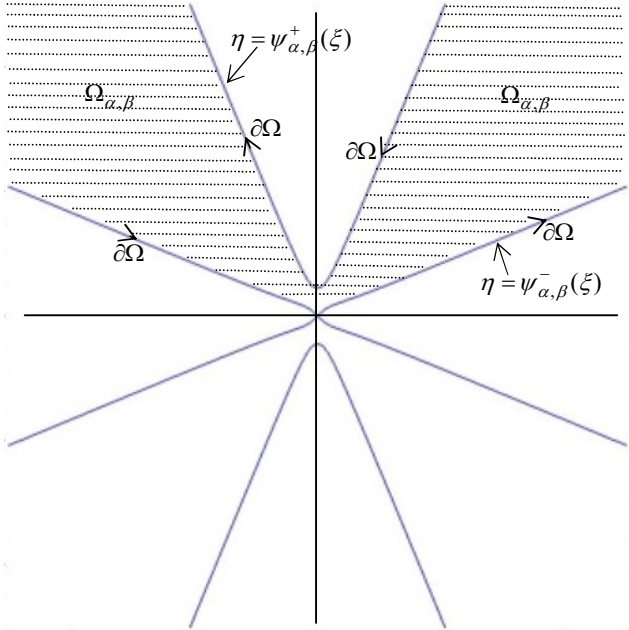

**Fig.7** The set $\Omega_{\alpha,\beta}$ with $\alpha = 1$ and $\beta = 5$.

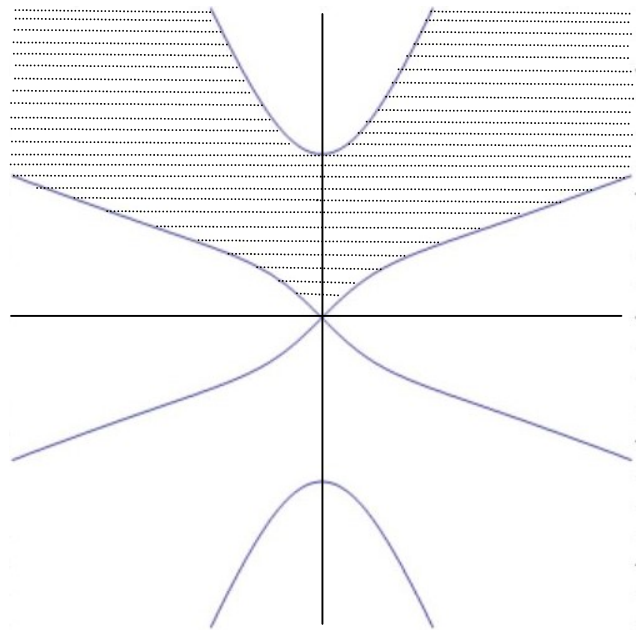

**Fig.8** The set $\Omega_{\alpha,\beta}$ with $\alpha = 7$ and $\beta = 1$.

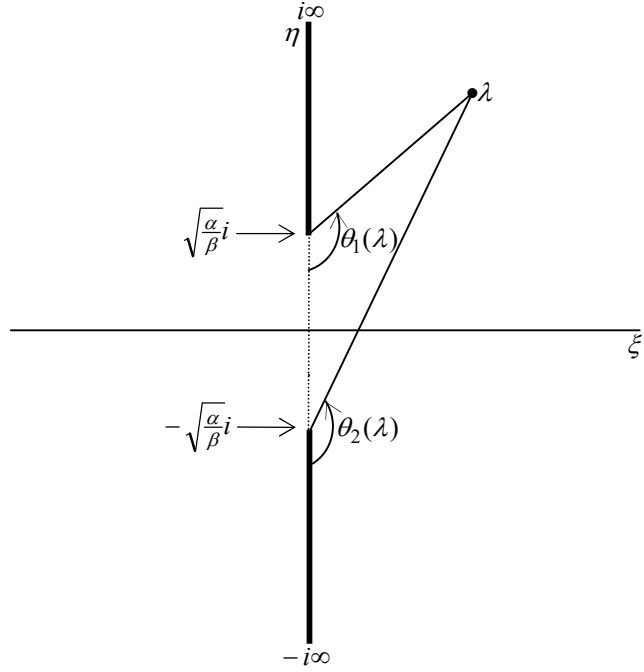

**Fig.9** With $-\pi < \theta_1(\lambda) < \pi$ and $0 < \theta_2(\lambda) < 2\pi$, the function $\rho_{\alpha,\beta}(\lambda)$ is analytic for $\lambda \in \Theta_{\alpha,\beta}$.

***The equation in divergence form*** The equation $\partial_t U = \alpha \partial_x^2 U - \beta \partial_x^4 U + f$ can be written in the form:

$$\frac{\partial}{\partial t}[e^{-i\lambda x + \omega(\lambda)t}U(x,t)] - \frac{\partial}{\partial x}\{e^{-i\lambda x + \omega(\lambda)t}[-\beta U_{xxx} - \beta i\lambda U_{xx} + (\beta\lambda^2 + \alpha)U_x + i\lambda(\beta\lambda^2 + \alpha)U]\} = e^{-i\lambda x + \omega(\lambda)t}f(x,t),$$

and Green's theorem gives the global relation:

$$\frac{1}{\beta}\hat{U}(\lambda,t)e^{\omega(\lambda)t} - \frac{1}{\beta}\hat{u}_0(\lambda) - [\widetilde{g}_3(\omega(\lambda),t) + i\lambda\widetilde{g}_2(\omega(\lambda),t)]$$

$$+ \left(\lambda^2 + \frac{\alpha}{\beta}\right)\widetilde{g}_1(\omega(\lambda),t) + i\lambda\left(\lambda^2 + \frac{\alpha}{\beta}\right)\widetilde{g}_0(\omega(\lambda),t) - \frac{1}{\beta}\widetilde{\hat{f}}(\lambda,\omega(\lambda),t) = 0, \ \operatorname{Im}\lambda \le 0. \quad (5.2)$$

Multiplying (5.2) by $e^{i\lambda x - \omega(\lambda)t}$ and integrating, we obtain

$$\frac{2\pi}{\beta}U(x,t) = \frac{1}{\beta}\int_{-\infty}^{\infty}e^{-i\lambda x + \omega(\lambda)t}\hat{u}_0(\lambda)d\lambda - \int_{\partial\Omega}e^{-i\lambda x + \omega(\lambda)t}\rho^2[i\lambda\widetilde{g}_0(\omega(\lambda),t) + \widetilde{g}_1(\omega(\lambda),t)]d\lambda$$

$$+ \frac{1}{\beta}\int_{-\infty}^{\infty}e^{-i\lambda x + \omega(\lambda)t}\widetilde{\hat{f}}(\lambda,\omega(\lambda),t)d\lambda + \int_{\partial\Omega}e^{-i\lambda x + \omega(\lambda)t}[\widetilde{g}_3(\omega(\lambda),t) + i\lambda\widetilde{g}_2(\omega(\lambda),t)]d\lambda. \quad (5.3)$$

Setting $-\lambda$ and $i\rho(\lambda)$, in place of $\lambda$, in (5.2), we obtain

$$\widetilde{g}_3(\omega(\lambda),t) - i\lambda\widetilde{g}_2(\omega(\lambda),t) = \frac{1}{\beta}\hat{U}(-\lambda,t)e^{\omega(\lambda)t} - \frac{1}{\beta}\hat{u}_0(-\lambda)$$

$$+ \rho^2\widetilde{g}_1(\omega(\lambda),t) - i\lambda\rho^2\widetilde{g}_0(\omega(\lambda),t) - \frac{1}{\beta}\widetilde{\hat{f}}(-\lambda,\omega(\lambda),t), \ \operatorname{Im}\lambda \ge 0, \quad (5.4)$$

and

$$\widetilde{g}_3(\omega(\lambda),t) - \rho(\lambda)\lambda\widetilde{g}_2(\omega(\lambda),t) = \frac{1}{\beta}\hat{U}(i\rho(\lambda),t)e^{\omega(\lambda)t} - \frac{1}{\beta}\hat{u}_0(i\rho(\lambda))$$

$$- \lambda^2\widetilde{g}_1(\omega(\lambda),t) + i\lambda\rho^2\widetilde{g}_0(\omega(\lambda),t) - \frac{1}{\beta}\widetilde{\hat{f}}(i\rho(\lambda),\omega(\lambda),t), \ \operatorname{Im}[i\rho(\lambda)] \le 0. \quad (5.5)$$

Solving the system of equations (5.4) and (5.5) for $\widetilde{g}_2(\omega(\lambda),t)$ and $\widetilde{g}_3(\omega(\lambda),t)$, with $\lambda \in \partial\Omega$, and substituting their values in (5.3), we obtain the solution (5.6), below.

**Solution** With the above notation for $\omega = \omega_{\alpha,\beta}$, $\Omega = \Omega_{\alpha,\beta}$ and $\rho = \rho_{\alpha,\beta}$, the Fokas method solution of the IBVP (5.1) is given by the following formula: For $x > 0$ and $t > 0$,

$$U(x,t) = \frac{1}{2\pi}\int_{-\infty}^{\infty}e^{i\lambda x - \omega(\lambda)t}\hat{u}_0(\lambda)d\lambda + \frac{1}{2\pi}\int_{\partial\Omega}e^{i\lambda x - \omega(\lambda)t}\left[\frac{i\lambda + \rho(\lambda)}{i\lambda - \rho(\lambda)}\hat{u}_0(-\lambda) - \frac{2i\lambda}{i\lambda - \rho(\lambda)}\hat{u}_0(i\rho(\lambda))\right]d\lambda$$

$$-\frac{\beta i}{\pi}\int_{\partial\Omega}e^{i\lambda x - \omega(\lambda)t}\lambda[i\lambda + \rho(\lambda)][\rho(\lambda)\widetilde{g}_0(\omega(\lambda),t) - \widetilde{g}_1(\omega(\lambda),t)]d\lambda$$

$$+\frac{1}{2\pi}\int_{-\infty}^{\infty}e^{i\lambda x - \omega(\lambda)t}\widetilde{\widehat{f}}(\lambda,\omega(\lambda),t)d\lambda + \frac{1}{2\pi}\int_{\partial\Omega}e^{i\lambda x - \omega(\lambda)t}\left[\frac{i\lambda + \rho(\lambda)}{i\lambda - \rho(\lambda)}\widetilde{\widehat{f}}(-\lambda,\omega(\lambda),t) - \frac{2i\lambda}{i\lambda - \rho(\lambda)}\widetilde{\widehat{f}}(i\rho(\lambda),\omega(\lambda),t))\right]d\lambda \ . \ (5.6)$$

The above formula can be written also in the following way: For $x > 0$ and $t > 0$,

$$U(x,t) = \frac{1}{2\pi}\int_{-\infty}^{\infty}e^{i\lambda x - \omega(\lambda)t}\hat{u}_0(\lambda)d\lambda + \frac{1}{2\pi}\int_{-\infty}^{\infty}e^{i\lambda x - \omega(\lambda)t}\left[\frac{i\lambda + \rho(\lambda)}{i\lambda - \rho(\lambda)}\hat{u}_0(-\lambda) - \frac{2i\lambda}{i\lambda - \rho(\lambda)}\hat{u}_0(i\rho(\lambda))\right]d\lambda$$

$$+\frac{1}{2\pi}\int_{-\infty}^{\infty}e^{\rho(\mu)x - \omega(\mu)t}\left[\frac{i\mu - \rho(\mu)}{i\mu + \rho(\mu)}\hat{u}_0(i\rho(\mu)) + \frac{2\rho(\mu)}{i\mu + \rho(\mu)}\hat{u}_0(\mu)\right]\frac{\mu d\mu}{i\rho(\mu)}$$

$$-\frac{\beta i}{\pi}\int_{-\infty}^{\infty}e^{i\lambda x - \omega(\lambda)t}\lambda[i\lambda + \rho(\lambda)][\rho(\lambda)\widetilde{g}_0(\omega(\lambda),t) - \widetilde{g}_1(\omega(\lambda),t)]d\lambda$$

$$-\frac{\beta i}{\pi}\int_{-\infty}^{\infty}e^{\rho(\mu)x - \omega(\mu)t}\mu[i\mu - \rho(\mu)][i\mu\widetilde{g}_0(\omega(\mu),t) + \widetilde{g}_1(\omega(\mu),t)]d\mu$$

$$+\frac{1}{2\pi}\int_{-\infty}^{\infty}e^{i\lambda x - \omega(\lambda)t}\widetilde{\widehat{f}}(\lambda,\omega(\lambda),t)d\lambda + \frac{1}{2\pi}\int_{-\infty}^{\infty}e^{i\lambda x - \omega(\lambda)t}\left[\frac{i\lambda + \rho(\lambda)}{i\lambda - \rho(\lambda)}\widetilde{\widehat{f}}(-\lambda,\omega(\lambda),t) - \frac{2i\lambda}{i\lambda - \rho(\lambda)}\widetilde{\widehat{f}}(i\rho(\lambda),\omega(\lambda),t))\right]d\lambda$$

$$+\frac{1}{2\pi}\int_{-\infty}^{\infty}e^{\rho(\mu)x - \omega(\mu)t}\left[\frac{i\mu - \rho(\mu)}{i\mu + \rho(\mu)}\widetilde{\widehat{f}}(i\rho(\mu),\omega(\mu),t) + \frac{2\rho(\mu)}{i\mu + \rho(\mu)}\widetilde{\widehat{f}}(\mu,\omega(\mu),t))\right]\frac{\mu d\mu}{i\rho(\mu)} \ . \ (5.7)$$

*Proof of (5.7)* Let us write $\partial\Omega = C_{lower} + C_{upper}$ where $C_{lower}$ is the connected component of $\partial\Omega$ which passes through 0 and $C_{upper}$ is the connected component of $\partial\Omega$ which passes through the point $\sqrt{\frac{\alpha}{\beta}}i$. Then, by Cauchy's theorem and Jordan's lemma,

$$\int_{C_{lower}}\mathcal{I}_{u_0} = \int_{-\infty}^{\infty}\mathcal{I}_{u_0} \ , \quad \int_{C_{lower}}\mathcal{I}_{g_0,g_1} = \int_{-\infty}^{\infty}\mathcal{I}_{g_0,g_1} \text{ and } \int_{C_{lower}}\mathcal{I}_f = \int_{-\infty}^{\infty}\mathcal{I}_f \ , \tag{5.8}$$

where $\mathcal{I}_{u_0}$, $\mathcal{I}_{g_0,g_1}$ and $\mathcal{I}_f$ are the integrands of the integrals $\int_{\partial\Omega}\cdots$ in (5.6) which contain $u_0$, $(g_0,g_1)$ and $f$, respectively.

On the other hand,

$$\int_{C_{upper}^+}\mathcal{I}_{u_0} = \int_{\sqrt{\alpha/\beta}i}^{i\infty}\mathcal{I}_{u_0} \ , \quad \int_{C_{upper}^+}\mathcal{I}_{g_0,g_1} = \int_{\sqrt{\alpha/\beta}i}^{i\infty}\mathcal{I}_{g_0,g_1} \text{ and } \int_{C_{upper}^+}\mathcal{I}_f = \int_{\sqrt{\alpha/\beta}i}^{i\infty}\mathcal{I}_f \ , \tag{5.9}$$

where $C_{upper}^+ = C_{upper} \cap \{\text{Re}\,\lambda \geq 0\}$ and $\rho(\lambda)$, for $\lambda \in [\sqrt{\frac{\alpha}{\beta}}i, i\infty]$, is the continuous extension of $\rho(\lambda)$, from the side of the quandrant $\{\text{Re}\,\lambda > 0, \text{Im}\,\lambda > 0\}$ .

Similarly,

$$\int_{C_{upper}^-}\mathcal{I}_{u_0} = \int_{\sqrt{\alpha/\beta}i}^{i\infty}\mathcal{I}_{u_0} \ , \quad \int_{C_{upper}^-}\mathcal{I}_{g_0,g_1} = \int_{\sqrt{\alpha/\beta}i}^{i\infty}\mathcal{I}_{g_0,g_1} \text{ and } \int_{C_{upper}^-}\mathcal{I}_f = \int_{\sqrt{\alpha/\beta}i}^{i\infty}\mathcal{I}_f \ , \tag{5.10}$$

where $C_{upper}^- = C_{upper} \cap \{\text{Re}\,\lambda \leq 0\}$ and, *this time*, $\rho(\lambda)$, for $\lambda \in [\sqrt{\frac{\alpha}{\beta}}i, i\infty]$, is the continuous extension of $\rho(\lambda)$, from the side of the quandrant $\{\text{Re}\,\lambda < 0, \text{Im}\,\lambda > 0\}$ .

Now, setting $\mu = i\rho(\lambda)$ in the integrals (with the appropriate choice of $\rho$, in each case)

$$\int\limits_{\sqrt{\alpha/\beta i}}^{i\infty}\mathcal{I}_{u_0}\ ,\quad \int\limits_{\sqrt{\alpha/\beta i}}^{i\infty}\mathcal{I}_{g_0,g_1}\ ,\quad \int\limits_{\sqrt{\alpha/\beta i}}^{i\infty}\mathcal{I}_f\ ,$$

of both (5.9) and (5.10), and substituting the resulting equations (5.9) and (5.10) in (5.6), we obtain (5.7), also in combination with (5.8).

As in the previous cases, a theorem, analogous to Theorem 4, can be stated and proved also in this case. In particular, the function $U^*(x,t)$, defined by

$$U^*(x,t)=\begin{cases}U(x,t)\ \text{for}\ (x,t)\in Q\\ u_0(x)\ \text{for}\ t=0\ \text{and}\ x>0\\ g_0(t)\ \text{for}\ x=0\ \text{and}\ t=0,\end{cases}\tag{5.11}$$

is $C^\infty$ for $(x,t)\in \overline{Q}-\{(0,0)\}$.

The next theorem gives sufficient conditions on the data for the existence of the limit of the derivatives of the solution (5.6), as $(x,t)\to(0,0)$.

**Theorem 6** Let $U^*(x,t)$ be the function, defined by (5.11). Then

$1^{st}$ If $u_0(0)=g_0(0)=0$ then the limit $\lim\limits_{\substack{(x,t)\to(0,0)\\(x,t)\in\overline{Q}-\{(0,0)\}}}U^*(x,t)$ exists.

$2^{nd}$ If $u_0(0)=g_0(0)=0$ and $u_0{}'(0)=g_1(0)=0$ then the limit $\lim\limits_{\substack{(x,t)\to(0,0)\\(x,t)\in\overline{Q}-\{(0,0)\}}}\dfrac{\partial U^*(x,t)}{\partial x}$ exists.

$3^{rd}$ If $u_0(0)=g_0(0)=0$, $u_0{}'(0)=g_1(0)=0$ and $u_0{}''(0)=0$ then the limit $\lim\limits_{\substack{(x,t)\to(0,0)\\(x,t)\in\overline{Q}-\{(0,0)\}}}\dfrac{\partial^2 U^*(x,t)}{\partial x^2}$ exists.

$4^{th}$ If $u_0(0)=g_0(0)=0$, $u_0{}'(0)=g_1(0)=0$ and $u_0{}''(0)=u_0{}'''(0)=0$ then the limit $\lim\limits_{\substack{(x,t)\to(0,0)\\(x,t)\in\overline{Q}-\{(0,0)\}}}\dfrac{\partial^3 U^*(x,t)}{\partial x^3}$ exists.

$5^{th}$ If $u_0(0)=g_0(0)=0$, $u_0{}'(0)=g_1(0)=0$, $u_0{}''(0)=u_0{}'''(0)=0$, $g_0{}'(0)=0$ and $\beta u_0{}''''(0)=f(0,0)$ then the limits

$$\lim\limits_{\substack{(x,t)\to(0,0)\\(x,t)\in\overline{Q}-\{(0,0)\}}}\dfrac{\partial^k U^*(x,t)}{\partial x^k}, \text{for}\ k=0,1,2,3,4,\ \text{and}\ \lim\limits_{\substack{(x,t)\to(0,0)\\(x,t)\in\overline{Q}-\{(0,0)\}}}\dfrac{\partial U^*(x,t)}{\partial t}$$

exist.

**Proof** The first four parts are easily established, using (5.6), expanding $\hat{u}_0$, $\widetilde{g}_j$ and $\widetilde{\tilde{f}}$, according to the formulas (1.4), (2.5), (3.8) and (3.9), and using the compatibility assumptions on the data, at the origin.

*Proof of $5^{th}$ part* This time we have to use (5.7). We introduce the following notation: For two functions $Z(x,t)$ and $W(x,t)$, defined for $x>0$ and $t>0$, we will write $Z(x,t)\approx W(x,t)$ if and only if the limits

$$\lim\limits_{Q\ni(x,t)\to(0,0)}\big[Z(x,t)-W(x,t)\big],\ \lim\limits_{x\to0^+}\left(\lim\limits_{\substack{t\to0^+\\x>0}}\big[Z(x,t)-W(x,t)\big]\right),\ \lim\limits_{t\to0^+}\left(\lim\limits_{\substack{x\to0^+\\t>0}}\big[Z(x,t)-W(x,t)\big]\right),$$

exist and are equal.

With this notation we have: For $x>0$ and $t>0$,

$$\dfrac{\partial^4}{\partial x^4}\left\{\int_0^\infty e^{i\lambda x-\omega(\lambda)t}\hat{u}_0(\lambda)d\lambda+\int_0^\infty e^{i\lambda x-\omega(\lambda)t}\left[\dfrac{i\lambda+\rho(\lambda)}{i\lambda-\rho(\lambda)}\hat{u}_0(-\lambda)-\dfrac{2i\lambda}{i\lambda-\rho(\lambda)}\hat{u}_0(i\rho(\lambda))\right]d\lambda\right\}$$

$$=\int_0^\infty (i\lambda)^4 e^{i\lambda x-\omega(\lambda)t}\hat{u}_0(\lambda)d\lambda+\int_0^\infty (i\lambda)^4 e^{i\lambda x-\omega(\lambda)t}\left[\dfrac{i\lambda+\rho(\lambda)}{i\lambda-\rho(\lambda)}\hat{u}_0(-\lambda)-\dfrac{2i\lambda}{i\lambda-\rho(\lambda)}\hat{u}_0(i\rho(\lambda))\right]d\lambda$$

$$\approx u_0^{(4)}(0)\left\{\int_1^\infty (i\lambda)^4 e^{i\lambda x - \omega(\lambda)t}\frac{d\lambda}{(i\lambda)^5} + \int_1^\infty (i\lambda)^4 e^{i\lambda x - \omega(\lambda)t}\frac{i\lambda + \rho(\lambda)}{i\lambda - \rho(\lambda)}\frac{d\lambda}{(-i\lambda)^5} - \int_1^\infty (i\lambda)^4 e^{i\lambda x - \omega(\lambda)t}\frac{2i\lambda}{i\lambda - \rho(\lambda)}\frac{d\lambda}{[i(i\rho(\lambda)]^5}\right\}$$

$$\approx u_0^{(4)}(0)\left\{\int_1^\infty e^{i\lambda x - \omega(\lambda)t}\frac{d\lambda}{i\lambda} + \int_1^\infty e^{i\lambda x - \omega(\lambda)t}\frac{i-1}{i+1}\frac{d\lambda}{(-i\lambda)} - \int_1^\infty e^{i\lambda x - \omega(\lambda)t}\frac{2i}{i+1}\frac{d\lambda}{\lambda}\right\} \approx -2(1+i)u_0^{(4)}(0)\int_1^\infty e^{i\lambda x - \omega(\lambda)t}\frac{d\lambda}{\lambda}, \quad (5.12)$$

where we used also the fact that, for $\lambda > 0$, $\rho(\lambda) = -\sqrt{\lambda^2 + \frac{\alpha}{\beta}} = -\lambda + O(1/\lambda)$, for $\lambda \to +\infty$.

Similarly,

$$\frac{\partial^4}{\partial x^4}\left\{\int_0^\infty e^{\rho(\mu)x - \omega(\mu)t}\left[\frac{i\mu - \rho(\mu)}{i\mu + \rho(\mu)}\hat{u}_0(i\rho(\mu)) + \frac{2\rho(\mu)}{i\mu + \rho(\mu)}\hat{u}_0(\mu)\right]\frac{\mu d\mu}{i\rho(\mu)}\right\} \approx (2+i)u_0^{(4)}(0)\int_1^\infty e^{\rho(\mu)x - \omega(\mu)t}\frac{d\mu}{\mu}, \quad (5.13)$$

$$\frac{\partial^4}{\partial x^4}\left\{\int_0^\infty e^{i\lambda x - \omega(\lambda)t}\widetilde{\tilde{f}}(\lambda, \omega(\lambda), t)d\lambda + \frac{1}{2\pi}\int_{-\infty}^\infty e^{i\lambda x - \omega(\lambda)t}\left[\frac{i\lambda + \rho(\lambda)}{i\lambda - \rho(\lambda)}\widetilde{\tilde{f}}(-\lambda, \omega(\lambda), t) - \frac{2i\lambda}{i\lambda - \rho(\lambda)}\widetilde{\tilde{f}}(i\rho(\lambda), \omega(\lambda), t))\right]d\lambda\right\}$$

$$\approx 2(1+i)\frac{f(0,0)}{\beta}\int_1^\infty e^{i\lambda x - \omega(\lambda)t}\frac{d\lambda}{\lambda}, \quad (5.14)$$

$$\int_0^\infty e^{\rho(\mu)x - \omega(\mu)t}\left[\frac{i\mu - \rho(\mu)}{i\mu + \rho(\mu)}\widetilde{\tilde{f}}(i\rho(\mu), \omega(\mu), t) + \frac{2\rho(\mu)}{i\mu + \rho(\mu)}\widetilde{\tilde{f}}(\mu, \omega(\mu), t))\right]\frac{\mu d\mu}{i\rho(\mu)}$$

$$\approx -(2+i)\frac{f(0,0)}{\beta}\int_1^\infty e^{\rho(\mu)x - \omega(\mu)t}\frac{d\mu}{\mu}. \quad (5.15)$$

Deriving also the analogues of (5.12) – (5.15) for the integrals $\int_{-\infty}^0 \cdots$ (keeping in mind that for $\lambda < 0$, $\rho(\lambda) = -\sqrt{\lambda^2 + \frac{\alpha}{\beta}} = \lambda + O(1/\lambda)$, for $\lambda \to -\infty$), we easily find that

$$2\pi\frac{\partial^4 U^*(x,t)}{\partial x^4} \approx \left[u_0^{(4)}(0) - \frac{f(0,0)}{\beta}\right]\left\{-2(1+i)\int_1^\infty e^{i\lambda x - \omega(\lambda)t}\frac{d\lambda}{\lambda} + (2+i)\int_1^\infty e^{\rho(\lambda)x - \omega(\lambda)t}\frac{d\lambda}{\lambda}\right.$$

$$\left. + 2(1-i)\int_{-\infty}^{-1} e^{i\lambda x - \omega(\lambda)t}\frac{d\lambda}{\lambda} + (-2+i)\int_{-\infty}^{-1} e^{\rho(\lambda)x - \omega(\lambda)t}\frac{d\lambda}{\lambda}\right\},$$

and the 5th conclusion of the theorem follows.

**Theorem 7** *Assuming*

$$u_0(0) = g_0(0) = 0, \ u_0'(0) = g_1(0) = 0, \ u_0''(0) = u_0'''(0) = 0, \ g_0'(0) = 0 \ and \ \beta u_0''''(0) = f(0,0), \quad (5.16)$$

*the solution $U(x,t)$ of (5.1), given by (5.6), is unique in the following sense: If $V(x,t)$ is $C^4$ in $\overline{Q} - \{(0,0)\}$, satisfies (5.1),*

$$\lim_{x \to \infty} V(x,t) = \lim_{x \to \infty} V_x(x,t) = 0, \ \sup_{x \geq 1}|V_{xx}(x,t)| < \infty, \ \sup_{x \geq 1}|V_{xxx}(x,t)| < \infty \ (\forall t > 0), \quad (5.17)$$

*and, for every $T > 0$, the functions $|V(x,t)|^2$, $|V_t(x,t)|^2$, $|V_x(x,t)|^2$, $|V_{xx}(x,t)|^2$, $|V_{xxx}(x,t)|^2$, are, uniformly for $0 < t \leq T$, integrable with respect to $x \in [0, \infty)$, i.e., there exists a positive function $B_T(x)$ such that $\int_0^\infty B_T(x)dx < +\infty$ and, for $0 < t \leq T$,*

$$|V(x,t)|^2 \leq B_T(x), \ |V_t(x,t)|^2 \leq B_T(x), \ |V_x(x,t)|^2 \leq B_T(x), \ |V_{xx}(x,t)|^2 \leq B_T(x), \ |V_{xxx}(x,t)|^2 \leq B_T(x) \ (x > 0), \quad (5.18)$$

*then $V \equiv U$.*

**Proof** We may assume that the data $u_0$, $g_0$, $g_1$ and $f$, and the functions $U$ and $V$ are real-valued. Let us set

$$W(x,t) := V(x,t) - U(x,t).$$

Then $W(x,t)$ is $C^4$ in $\overline{Q} - \{(0,0)\}$,

$$W_t = \alpha W_{xx}(x,t) - \beta W_{xxxx} \text{ for } (x,t) \in \overline{Q} - \{0,0\}, \ W(x,0) = 0 \text{ for } x > 0, \text{ and } W(0,t) = W_x(0,t) = 0 \text{ for } t > 0,$$

and (5.17) and (5.18) hold with $V$ replaced by $U$ (by taking $M_T$ larger – if necessary). Indeed, this last assertion for $U$ follows from the analogues – in this case – of Theorems 4 and 5, and Theorem 6, also in view of the assumption (5.16). Thus, (5.17) and (5.18) hold also with $V$ replaced by $W$ (by taking $M_T$ larger – if necessary).

Then the equation $W_t = \alpha W_{xx}(x,t) - \beta W_{xxxx}$ implies that

$$\int_{x=0}^{\infty} W(x,t) W_t(x,t) dx = \int_{x=0}^{\infty} W(x,t)[\alpha W_{xx}(x,t) - \beta W_{xxxx}] dx \text{, for every } t > 0 . \tag{5.19}$$

We note that, since (5.18) holds with $V$ replaced by $W$, the integrals in (5.19) are absolutely convergent. This also implies that, by Lebesgue's dominated convergence theorem,

$$\int_{x=0}^{\infty} W(x,t) W_t(x,t) dx = \frac{1}{2} \frac{d}{dt}\left[\int_{x=0}^{\infty}[W(x,t)]^2 dx\right] \text{, for every } t > 0 . \tag{5.20}$$

Integrating by parts, we obtain that, for $t > 0$ and $A > 0$,

$$\int_{x=0}^{A} W(x,t)[\alpha W_{xx}(x,t) - \beta W_{xxxx}] dx = \alpha [W(x,t) W_x(x,t)]\big|_{x=0}^{x=A} - \alpha \int_{x=0}^{A}[W_x(x,t)]^2 dx$$

$$- \beta [W(x,t) W_{xxx}(x,t)]\big|_{x=0}^{x=A} + \beta [W_x(x,t) W_{xx}(x,t)]\big|_{x=0}^{x=A} - \beta \int_{x=0}^{A}[W_{xx}(x,t)]^2 dx . \tag{5.21}$$

Letting $A \to \infty$ and taking into consideration that (5.17) holds with $V$ replaced by $W$, we see that (5.21) implies that

$$\int_{x=0}^{\infty} W(x,t)[\alpha W_{xx}(x,t) - \beta W_{xxxx}] dx = -\alpha \int_{x=0}^{\infty}[W_x(x,t)]^2 dx - \beta \int_{x=0}^{\infty}[W_{xx}(x,t)]^2 dx . \tag{5.22}$$

It follows from (5.19), (5.20) and (5.22), that

$$\frac{d}{dt}\left[\int_{x=0}^{\infty}[W(x,t)]^2 dx\right] \leq 0 \text{, for every } t > 0 ,$$

and, therefore,

$$\int_{x=0}^{\infty}[W(x,s)]^2 dx \leq \int_{x=0}^{\infty}[W(x,t)]^2 dx \text{, for every } s > t > 0 . \tag{5.23}$$

Since

$$|V(x,t)|^2 \leq M_T(x), \ 0 < t \leq T \text{, and } \int_0^{\infty} M_T(x) dx < +\infty ,$$

Lebesgue's dominated convergence theorem implies that

$$\lim_{t \to 0^+} \int_{x=0}^{\infty}[W(x,t)]^2 dx = \int_{x=0}^{\infty}[W(x,0)]^2 dx = 0 . \tag{5.24}$$

(The last equation in (5.24) follows from the fact that $W(x,0) = 0$, for $x > 0$.)

Thus, letting $t \to 0^+$ in (5.23), we obtain, in view of (5.24), that

$$\int_{x=0}^{\infty}[W(x,s)]^2 dx \leq 0 \ (\forall s > 0) \implies W(x,s) \equiv 0 .$$

This completes the proof.

***Note*** The above theorem holds for $\alpha \geq 0$, i.e., it includes also the case $\alpha = 0$.

***Example of non-uniqueness*** Let us consider the solution $U(x,t)$ of (5.1), given by (5.6), with data: $u_0(x) \equiv 0$, $g_0(t) \equiv 1$, $g_1(t) \equiv 0$ and $f(x,t) \equiv 0$, i.e.,

$$U(x,t) = \frac{\beta i}{\pi} \int_{\partial\Omega} [e^{i\lambda x - \omega(\lambda)t} - e^{i\lambda x}] \frac{\lambda[i\lambda + \rho(\lambda)]\rho(\lambda)}{\omega(\lambda)} d\lambda \text{, for } x > 0, \ t > 0 .$$

By the analogue – in this case – of Theorem 5, $U(x,t)$ is $C^{\infty}$ in $\overline{Q} - \{(0,0)\}$, $\partial_t U = \alpha \partial_x^2 U - \beta \partial_x^4 U$ (for $(x,t) \in Q$) and

$$U(x,0) = 0, \ U(0,t) = 1 \text{ and } U_x(0,t) = 0 .$$

By Cauchy's theorem and Jordan's lemma,

$$U(x,t) = \frac{\beta i}{\pi} \int_{\Gamma} \left[ e^{i\lambda x - \omega(\lambda)t} - e^{i\lambda x} \right] \frac{\lambda[i\lambda + \rho(\lambda)]\rho(\lambda)}{\omega(\lambda)} d\lambda = \frac{\beta i}{\pi} \int_{\Gamma} e^{i\lambda x - \omega(\lambda)t} \frac{\lambda[i\lambda + \rho(\lambda)]\rho(\lambda)}{\omega(\lambda)} d\lambda ,$$

where $\Gamma = \{\lambda \in \partial\Omega : |\lambda| \geq r\} + \{\lambda \in \Omega : |\lambda| = r\}$ with $r > \sqrt{\alpha/\beta}$ .

Let us also consider the $t$–derivative of $U(x,t)$ ,

$$v(x,t) := \frac{\partial U(x,t)}{\partial t} = \frac{\beta i}{\pi} \frac{\partial}{\partial t} \left[ \int_{\Gamma} e^{i\lambda x - \omega(\lambda)t} \frac{\lambda[i\lambda + \rho(\lambda)]\rho(\lambda)}{\omega(\lambda)} d\lambda \right] = -\frac{\beta i}{\pi} \int_{\Gamma} e^{i\lambda x - (\alpha\lambda^2 + \beta\lambda^4)t} \lambda[i\lambda + \rho(\lambda)]\rho(\lambda) d\lambda .$$

Again, by the analogue – in this case – of Theorem 5, $v(x,t)$ is $C^{\infty}$ in $\overline{Q} - \{(0,0)\}$ , $\partial_t v = \alpha\partial_{xx}^2 v - \beta\partial_{xxxx}^4 v$ (for $(x,t) \in Q$ ) and

$$v(0,t) = 0 \text{ and } v_x(0,t) = 0 .$$

Also,

$$v(x,0) = -\frac{\beta i}{\pi} \int_{\Gamma} e^{i\lambda x} \lambda[i\lambda + \rho(\lambda)]\rho(\lambda) d\lambda = 0 ,$$

since

$$\lim_{R \to \infty} \int_{\Omega \cap \{|\lambda| = R\}} e^{i\lambda x} \lambda[i\lambda + \rho(\lambda)]\rho(\lambda) d\lambda = 0 .$$

We claim that

$$v(x,t) \not\equiv 0 .$$

Indeed, to reach a contradiction, let us suppose that $v(x,t) \equiv 0$ . Then $\frac{\partial U(x,t)}{\partial t} \equiv 0$ , which would imply that $U(x,t) \equiv U(x,0) \equiv 0$ . But this contradicts the equation $U(0,t) \equiv 1$ .

Thus, the IBVP (5.1) with data $u_0(x) \equiv 0$ , $g_0(x) \equiv 0$ , $g_1(x) \equiv 0$ and $f(x,t) \equiv 0$ , does not have a unique solution – in fact it has an infinity of solutions.

**The Fokas method solution as a canonical solution** The solution of (5.1), given by (5.6), satisfies the conditions of Theorem 7, namely:

$$V(x,t) \in C^4(\overline{Q} - \{(0,0)\}) , \text{(5.17) and (5.18)}, \tag{*}$$

provided that the data satisfy (5.16), and, therefore, it is the unique solution among the solutions which satisfy (*). In this sense, the Fokas method solution (5.6) is a *canonical* solution of (5.1), when the data satisfy, in addition, (5.16).

**Theorem 8** (*lack of null controllability*) *Suppose that* $u_0$ *and* $f$ *satisfy* $u_0(0) = u_0'(0) = u_0''(0) = u_0'''(0) = 0$ , $\beta u_0''''(0) = f(0,0)$ , *and that* $\hat{u}_0(\lambda)$ *does not extend to an analytic function of* $\lambda \in \mathbb{C}$ *(in view of [89], the latter holds generically). Then for every* $T > 0$ , *so that* $\widetilde{\widetilde{f}}(\lambda, \omega(\lambda), T)$ *extends to an entire function, i.e. analytic in* $\lambda \in \mathbb{C}$ , *and for every* $g_0$ , $g_1$ , *with* $g_0(0) = g_1(0) = g_0'(0) = 0$ , *we have that the solution of (5.1), given by (5.6), satisfies* $U(x,T) \not\equiv 0$ .
*In particular, the homogenous version of problem (5.1) (i.e., with* $f \equiv 0$ *) is not null controllable.*

**Proof** Let $T$ , $g_0$ , $g_1$ and $f$ be as in the statement of the theorem. Suppose – to reach a contradiction – that $U(x,T) \equiv 0$ . Then $\hat{U}(\lambda, T) = 0$ , and (5.2) becomes

$$\frac{1}{\beta} \hat{u}_0(\lambda) = -[\widetilde{g}_3(\omega(\lambda), T) + i\lambda\widetilde{g}_2(\omega(\lambda), T)]$$

$$+ \left(\lambda^2 + \frac{\alpha}{\beta}\right)\widetilde{g}_1(\omega(\lambda), T) + i\lambda\left(\lambda^2 + \frac{\alpha}{\beta}\right)\widetilde{g}_0(\omega(\lambda), T) - \frac{1}{\beta}\widetilde{\widetilde{f}}(\lambda, \omega(\lambda), T) = 0 , \text{ Im} \lambda \leq 0 . \tag{5.25}$$

We emphasize that it is because of the assumptions on the data that the derivation of (5.25) is rigorous. Indeed, this follows from Theorems 5 and 6, which guarantee the behavior of the solution $U(x,t)$ on the boundary of $Q$ , i.e., when $x = 0$ or $t = 0$ , including the point $(0,0)$ , and when $x \to \infty$ . In particular the existence of the limits

$$\lim_{\substack{(x,t)\to(0,0)\\(x,t)\in\widetilde{Q}-\{(0,0)\}}} \frac{\partial^k U^*(x,t)}{\partial x^k}, \text{ for } k=2 \text{ and } 3,$$

as this is described in Theorem 6, imply that the functions

$$\widetilde{g}_k(\omega(\lambda),t) = \int_{\tau=0}^{t} e^{\omega(\lambda)\tau} \left[ \frac{\partial^k U(x,\tau)}{\partial x^k} \bigg|_{x=0} \right] d\tau, \text{ for } k=2 \text{ and } 3,$$

are entire functions of $\lambda$ (for every $t$). This holds also for $k=0,\ 1$, and for the function $\widetilde{f}(\lambda,\omega(\lambda),T)$.

Therefore, (5.2) implies that the function $\hat{u}_0(\lambda)$ extends to an analytic function of $\lambda\in\mathbb{C}$. This contradicts our assumption and proves the conclusion of the theorem.

**Remark** Examples of functions $u_0$, for which $\hat{u}_0(\lambda)$ does not extend to an entire function of $\lambda$, are those with $u_0(x) = e^{-Ax}$ for $x \geq x_0$, for some $A > 0$ and $x_0 > 0$. Indeed,

$$\hat{u}_0(\lambda) = \frac{1}{i\lambda + A} + an\ entire\ function.$$

Examples of functions $f$ such that $\widetilde{f}(\lambda,\omega(\lambda),T)$ extends to an entire function of $\lambda$, are those which can be written as $f(x,t) = e^{-Ax^2}\phi(t)$ for $x \geq x_0$ and $t \geq t_0$, for some $A > 0$ and $x_0, t_0 > 0$.

**Comment** (restatement of Theorem 8) With the notation and the assumptions as in Theorem 8, $V(x,T) \not\equiv 0$, where $V$ is assumed to be as in Theorem 7. Indeed, this follows immediately from Theorem 8, since – in view of Theorem 7 – we have $V \equiv U$.

## 6. The equation $\partial_t U = -\alpha\partial_x^2 U - \beta\partial_x^4 U + f$

In light of all the above expository analysis, we are finally ready to formulate our novel closed-form solution of the celebrated C-H eq. with arbitrary forcing, initial and boundary data.

**Problem** With $\alpha > 0$ and $\beta > 0$, solve

$$\begin{cases} \partial_t U = -\alpha\partial_x^2 U - \beta\partial_x^4 U + f, \ (x,t)\in Q := \mathbb{R}^+\times\mathbb{R}^+, \\ U(x,0) = u_0(x), \ x\in\mathbb{R}^+, \\ U(0,t) = g_0(t), \ t\in\mathbb{R}^+, \\ U_x(0,t) = g_1(t), \ t\in\mathbb{R}^+. \end{cases} \tag{6.1}$$

**Solution** The solution of (6.1) is given by formula (5.6), **but** with the following notation for $\omega$, $\Omega$ and $\rho$:

$$\omega = \omega_{\alpha,\beta} = -\alpha\lambda^2 + \beta\lambda^4,$$

$$\Omega = \Omega_{\alpha,\beta} = \{\lambda : \text{Im}\,\lambda \geq 0\ \&\ \text{Re}\,\omega_{\alpha,\beta}(\lambda) = \text{Re}(-\alpha\lambda^2 + \beta\lambda^4) \leq 0\}, \tag{6.2}$$

and

$$\rho(\lambda) = \rho_{\alpha,\beta}(\lambda) = \begin{cases} -\sqrt{\left|\lambda^2 - \dfrac{\alpha}{\beta}\right|}\exp\{i[\theta_1(\lambda)+\theta_2(\lambda)+\pi]/2\} & for\ \text{Re}\,\lambda > 0 \\ \sqrt{\left|\lambda^2 - \dfrac{\alpha}{\beta}\right|}\exp\{i[\theta_1(\lambda)+\theta_2(\lambda)+\pi]/2\} & for\ \text{Re}\,\lambda < 0. \end{cases} \tag{6.3}$$

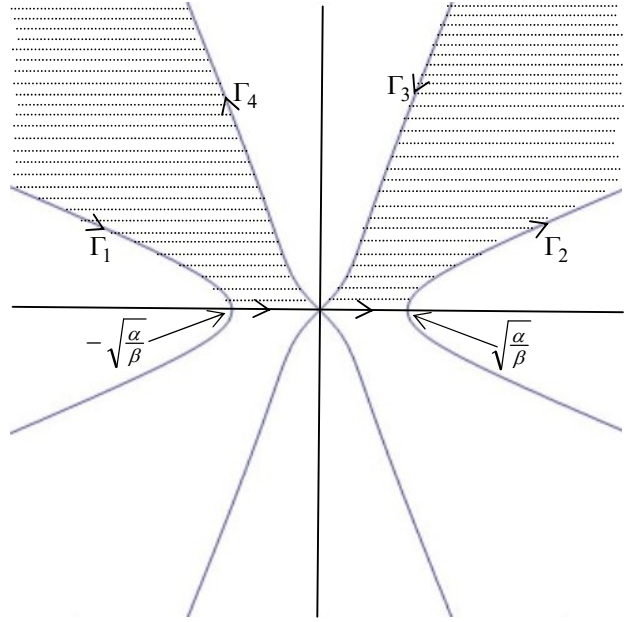

**Fig. 10** The set $\Omega_{\alpha,\beta}$, defined by (6.2).

The boundary of $\Omega_{\alpha,\beta}$ is

$$\partial\Omega_{\alpha,\beta} = \left[\Gamma_1 \cup \left[-\sqrt{\tfrac{\alpha}{\beta}}, 0\right] \cup \Gamma_4\right] \cup \left[\Gamma_2 \cup \left[0, \sqrt{\tfrac{\alpha}{\beta}}\right] \cup \Gamma_3\right] \text{ (see fig. 10)}.$$

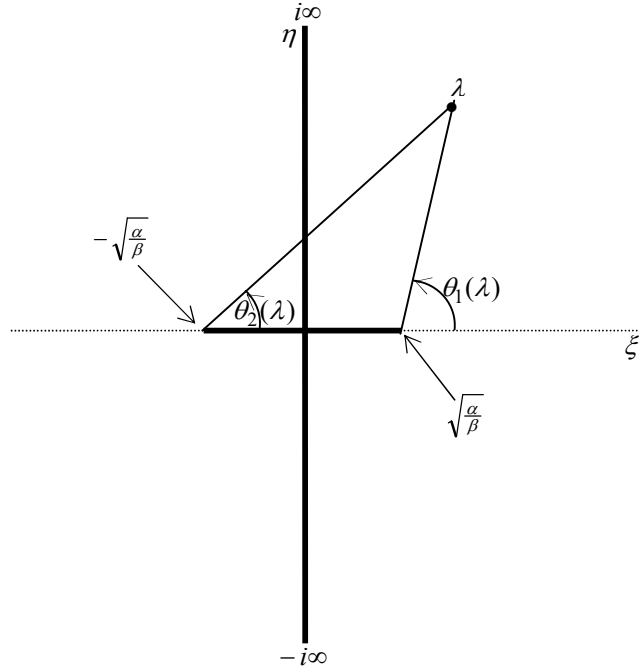

**Fig. 11** With $-\pi < \theta_1(\lambda) < \pi$ and $-\pi < \theta_2(\lambda) < \pi$, the function $\rho(\lambda) = \rho_{\alpha,\beta}(\lambda)$ is defined by (6.3).

The function $\rho(\lambda) = \rho_{\alpha,\beta}(\lambda)$, defined by (6.3), is analytic in $\Theta_{\alpha,\beta} := \mathbb{C} - \{[-\sqrt{\tfrac{\alpha}{\beta}}, \sqrt{\tfrac{\alpha}{\beta}}] \cup i\mathbb{R}\}$ (see fig. 11) and the restriction $\rho|_{\Theta_{\alpha,\beta} \cap \{\operatorname{Im}\lambda > 0\}}$ extends continuously to $[\Theta_{\alpha,\beta} \cap \{\operatorname{Im}\lambda > 0\}] \cup (\mathbb{R} - \{0\})$.

The rigorous justification should by now be straightforward for the readers and is hence omitted for the sake of economy of space and avoiding repetitions.

Finally, theorems analogous to Theorems 5 and 6 can be proved in this case too. However, in this case, we do not have a uniqueness theorem analogous to Theorem 7, due to the sign of one of the parameters which prevents applicability of the energy estimate argument.

## 7. Eventual periodicity of the solution

**Theorem 9** *Let* $T > 0$ *be fixed. Suppose that, in addition to (1.2),* $g_0(t)$, $g_1(t)$ *and* $f(x,t)$ *are* $T$−*periodic, i.e.,* $g_k(t+T) = g_k(t)$ $(k = 0,1)$ *and* $f(x,t) = f(x,t+T)$ *for every* $t \geq 0$ *and* $x \geq 0$*. Then, the solution* $U(x,t)$ *of problem (5.1), defined by (5.6), is eventually* $T$−*periodic, i.e., for every fixed* $x_0 > 0$,

$$\lim_{t \to \infty}[U(x_0, t+T) - U(x_0, t)] = 0 . \tag{7.1}$$

*More precisely,*

$$U(x_0, t+T) - U(x_0, t) = \mathfrak{C}_1 \frac{1}{t^{1/2}} + \mathfrak{C}_2 \frac{1}{t^{3/2}} + O(1/t^{5/2}), \ as \ t \to \infty, \tag{7.2}$$

*where* $\mathfrak{C}_1$, $\mathfrak{C}_2$ *are constants which are explicitly expressed in terms of the data.*

**Proof** We start with some remarks. Firstly, with notation as in Section 5 and by the periodicity assumptions, we have

$$e^{-\omega(\lambda)(t+T)}\widetilde{g}_k(\omega(\lambda), t+T) - e^{-\omega(\lambda)t}\widetilde{g}_k(\omega(\lambda), t)$$

$$= e^{-\omega(\lambda)(t+T)}\int_{\tau=0}^{t+T} e^{\omega(\lambda)\tau} g_k(\tau) d\tau - e^{-\omega(\lambda)t}\int_{\tau=0}^{t} e^{\omega(\lambda)\tau} g_k(\tau) d\tau = e^{-\omega(\lambda)t}\int_{\tau=-T}^{0} e^{\omega(\lambda)\tau} g_k(\tau) d\tau, \tag{7.3}$$

$\hat{f}(\lambda, t+T) = \hat{f}(\lambda, t)$ and

$$e^{-\omega(t+T)}\widetilde{\hat{f}}(\lambda, \omega(\lambda), t+T) - e^{-\omega t}\widetilde{\hat{f}}(\lambda, \omega(\lambda), t) = e^{-\omega t}\int_{\tau=-T}^{0} e^{\omega(\lambda)\tau} \hat{f}(\lambda, \tau) d\tau . \tag{7.4}$$

Secondly, we point out that the contours of the integrals in (5.6) can be deformed, in such a way that the values of the integrals remain unchanged. For example, with appropriate integrands,

$$\int_{\partial\Omega} = \int_{E_1} + \int_{E_2} + \int_{E_3} + \int_{E_4} , \tag{7.5}$$

where $E_1$, $E_2$, $E_3$, $E_4$ are half-lines as the ones depicted in Fig.12. (See also Fig.7.)

Also, with appropriate integrands,

$$\int_{E_4} + \int_{E_1} = \int_{E_4 \cap \{|\lambda| \geq |A|\}} + \int_{[A,-1]} + \int_{-1}^{1} + \int_{[1,B]} + \int_{E_1 \cap \{|\lambda| \geq |B|\}} = \int_{Z} + \int_{-1}^{1} \tag{7.6}$$

where

$A := E_4 \cap \{\text{Re}\,\lambda = -1\}$, $B := E_1 \cap \{\text{Re}\,\lambda = 1\}$ and $Z := (E_4 \cap \{|\lambda| \geq |A|\}) + [A,-1] + [1,B] + (E_1 \cap \{|\lambda| \geq |B|\})$.

(See Fig.13.)

Thirdly, with the notation as in the proof of (5.7), the integrals

$$\int_{C^{\pm}_{upper}}\mathcal{I}_{u_0} , \ \int_{C^{\pm}_{upper}}\mathcal{I}_{g_0, g_1} \ and \ \int_{C^{\pm}_{upper}}\mathcal{I}_{f} , \tag{7.7}$$

can be tranformed to integrals over the contours $C^{\pm}_{lower}$, by appropriate change of variable, and the contours $C^{\pm}_{lower}$ can, in turn, be deformed to the half-lines $E_4$ and $E_1$.

Keeping in mind the above remarks, we proceed with the proof.

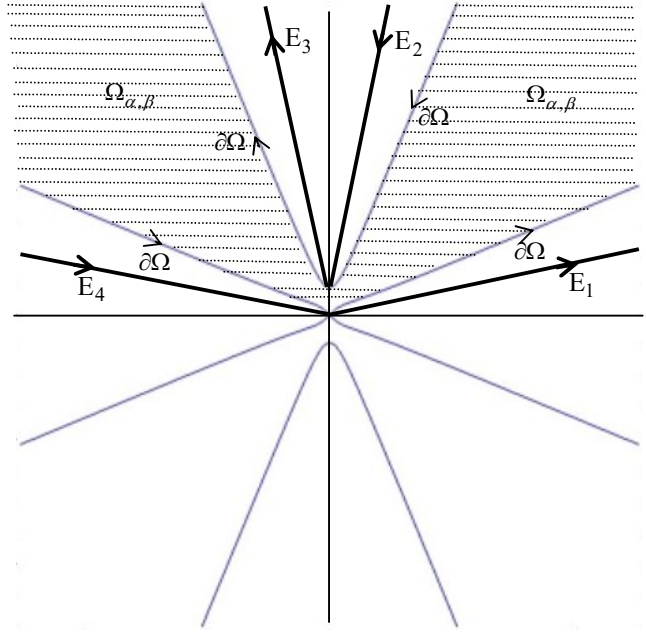

**Fig.12** The half-lines $\mathrm{E}_j$ are chosen so that $\mathrm{Im}\,\lambda > 0$ and $\mathrm{Re}\,\omega(\lambda) > 0$, for $\lambda \in \mathrm{E}_j - \{0\}$.

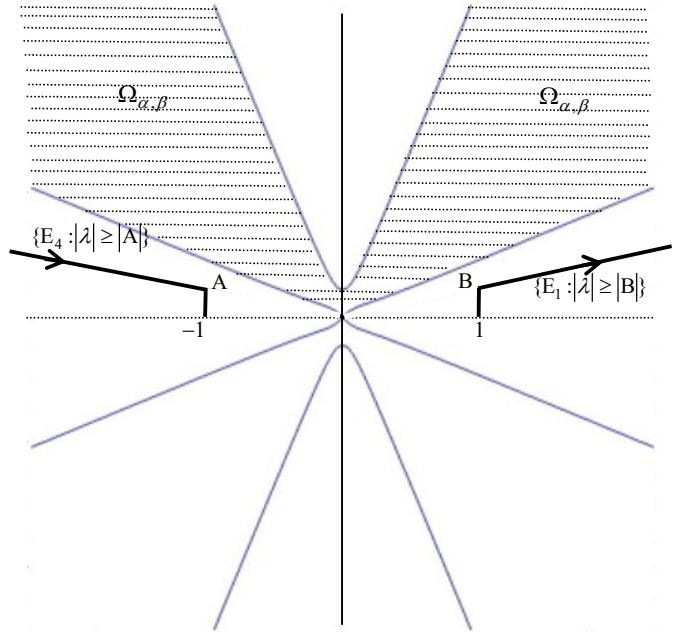

**Fig.13** The contour $Z$.

***Step 1*** We claim that

$$\lim_{t \to \infty} \int_{\partial\Omega} \lambda[i\lambda + \rho(\lambda)]\rho(\lambda)e^{i\lambda x}[e^{-\omega(\lambda)(t+T)}\widetilde{g}_0(\omega(\lambda), t+T) - e^{-\omega(\lambda)t}\widetilde{g}_0(\omega(\lambda), t)]d\lambda = 0. \qquad (7.8)$$

Indeed, in view of (7.3), (7.5) is written as follows:

$$\lim_{t \to \infty} \int_{\partial\Omega} \lambda[i\lambda + \rho(\lambda)]\rho(\lambda)e^{i\lambda x}e^{-\omega(\lambda)t} \int_{\tau=-T}^{0} e^{\omega(\lambda)\tau} g_0(\tau)d\tau d\lambda = 0. \qquad (7.9)$$

But for fixed $\lambda \in \mathrm{E}_j$,

$$\lim_{t \to \infty}\left\{\lambda[i\lambda + \rho(\lambda)]\rho(\lambda)e^{i\lambda x}e^{-\omega(\lambda)t} \int_{\tau=-T}^{0} e^{\omega(\lambda)\tau} g_0(\tau)d\tau\right\} = 0,$$

since $\operatorname{Re}\omega(\lambda) > 0$ for $\lambda \in \mathrm{E}_j - \{0\}$.

Also, for some $\varepsilon > 0$,

$$\left| \lambda[i\lambda + \rho(\lambda)]\rho(\lambda)e^{i\lambda x}e^{-\omega(\lambda)t}\int_{\tau=-T}^{0}e^{\omega(\lambda)\tau}g_0(\tau)d\tau \right| \leq e^{-\varepsilon t|\lambda|}\left|\lambda[i\lambda + \rho(\lambda)]\rho(\lambda)\right|\int_{\tau=-T}^{0}|g_0(\tau)|d\tau, \text{ for } \lambda \in \mathrm{E}_j.$$

Thus, in view of (7.5), (7.8) follows from Lebesgue's dominated convergence theorem.

***Step 2*** Similarly, an equation analogous to (7.8) can be obtained also for the part of $U(x_0, t+T) - U(x_0, t)$ which involves the function $g_1$.

***Step 3*** We claim that

$$\lim_{t\to\infty}\int_{-\infty}^{\infty}e^{i\lambda x - \omega(\lambda)t}\hat{u}_0(\lambda)d\lambda = 0, \text{ and, therefore, } \lim_{t\to\infty}\left[\int_{-\infty}^{\infty}e^{i\lambda x - \omega(\lambda)(t+T)}\hat{u}_0(\lambda)d\lambda - \int_{-\infty}^{\infty}e^{i\lambda x - \omega(\lambda)t}\hat{u}_0(\lambda)d\lambda\right] = 0. \quad (7.10)$$

This follows from Lebesgue's dominated convergence theorem, since

$$\lim_{t\to\infty}e^{i\lambda x - \omega(\lambda)t}\hat{u}_0(\lambda) = 0, \text{ for every fixed } \lambda \in \mathbb{R} - \{0\}, \text{ and } \sup_{\lambda\in\mathbb{R}}\left|e^{i\lambda x - \omega(\lambda)t}\hat{u}_0(\lambda)\right| \preceq e^{-\alpha\lambda^2}, \ \forall t.$$

***Step 4*** Similarly, using also the analogue of (7.5) for the integral over $\partial\Omega$, equations analogous to (7.10) can be obtained also for the parts of $U(x_0, t+T) - U(x_0, t)$ which involve the quantities $\hat{u}_0(-\lambda)$ and $\hat{u}_0(i\rho(\lambda))$.

***Step 5*** Since

$$\sup\left\{\left|\int_{\tau=-T}^{0}e^{\omega(\lambda)\tau}\hat{f}(\lambda,\tau)d\tau\right|: \ \lambda\in\mathbb{R}\right\} < \infty,$$

it is easy to see that

$$\lim_{t\to\infty}\int_{-\infty}^{\infty}e^{i\lambda x - \omega(\lambda)t}\int_{\tau=-T}^{0}e^{\omega(\lambda)\tau}\hat{f}(\lambda,\tau)d\tau d\lambda = 0.$$

Similarly, using also the analogue of (7.5) for the integral over $\partial\Omega$,

$$\lim_{t\to\infty}\int_{\partial\Omega}e^{i\lambda x - \omega(\lambda)t}\left[\frac{i\lambda + \rho(\lambda)}{i\lambda - \rho(\lambda)}\int_{\tau=-T}^{0}e^{\omega(\lambda)\tau}\hat{f}(-\lambda,\tau)d\tau - \frac{2i\lambda}{i\lambda - \rho(\lambda)}\int_{\tau=-T}^{0}e^{\omega(\lambda)\tau}\hat{f}(i\rho(\lambda),\tau)d\tau\right]d\lambda = 0.$$

Thus, the part of $U(x, t+T) - U(x, t)$ which involves the quantities $\hat{f}(\lambda,\tau)$, $\hat{f}(-\lambda,\tau)$, $\hat{f}(i\rho(\lambda),\tau)$, tends to zero, as $t\to\infty$, also in view of (7.4).

*Proof of (7.1)* It follows from the conclusions of *Steps 1-5*.

***Step 6*** By (7.3),

$$\int_{\partial\Omega}[\mathcal{I}_{g_0,g_1}(x_0, t+T) - \mathcal{I}_{g_0,g_1}(x_0, t)] = \int_{\partial\Omega}\mathcal{J}_{g_0,g_1}(x_0, t),$$

where

$$\mathcal{J}_{g_0,g_1}(x_0, t) := -\frac{\beta i}{\pi}e^{i\lambda x_0 - \omega(\lambda)t}\lambda[i\lambda + \rho(\lambda)]\left\{\int_{\tau=-T}^{0}e^{\omega(\lambda)\tau}[\rho(\lambda)g_0(\tau) - g_1(\tau)]d\tau\right\}d\lambda.$$

Also, using (7.6),

$$\left(\int_{\mathrm{E}_4} + \int_{\mathrm{E}_1}\right)\mathcal{J}_{g_0,g_1}(x_0, t) = \int_{\mathrm{Z}}\mathcal{J}_{g_0,g_1}(x_0, t) + \int_{-1}^{1}\mathcal{J}_{g_0,g_1}(x_0, t). \quad (7.11)$$

Writing

$$\mathcal{J}_{g_0,g_1}(x_0, t) = e^{i\lambda x_0 - \omega(\lambda)t}\cdots d\lambda = \frac{1}{t}\left(\frac{d}{d\lambda}e^{-\omega(\lambda)t}\right)\left\{\frac{1}{[-\omega(\lambda)]'}e^{i\lambda x_0}\cdots\right\}d\lambda$$

and integrating by parts, repeatedly, we see that

$$\int_{\mathrm{Z}}\mathcal{J}_{g_0,g_1}(x_0, t) = \mathrm{O}(1/t^{\ell}), \text{ as } t\to\infty, \ \forall\ell. \quad (7.12)$$

On the other hand,

$$\int_{-1}^{1} \mathcal{J}_{g_0, g_1}(x_0, t) = \mathfrak{C}_1' \frac{1}{t^{1/2}} + \mathfrak{C}_2' \frac{1}{t^{3/2}} + O(1/t^{5/2}), \text{ as } t \to \infty, \tag{7.13}$$

for some constants $\mathfrak{C}_1'$, $\mathfrak{C}_2'$, which can be explicitly computed. Indeed, this follows from the following version of Laplace's asymptotic expansion formula: If $I \subset \mathbb{R}$ is a bounded open interval, $\lambda_0 \in I$, $\phi'(\lambda_0) = 0$, $\phi'(\lambda) \neq 0$ for $\lambda \in \bar{I} - \{\lambda_0\}$, and $\phi''(\lambda_0) < 0$ then, as $t \to \infty$,

$$\int_I e^{t\phi(\lambda)} h(\lambda) d\lambda = e^{t\phi(\lambda_0)} \frac{\sqrt{2\pi}}{\sqrt{-\phi''(\lambda_0)}} \Bigg\{ h(\lambda_0) \frac{1}{t^{1/2}}$$

$$+ \Bigg[ -\frac{h^{(2)}(\lambda_0)}{2\phi^{(2)}(\lambda_0)} + \frac{h(\lambda_0)\phi^{(4)}(\lambda_0)}{8[\phi^{(2)}(\lambda_0)]^2} + \frac{h^{(1)}(\lambda_0)\phi^{(3)}(\lambda_0)}{2[\phi^{(2)}(\lambda_0)]^2} - \frac{5h(\lambda_0)[\phi^{(3)}(\lambda_0)]^2}{24[\phi^{(2)}(\lambda_0)]^3} \Bigg] \frac{1}{t^{3/2}} + O(1/t^{5/2}) \Bigg\}.$$

The functions $\phi$ and $h$ are assumed to be sufficiently smooth in $\bar{I} = closure(I)$.

Combining (7.11), (7.12) and (7.13), we obtain

$$\left( \int_{E_4} + \int_{E_1} \right) [\mathcal{I}_{g_0, g_1}(x_0, t+T) - \mathcal{I}_{g_0, g_1}(x_0, t)] = \mathfrak{C}_1' \frac{1}{t^{1/2}} + \mathfrak{C}_2' \frac{1}{t^{3/2}} + O(1/t^{5/2}), \text{ as } t \to \infty. \tag{7.14}$$

*Proof of (7.2)* This can be carried out by deriving equations analogous to (7.14) for all the integrals in (5.6).

**Theorem 10** *Let $T > 0$ be fixed. Suppose that, in addition to (1.2), $g_0(t)$, $g_1(t)$ and $f(x,t)$ are weakly $T-$ periodic in the following sense:*

$$\int_0^{\infty} |g_k(\tau + T) - g_k(\tau)| d\tau < +\infty, \ k = 0, 1, \ (7.15), \qquad \int_0^{\infty} |f(0, \tau + T) - f(0, \tau)| d\tau < +\infty, \ (7.16)$$

$$\int_0^{\infty} |f_x(0, \tau + T) - f_x(0, \tau)| d\tau < +\infty, \ (7.17) \quad and \quad \int_0^{\infty}\int_0^{\infty} |f(x, \tau + T) - f(x, \tau)| d\tau dy < +\infty, \ (7.18)$$

$$\int_0^{\infty}\int_0^{\infty} |f_{xx}(x, \tau + T) - f_{xx}(x, \tau)| d\tau dy < +\infty. \ (7.19)$$

*Then, the solution $U(x,t)$ of problem (5.1), defined by (5.6), is eventually $T-$ periodic, i.e., for every fixed $x_0 > 0$,*

$$\lim_{t \to \infty} [U(x_0, t+T) - U(x_0, t)] = 0.$$

**Proof** Since

$$e^{-\omega(t+T)} \tilde{g}_0(\omega(\lambda), t+T) - e^{-\omega t} \tilde{g}_0(\omega(\lambda), t) = e^{-\omega(t+T)} \int_0^{t+T} e^{\omega(\lambda)\tau} g_0(\tau) d\tau - e^{-\omega t} \int_0^t e^{\omega(\lambda)\tau} g_0(\tau) d\tau$$

$$= e^{-\omega t} \int_0^t e^{\omega(\lambda)\tau} [g_0(\tau + T) - g_0(\tau)] d\tau + e^{-\omega t} \int_{-T}^0 e^{\omega(\lambda)\tau} g_0(\tau + T) d\tau,$$

we have

$$\int_{\partial\Omega} e^{i\lambda x_0 - \omega(t+T)} \lambda[i\lambda + \rho(\lambda)] \rho(\lambda) \tilde{g}_0(\omega, t+T) d\lambda - \int_{\partial\Omega} e^{i\lambda x - \omega t} \lambda[i\lambda + \rho(\lambda)] \rho(\lambda) \tilde{g}_0(\omega, t) d\lambda$$

$$= \int_{\partial\Omega} e^{i\lambda x - \omega(\lambda)t} \lambda[i\lambda + \rho(\lambda)] \rho(\lambda) \left\{ \int_0^t e^{\omega(\lambda)\tau} [g_0(\tau + T) - g_0(\tau)] d\tau \right\} d\lambda$$

$$+ \int_{\partial\Omega} e^{i\lambda x - \omega(\lambda)t} \left[ \lambda[i\lambda + \rho(\lambda)] \rho(\lambda) \int_{-T}^0 e^{\omega(\lambda)\tau} g_0(\tau + T) d\tau \right] d\lambda. \tag{7.20}$$

The last integral in (7.20) tends to zero, as $t \to \infty$, as in the proof of Theorem 9.

Next we claim that, with the contours $E_j$ as in (7.5),

$$\lim_{t\to\infty}\int_{E_j} e^{i\lambda x-\omega(\lambda)t}\lambda[i\lambda+\rho(\lambda)]\rho(\lambda)\left\{\int_0^t e^{\omega(\lambda)\tau}[g_0(\tau+T)-g_0(\tau)]d\tau\right\}d\lambda=0\,. \tag{7.21}$$

First we show that, for every (fixed) $\lambda\in E_j$ ($\lambda\neq 0$),

$$\lim_{t\to\infty}\int_0^t e^{-\omega(\lambda)(t-\tau)}[g_0(\tau+T)-g_0(\tau)]d\tau=0\,. \tag{7.22}$$

Let us write the above integral in the form

$$\int_0^\infty \chi_{[0,t]}(\tau)e^{-\omega(\lambda)(t-\tau)}[g_0(\tau+T)-g_0(\tau)]d\tau\,,$$

where $\chi_{[0,t]}$ is the characteristic function of the interval $[0,t]$. For any fixed $\tau$,

$$\lim_{t\to\infty}\left\{\chi_{[0,t]}(\tau)e^{-\omega(\lambda)(t-\tau)}[g_0(\tau+T)-g_0(\tau)]\right\}=0\quad(\lambda\in E_j,\ \lambda\neq 0)\,.$$

Furthermore,

$$\left|\chi_{[0,t]}(\tau)e^{-\omega(\lambda)(t-\tau)}[g_0(\tau+T)-g_0(\tau)]\right|\leq\left|g_0(\tau+T)-g_0(\tau)\right|,\text{ for every }\tau\text{ and }t\,.$$

Thus, in view of (7.15) and Lebesgue's dominated convergence theorem, (7.22) follows.

Also, there exists $\varepsilon>0$ so that for $\lambda\in E_j$,

$$\left|e^{i\lambda x-\omega(\lambda)t}\lambda[i\lambda+\rho(\lambda)]\rho(\lambda)\left\{\int_0^t e^{\omega(\lambda)\tau}[g_0(\tau+T)-g_0(\tau)]d\tau\right\}\right|\leq e^{-\varepsilon|\lambda|}\left|\lambda[i\lambda+\rho(\lambda)]\rho(\lambda)\right|\int_0^\infty\left|g_0(\tau+T)-g_0(\tau)\right|d\tau\,.$$

Therefore, again by Lebesgue's dominated convergence theorem, (7.22) implies (7.21).

Next, we have

$$e^{-\omega(t+T)}\widetilde{\widehat{f}}(\lambda,\omega(\lambda),t+T)-e^{-\omega t}\widetilde{\widehat{f}}(\lambda,\omega(\lambda),t)=e^{-\omega t}\int_0^t e^{\omega(\lambda)\tau}[\hat{f}(\lambda,\tau+T)-\hat{f}(\lambda,\tau)]d\tau+e^{-\omega t}\int_{-T}^0 e^{\omega(\lambda)\tau}\hat{f}(\lambda,\tau+T)d\tau\,,$$

and, therefore,

$$\int_{-\infty}^\infty e^{i\lambda x-\omega(\lambda)(t+T)}\widetilde{\widehat{f}}(\lambda,\omega(\lambda),t+T)d\lambda-\int_{-\infty}^\infty e^{i\lambda x-\omega(\lambda)t}\widetilde{\widehat{f}}(\lambda,\omega(\lambda),t)d\lambda$$

$$=\int_{-\infty}^\infty e^{i\lambda x}\int_0^t e^{-\omega(\lambda)(t-\tau)}[\hat{f}(\lambda,\tau+T)-\hat{f}(\lambda,\tau)]d\tau d\lambda+\int_{-\infty}^\infty e^{i\lambda x-\omega(\lambda)t}\int_{-T}^0 e^{\omega(\lambda)\tau}\hat{f}(\lambda,\tau+T)d\tau d\lambda\,. \tag{7.23}$$

Working as in the proof of Theorem 9, we show that the last integral in (7.23) tends to zero, as $t\to\infty$.

We claim also that

$$\lim_{t\to\infty}\int_{-\infty}^\infty e^{i\lambda x}\int_0^t e^{-\omega(\lambda)(t-\tau)}[\hat{f}(\lambda,\tau+T)-\hat{f}(\lambda,\tau)]d\tau d\lambda=0\,. \tag{7.24}$$

To prove this, first we write: For $\lambda\neq 0$,

$$\hat{f}(\lambda,\tau+T)-\hat{f}(\lambda,\tau)=\left[\frac{f(0,\tau+T)}{i\lambda}-\frac{f(0,\tau)}{i\lambda}\right]+\left[\frac{f_x(0,\tau+T)}{(i\lambda)^2}-\frac{f_x(0,\tau)}{(i\lambda)^2}\right]+\left[\frac{(f_{xx})\hat{}(\lambda,\tau+T)}{(i\lambda)^2}-\frac{(f_{xx})\hat{}(\lambda,\tau)}{(i\lambda)^2}\right]\,. \tag{7.25}$$

In view of (7.25), the integral in (7.24) is equal to

$$\int_{-1}^1 e^{i\lambda x}\int_0^t e^{-\omega(\lambda)(t-\tau)}[\hat{f}(\lambda,\tau+T)-\hat{f}(\lambda,\tau)]d\tau d\lambda$$

$$+\int_Z e^{i\lambda x}\int_0^t e^{-\omega(\lambda)(t-\tau)}\left[\frac{f(0,\tau+T)}{i\lambda}-\frac{f(0,\tau)}{i\lambda}\right]d\tau d\lambda+\int_Z e^{i\lambda x}\int_0^t e^{-\omega(\lambda)(t-\tau)}\left[\frac{f_x(0,\tau+T)}{(i\lambda)^2}-\frac{f_x(0,\tau)}{(i\lambda)^2}\right]d\tau d\lambda$$

$$\left(\int_{-\infty}^{-1}+\int_1^\infty\right)\left\{e^{i\lambda x}\int_0^t e^{-\omega(\lambda)(t-\tau)}\left[\frac{(f_{xx})\hat{}(\lambda,\tau+T)}{(i\lambda)^2}-\frac{(f_{xx})\hat{}(\lambda,\tau)}{(i\lambda)^2}\right]d\tau d\lambda\right\}\,. \tag{7.26}$$

The first integral in (7.26) is equal to

$$\int_{\lambda=-1}^1\int_{\tau=0}^\infty\int_{y=0}^\infty e^{i\lambda x}e^{-i\lambda y}\chi_{[0,t]}(\tau)e^{-\omega(\lambda)(t-\tau)}[f(y,\tau+T)-f(y,\tau)]dyd\tau d\lambda\,,$$

which tends to zero, as $t\to\infty$, in view of (7.18).

The second and third integrals in (7.26) tend to zero, as $t \to \infty$, in view of (7.16) and (7.17).

Finally, the last integral in (7.26) can be written as follows:

$$\left(\int_{\lambda=-\infty}^{-1} + \int_{\lambda=1}^{\infty}\right) \int_{\tau=0}^{\infty} \int_{y=0}^{\infty} \frac{1}{(i\lambda)^2} e^{i\lambda x} e^{-i\lambda y} \chi_{[0,t]}(\tau) e^{-\omega(\lambda)(t-\tau)} [f_{xx}(y, \tau+T) - f_{xx}(y, \tau)] dy \, d\tau \, d\lambda \, ,$$

and this tends to zero, as $t \to \infty$, by Lebesgue's dominated convergence theorem and assumption (7.19).

Thus, the above computations imply (7.24).

The proof of the theorem can be completed by deriving equation analogous to (7.21) and (7.24) for all the integrals in (5.6).

## Declarations.

*Ethical Approval*: Not applicable.    *Competing Interests*: Not applicable.    *Availability of data and materials*: Not applicable.

*Acknowledgment*. Andreas Chatziafratis gratefully acknowledges partial funding, at different stages of this project, from the: Academy of Athens, State Scholarships Foundation (IKY), Hellenic Foundation for Research and Innovation, and European Research Council (ERC grant 101078061 SINGinGR). A. C. also wishes to express his thankfulness to Professors: N. D. Alikakos, G. Barbatis, T. Bountis, C. M. Dafermos, G. Dassios, G. Fournodavlos, L. Grafakos, M. Grillakis, T. Hatziafratis, A. A. Himonas, S. Kamvissis, N. I. Karachalios, G. Kastis, P. G. Kevrekidis, A. Konstantinidis, D. Mitsotakis, D. T. Papageorgiou, P. Smyrnelis, I. G. Stratis, N. Stylianopoulos, C. E. Synolakis, N. L. Tsitsas and A. Vidras, for providing encouragement and academic support.

## References

[1] J.W. Cahn, J.E. Hilliard, Free energy of a nonuniform system. I. Interfacial free energy, J. Chem. Phys. 28, 258-267, (1958); J.W. Cahn, Free energy of a nonuniform system. II. Thermodynamic basis, J. Chem. Phys. 30, 1121-24, (1959).

[2] H. Cook, Brownian motion in spinodal decomposition, Acta Metallurgica, 18, 297– 306 (1970).

[3] J. S. Langer, Theory of spinodal decomposition in alloys, Ann. Phys. 65, 53–86 (1971).

[4] C. M. Elliott, S. Zheng, On the Cahn-Hillard equation, Arch. Rational Anal. 9, 339-357 (1986).

[5] J. Carr, M. Gurtin, and M. Slemrod, Structured phase transitions on a finite interval, Arch. Rational Mech. Anal. 86 (1984).

[6] V. Alexiades, E. C. Aifantis, On the thermodynamic theory of fluid interfaces, infinite intervals, equilibrium solutions. and minimizers, J. Colloid. Interface Sci. 111, 119-132, (1986).

[7] C. P. Grant, Spinodal decomposition for the Cahn-Hilliard equation, Comm. PDE, 18, 453-490 (1993).

[8] N. Alikakos, P. Bates, G. Fusco, Slow motion for the Cahn-Hilliard equation in one space dimension, J Diff Eq 90 (1991).

[9] P. Bates, J.Xun, Metastable patterns for the Cahn-Hilliard equation: Part I, JDE 111, 421-457 (1995); P. Bates, J. Xun, Metastable patterns for the Cahn- Hilliard equation: Part II, JDE, 117, 165-216 (1995).

[10] N. Alikakos, G. Fusco, The spectrum of the Cahn-Hilliard operator for generic interface in higher space dimensions, Indiana Mathematics Journal 42, 637-674 (1993); N.D. Alikakos, G. Fusco, Slow dynamics for the Cahn-Hilliard equation in higher space dimensions: The motion of bubbles, Arch. Rat. Mech. Anal., 141, 1-61 (1998).

[11] X. Chen, Spectrum for the Allen-Cahn, Cahn-Hilliard, and phase-field equations for generic interfaces, Comm. PDE 19, 1371-1395 (1994).

[12] N.D. Alikakos, G. Fusco, G. Karali, Ostwald ripening in two dimensions- The rigorous derivation of the equations from Mullins-Sekerka dynamics, J Differ Equ, 205, 1–49 (2004); N. D. Alikakos, G. Fusco, G. Karali, The effect of the geometry of the particle distribution in Ostwald Ripening, Comm. Math. Phys. 238, 480-488 (2003); N. D. Alikakos, G. Fusco, G. Karali, Motion of bubbles towards the boundary for the Cahn-Hilliard equation, Eur J Appl Math 15, 103-124 (2004).

[13] D. Antonopoulou, D. Blomker, G. Karali, Front motion in the one-dimensional stochastic Cahn-Hilliard equation, SIAM J. Math. Anal. 44, 3242-3280 (2012).

[14] D. Antonopoulou, G. Karali, A. Millet, Existence and regularity of solution for astochastic Cahn-Hilliard / Allen-Cahn equation with unbounded noise diffusion, J Differ Equ 260, 2383-2417 (2016).

[15] A. Miranville, The Cahn-Hilliard equation: recent advances and applications, CBMS-NSF Regional Conference Series in Applied Mathematics 95, SIAM (2019).

[16] E.C. Aifantis, Internal length gradient (ILG) material mechanics across scales and disciplines, Adv. Appl. Mech. 49 (2016).

[17] E.C. Aifantis, Gradient extension of classical material models: From nuclear & condensed matter scales to earth & cosmological scales, In: Ghavanloo, E., Fazelzadeh, S.A., Marotti de Sciarra, F. (eds), Size-Dependent Continuum Mechanics Approaches. Springer Tracts in Mechanical Engineering. Springer (2021).

[18] Van der Waals, J.D. (1873) Over de Continuiteit van den Gas-enVloeistoftoestand (On the Continuity of the Gas and Liquid State). Ph.D. Thesis, University of Leiden, Leiden; Van der Waals, J.D., Théorie thermodynamique de la capillarité dans l' hypothèse d' une variation continue de densité, Acta. Neerl. Sci. Exact. Nat. 28, 121-209 (1895); Translation: J.S. Rowlinson, The thermodynamic theory of capillarity under the hypothesis of a continuous variation of density, J. Stat. Phys. 20 (1979).

[19] E.C. Aifantis, J.B. Serrin, The mechanical theory of fluid interfaces and Maxwell's rule, J. Colloid Interf. Sci. 96, 517-529 (1983); E.C. Aifantis, J.B. Serrin, Equilibrium solutions in the mechanical theory of fluid microstructures, J. Colloid Interf. Sci. 96, 530-547 (1983).

[20] H.T. Davis, L.E. Scriven, Stress and structure in fluid interfaces, Adv. Chem. Phys. 49, 357-454 (1982).

[21] E.C. Aifantis, On the problem of diffusion in solids, Acta Mech.37, 265–296 (1980).

[22] J. C. Maxwell. On the dynamical theory of gases, Phil. Trans. Roy. Soc. 157, 49–88 (1867).

[23] K. Kuttler, E.C. Aifantis, Existence and uniqueness in nonclassical diffusion, Quart. Appl. Math. 45, 549-560 (1987).

[24] Y. Kuramoto, Diffusion-induced chaos in reaction systems, Progress Theor Phys 64, 346–367 (1978).

[25] G.I. Sivashinsky, On flame propagation under conditions of stoichiometry, SIAM J Appl Math 39, 67–82 (1980).

[26] L.A. Caffarelli and N.E. Müller, An L^infinity bound for solutions of the Cahn-Hilliard equation, Arch. Ration. Mech. Anal. 133, 129-144 (1995).

[27] A.S. Fokas, A unified transform method for solving linear and certain nonlinear PDEs, Proc. Roy. Soc. London Ser. A 453 (1997); A. S. Fokas, On the integrability of linear and nonlinear PDEs, J. Math. Phys. 41, 4188-4237 (2000).

[28] A.S. Fokas, I.M. Gelfand, Integrability of linear and nonlinear evolution equations and the associated nonlinear Fourier transforms, Lett. Math. Phys. 32, 189–210 (1994).

[29] A. S. Fokas, B. Pelloni, Integrable evolution equations in time-dependent domains, Inverse Problems 17, 919–935 (2001).

[30] A.S. Fokas, Integrable nonlinear evolution equations on the half-line, Commun. Math. Phys. 230, 1–39 (2002).

[31] A. S. Fokas, A new transform method for evolution partial differential equations, IMA J. Appl. Math. 67, 559-590 (2002).

[32] A. S. Fokas, M. Zyskin, The fundamental differential form and boundary-value problems, Q J Mech Appl Math 55 (2002).

[33] A.S. Fokas, P.F. Schultz, Long-time asymptotics of moving boundary problems using an Ehrenpreis-type representation and its Riemann-Hilbert nonlinearization, Commun Pure Appl Math. 56, 517-548 (2003).

[34] G. Dassios, From d' Alembert and Fourier to Gelfand and Fokas via Lax: Linearity revisited, In: Advances in Scattering and Biomedical Engineering, pp. 252-259, World Scientific (2004).

[35] A.S. Fokas, S. Kamvissis, Zero-dispersion limit for integrable equations on the half-line with linearisable data, Abstr. Appl. Anal. 2004:5, 361–370 (2004).

[36] A.S. Fokas, A.R. Its, L.-Y. Sung, The nonlinear Schrödinger equation on the half-line, Nonlinearity 18, 1771 (2005).

[37] A.S. Fokas, L.Y. Sung, Generalized Fourier transforms, their nonlinearization and the imaging of the brain, Notices Amer. Math Soc. 52, 1178-92 (2005).

[38] A.S. Fokas, A. Iserles, V. Marinakis, Reconstruction algorithm for single photon emission computed tomography and its numerical implementation, J. R. Soc. Interface. 3, 45–54 (2006).

[39] A.S. Fokas, From Green to Lax via Fourier, In: Recent Advances in Nonlinear Partial Differential Equations and Applications, Proc. Sympos. Appl. Math. 65 (Dedicated to P.D. Lax and L. Nirenberg), Amer. Math. Soc. (2007).

[40] A.S. Fokas, A Unified Approach to Boundary Value Problems, CBMS-NSF Series Appl Math 78, SIAM (2008); AS Fokas, B Pelloni (eds), Unified Transform for Boundary Value Problems: Applications & Advances, SIAM, Philadelphia (2015).

[41] G.M. Dujardin, Asymptotics of linear initial boundary value problems with periodic boundary data on the half-line and finite intervals, Proc R Soc London 465, 3341–60 (2009).

[42] A. S. Fokas, Lax Pairs: A novel type of separability (invited paper), Inverse Problems 25, 1-44 (2009).

[43] A.S. Fokas, N. Flyer, S.A. Smitheman, E.A. Spence, A semi-analytical numerical method for solving evolution and elliptic partial differential equations, J. Comp. Appl. Math. 227, 59-74 (2009).

[44] A. S. Fokas, E. A. Spence, Novel analytical and numerical methods for elliptic boundary value problems, In: Engquist B., Fokas A., Hairer E., Iserles A. (Eds.), Highly Oscillatory Problems, London Mathematical Society Lecture Note Series 366. Cambridge: Cambridge University Press (2009).

[45] A. S. Fokas, E. A. Spence, Synthesis, as opposed to separation, of variables, SIAM Review 54 (2012).

[46] B. Deconinck, T. Trogdon, V. Vasan, The method of Fokas for solving linear partial differential equations, SIAM Review 56, 159-186 (2014).

[47] A.S. Fokas and G. A. Kastis, Boundary value problems and medical imaging, Journal of Physics: Conf. Ser. 490 (2014).

[48] B. Pelloni, Advances in the study of boundary value problems for nonlinear integrable PDEs, *Nonlinearity* 28, 1-38 (2015).

[49] D. C. Antonopoulou, S. Kamvissis, On the Dirichlet-to-Neumann problem for the 1-dimensional cubic NLS equation on the half-line, Nonlinearity 28, 3073-3099 (2015); Addendum, Nonlinearity 29, 3206-3214 (2016).

[50] B. Deconinck, N.E. Sheils, D.A. Smith, The linear KdV equation with an interface, Comm Math. Phys. 347 (2016).

[51] A. S. Fokas, A. Alexandrou Himonas, D. Mantzavinos, The nonlinear Schrödinger equation on the half-line, Trans Amer Math Soc 369, 681-709 (2017).

[52] S.-F. Tian, Initial-boundary value problems for the general coupled nonlinear Schrödinger equation on the interval via the Fokas method, J Differ Equ 262, 506-558 (2017).

[53] A. Fernandez, D. Baleanu, A.S. Fokas, Solving PDEs of fractional order using the unified transform method. Appl. Math. Comput. 339, 738–749 (2018).

[54] E. Kesici, B. Pelloni, T. Pryer, D. A. Smith, A numerical implementation of the unified Fokas transform for evolution problems on a finite interval, Eur. J. Appl. Math. 29(3), 543-567 (2018).

[55] P. D. Miller, D. A. Smith, The diffusion equation with nonlocal data, J. Math. Anal. Appl. 466, 1119-1143 (2018).

[56] F.P.J. de Barros, M.J. Colbrook and A.S. Fokas, A hybrid analytical-numerical method for solving advection-dispersion problems on a half-line, Int. J. Heat Mass Tran. 139, 482-491 (2019).

[57] A.S. Fokas, T. Özsari, New rigorous developments regarding the Fokas method and an open problem, In: Solved and Unsolved Problems, M.T. Rassias, ed., EMS Newsletter 113, 60-61 (2019).

[58] A.A. Himonas, D. Mantzavinos, F. Yan, The Korteweg-de Vries equation on an interval. J. Math. Phys. 60, 1–26 (2019).

[59] T. Trogdon, G. Biondini, Evolution partial differential equations with discontinuous data, Quart. Appl. Math. 77 (2019).

[60] A. Batal, A.S. Fokas, T. Özsari, Fokas method for linear boundary value problems involving mixed spatial derivatives, Proc. R. Soc. A 476, 20200076 (2020).

[61] A. Alexandrou Himonas, D. Mantzavinos, Well-posedness of the nonlinear Schrödinger equation on the half-plane, Nonlinearity 33, 5567–5609 (2020).

[62] G Hwang, Initial-boundary value problems for the one-dimensional linear advection–dispersion equation with decay, Z Naturforschung A 75, 713-725 (2020).

[63] K. Kalimeris, T. Özsarı, An elementary proof of the lack of null controllability for the heat equation on the half line, Appl. Math. Lett. 104 (2020).

[64] Olver PJ, Sheils NE, Smith DA. Revivals and fractalisation in the linear free space Schrödinger equation. Q Appl Math. 78, 161-192 (2020).

[65] A.S. Fokas, M.C. van der Weele, The unified transform for evolution equations on the half-line with time-periodic boundary conditions, Stud Appl Math 147 (2021).

[66] B. Deconinck, A.S. Fokas, J. Lenells, The implementation of the unified transform to the nonlinear Schrödinger equation with periodic initial conditions, Lett Math Phys 111:17 (2021).

[67] A. Himonas, C. Madrid, F. Yan, The Neumann and Robin problems for the KdV equation on the half-line, J Math Phys 62, 111503 (2021).

[68] A. A. Himonas, F. Yan, The modified KdV system on the half-line, J Dyn Diff Equ (2023).

[69] A.S. Fokas, K. Kalimeris, Extensions of the d'Alembert formulae to the half line and the finite interval obtained via the unified transform, IMA J. Appl. Math. 87, 1010-42 (2022).

[70] A.S. Fokas, E Kaxiras, Modern Mathematical Methods for Computational Sciences & Engineering, World Scientific (2022).

[71] J.M. Lee, J. Lenells, The nonlinear Schrödinger equation on the half-line with homogeneous Robin boundary conditions, Proc. Lond. Math. Soc. 126, 334-89 (2023).

[72] A. Chatziafratis, Rigorous analysis of the Fokas method for linear evolution PDEs on the half-space, Thesis, Advisors: N. Alikakos, G. Barbatis, I.G. Stratis, National and Kapodistrian University of Athens (2019).

[73] A. Chatziafratis, S. Kamvissis, I.G. Stratis, Boundary behavior of the solution to the linear KdV equation on the half-line, Stud. Appl. Math. 150, 339-379 (2023).

[74] A. Chatziafratis, L. Grafakos, S. Kamvissis, Long-range instability of linear evolution PDEs on semi-bounded domains via the Fokas method, Dyn. PDE 21, 97-169 (2024); A. Chatziafratis, L. Grafakos, S. Kamvissis, I.G. Stratis, Instabilities of linear evolution PDE on the half-line via the Fokas method, In: Chaos, Fractals and Complexity, Eds.: T. Bountis et al, Springer Nature (2023).

[75] A. Chatziafratis, T. Ozawa, S.-F. Tian, Rigorous analysis of the unified transform method and long-range instabilities for the inhomogeneous time-dependent Schrödinger equation on the quarter-plane, Math. Ann. (2023)

[76] A. Chatziafratis, E.C. Aifantis, A. Carbery, A.S. Fokas, Integral representations for the double-diffusion system on the half-line, Z. Angew. Math. Phys. 75 (2024).

[77] A. Chatziafratis, A.S. Fokas, E.C. Aifantis, On Barenblatt's pseudoparabolic equation with forcing on the half-line via the Fokas method, Z. Angew. Math. Mech. 104 (2024).

[78] J.L. Bona, A. Chatziafratis, H. Chen, S. Kamvissis, The linear BBM-equation on the half-line revisited, Lett. Math. Phys. 114 (2023).

[79] A. Chatziafratis, T. Ozawa, New phenomena for Sobolev-type evolution equations: Asymptotic analysis for a celebrated pseudo-parabolic model on the quarter-plane, Partial Differ. Equ. 5 (2024).

[80] A. Chatziafratis, S. Kamvissis, Uniqueness theory and ill-posedness for linear evolution PDEs on the quarter-plane, preprint (2023).

[81] A. Chatziafratis et al., On controllability for linear evolution PDEs posed in the quarter-plane via the Fokas method, preprint (2024).

[82] A. Chatziafratis, E.C. Aifantis, Explicit Ehrenpreis-Palamodov-Fokas representations for the Sobolev-Barenblatt model and the Rubinshtein-Aifantis system on semi-strips and quarter-planes, Lobachevskii J Math (to appear).

[83] A. Chatziafratis et al., Analysis for the linearized Whitham-Broer-Kaup system on the quarter-plane, preprint (2024).

[84] A. Chatziafratis et al., On the Maxwell-Cattaneo-Vernotte and Moore-Gibson-Thompson equations on the half-line, preprint (2024).

[85] J.L. Bona, J. Wu, Temporal growth and eventual periodicity for dispersive wave equations in a quarter plane, Discrete Contin. Dyn. Syst 23, 1141–1168 (2009).

[86] J.L. Bona, J. Lenells, The KdV equation on the half-line: time-periodicity and mass transport, SIAM J Math Anal 52 (2020).

[87] J. Shen, J. Wu, J.M. Yuan, Eventual periodicity for the KdV equation on a half-line, Phys D: Nonlinear Phenom 227, 105-119 (2007).

[88] J.M.K. Hong, J. Wu, J.M. Yuan, A new solution representation for the BBM equation in a quarter plane and the eventual periodicity, Nonlinearity 22 (2009)

[89] A. Chatziafratis, T. Hatziafratis, Density in the half-line Schwartz space of functions whose Fourier-Laplace transform has natural boundary the real line, preprint (2024). arXiv:2403.03718 [math.CV]; The natural boundary of Fourier-Laplace transforms of functions in the half-line Schwartz space, Integr. Transf. Spec. Funct. (to appear).

[90] P. A. Clarkson, A. S. Fokas, M. J. Ablowitz, Hodograph transformations of linearizable partial differential equations, SIAM J. Appl. Math. 49, 1188-1209 (1989).

[91] A.S. Fokas, A.R. Its, The linearization of the initial-boundary value problem of the nonlinear Schrödinger equation, SIAM J. Math. Anal. 27, 738–64 (1996).

[92] J. Lenells, A.S. Fokas, Linearizable boundary value problems for the elliptic sine-Gordon and the elliptic Ernst equations, J. Nonlinear Math. Phys. 27, 337-356 (2020).